# Formal and Analytic Diagonalization of Operator functions

M. Stiefenhofer

ABSTRACT. We give conditions for local diagonalization of analytic operator families $L(\varepsilon)$, acting between real or complex Banach spaces, of the form $L(\varepsilon) = \psi(\varepsilon) \cdot \Delta(\varepsilon) \cdot \phi^{-1}(\varepsilon)$ with a diagonal operator polynomial $\Delta(\varepsilon)$ and analytic near identity bijections $\psi(\varepsilon)$ and $\phi(\varepsilon)$.

The transformation $\phi(\varepsilon)$ is constructed from an operator Töplitz matrix obtained from Jordan chains of increasing length. The basic assumption is given by stabilization of the Jordan chains at length $k$ in the sense that no root elements with finite rank above $k$ are allowed to exist. Jordan chains with infinite rank may appear. These assumptions ensure finite pole order equal to $k$ of the generalized inverse of $L(\varepsilon)$ at $\varepsilon = 0$. The Smith form arises immediately.

Smooth continuation of kernels and ranges towards appropriate limit spaces at $\varepsilon = 0$ is considered using associated families of analytic projection functions.

No Fredholm properties or other finiteness assumptions, besides the pole order, are assumed. Real and complex Banach spaces are treated without difference by elementary analysis of the system of undetermined coefficients.

Formal power series solutions of the system of undetermined coefficients are constructed, which are turning into convergent solutions, as soon as analyticity of $L(\varepsilon)$ and continuity of the projections is assumed. Along these lines, results concerning linear Artin approximation follow immediately, which are well known in finite dimensions. The main technical tool is given by a defining equation of Nakayama Lemma type.

**Keywords:** Diagonalization, Jordan chain, Generalized Inverse, Transversalization, Diagonal Operator Function, Toeplitz Matrix, Smith Form, Artin approximation

*Contents*

## *1. Introduction*

Given an analytic family of linear operators $L(\varepsilon) = \sum_{i=0}^{\infty} \varepsilon^i \cdot L_i$, $L \in C^\omega(U, L[B, \bar{B}])$, $U$ open with $B, \bar{B}$ real or complex Banach spaces, $0 \in U \subset \mathbb{K} = \mathbb{R}, \mathbb{C}$ and $L[B, \bar{B}]$ denoting bounded linear operators from $B$ to $\bar{B}$.

We give conditions for diagonalization of $L(\varepsilon)$ in the sense that families of analytic near identity transformations $\phi(\varepsilon)$ and $\psi(\varepsilon)$ of $B$ and $\bar{B}$ exist satisfying

$$\psi^{-1}(\varepsilon) \cdot L(\varepsilon) \cdot \phi(\varepsilon) = \Delta(\varepsilon) \tag{1.1}$$

with a diagonal operator polynomial of degree $k \geq 0$ of the form

$$\Delta(\varepsilon) = S_1 P_1 + \varepsilon \cdot S_2 P_2 + \cdots + \varepsilon^k \cdot S_{k+1} P_{k+1} \,. \tag{1.2}$$

In the sense of [Gohberg], the families $L(\varepsilon)$ and $\Delta(\varepsilon)$ are analytically equivalent and the following diagram commutes.

$$\begin{array}{ccc} B & \overset{\phi(\varepsilon)}{\leftarrow} & B \\ L(\varepsilon)\downarrow & & \downarrow \Delta(\varepsilon) \\ \bar{B} & \underset{\psi(\varepsilon)}{\leftarrow} & \bar{B} \end{array} \tag{1.3}$$

The operators in (1.2) are defined using direct sum decompositions of $B$ and $\bar{B}$ according to

$$\begin{array}{ccccccccccc} B & = & \overbrace{N_1^c}^{\downarrow P_1} & \oplus & \overbrace{N_2^c}^{\downarrow P_2} & \oplus & \cdots & \oplus & \overbrace{N_{k+1}^c}^{\downarrow P_{k+1}} & \oplus & N_{k+1} \\ & & \uparrow & & \uparrow & & & & \uparrow & & \\ & & \boxed{S_1} & & \boxed{S_2} & & \cdots & & \boxed{S_{k+1}} & & \\ & & \downarrow & & \downarrow & & & & \downarrow & & \\ \bar{B} & = & R_1 & \oplus & R_2 & \oplus & \cdots & \oplus & R_{k+1} & \oplus & R_{k+1}^c \end{array} \tag{1.4}$$

with $P_1, \ldots, P_{k+1}$ meaning bounded projections to $N_1^c, \ldots, N_{k+1}^c$ respectively. When restricted to subspaces $N_i^c$ and $R_i$, the operators $S_i \in L[B, \bar{B}]$ represent isomorphism $S_i \in GL[N_i^c, R_i]\,, i = 1, \ldots, k+1$. In addition, $k \geq 0$ is given by the maximal length of finite Jordan chains of $L(\varepsilon)$, i.e. in the sense of [Kaballo], we assume stabilization of the Jordan chains at $k \geq 0$. The subspaces in (1.4) are assumed to be closed.

The matrix representation of $\Delta(\varepsilon)$ with respect to the direct sums in (1.4) reads

$$\begin{array}{cc} & \begin{array}{cccc} \boxed{N_1^c} & \cdots & \boxed{N_{k+1}^c} & \boxed{N_{k+1}} \end{array} \\ \Delta(\varepsilon) = & \begin{pmatrix} \boxed{S_1} & & & \\ & \ddots & & \\ & & \boxed{\varepsilon^k \cdot S_{k+1}} & \\ & & & \boxed{0} \end{pmatrix} \begin{array}{c} \boxed{R_1} \\ \vdots \\ \boxed{R_{k+1}} \\ \boxed{R_{k+1}^c} \end{array} \end{array} \tag{1.5}$$

and diagonalization of $L(\varepsilon)$ occurs by analytic left and right transformation according to (1.1).

No Fredholm properties are needed. Finiteness merely occurs with respect to maximal finite length $k \geq 0$ of Jordan chains, implying finite pole order of corresponding generalized inverses. This aspect is also stressed in [Magnus] and [Kaballo].

If $L(\varepsilon)$ is given by a matrix function $L(\varepsilon) \in L[\mathbb{K}^n, \mathbb{K}^m]$, then $\Delta(\varepsilon)$ represents the Smith form of $L(\varepsilon)$.

The coefficients $\phi_i$ and $\psi_i$ of the transformations $\phi(\varepsilon) \coloneqq I_B + \sum_{i=1}^{\infty} \varepsilon^i \cdot \phi_i$ and $\psi(\varepsilon) \coloneqq I_{\bar{B}} + \sum_{i=1}^{\infty} \varepsilon^i \cdot \psi_i$ are calculated recursively in dependance of given coefficients $L_i$ of the family $L(\varepsilon) = \sum_{i=0}^{\infty} \varepsilon^i \cdot L_i$. Using these transformations, we define analytic families of subspaces

$$N(\varepsilon) \coloneqq \phi(\varepsilon) \cdot N_{k+1} \qquad and \qquad R(\varepsilon) \coloneqq \psi(\varepsilon) \cdot [\, R_1 \oplus \cdots \oplus R_{k+1} \,] \tag{1.6}$$

for $\varepsilon \in U$, motivated by the fact that kernels and ranges of the diagonal operator function $\Delta(\varepsilon)$ in (1.2) are given by

$$N[\, \Delta(\varepsilon) \,] \equiv N_{k+1} \qquad and \qquad R[\, \Delta(\varepsilon) \,] \equiv R_1 \oplus \cdots \oplus R_{k+1} \tag{1.7}$$

for $\varepsilon \neq 0$. Then, concerning kernels and ranges of $L(\varepsilon)$ we obtain from (1.1)

$$N(\varepsilon) = N[\, L(\varepsilon) \,] \qquad and \qquad R(\varepsilon) = R[\, L(\varepsilon) \,] \tag{1.8}$$

within the punctured neighbourhood $\varepsilon \in U \setminus \{0\}$ and we see that kernels and ranges of $L(\varepsilon)$ can analytically be continued to $\varepsilon = 0$ by $N(0) = N_{k+1} \subset N[L(0)]$ and $R(0) = R_1 \oplus \cdots \oplus R_{k+1} \supset R[L(0)]$, i.e. smoothing of kernels $N[L(\varepsilon)]$ and ranges $R[L(\varepsilon)]$ occurs by $N(\varepsilon)$ and $R(\varepsilon)$, as defined in (1.6).

In the following diagram, the constellation is qualitatively depicted in finite dimensions in case of $\mathbb{K} = \mathbb{R}$.

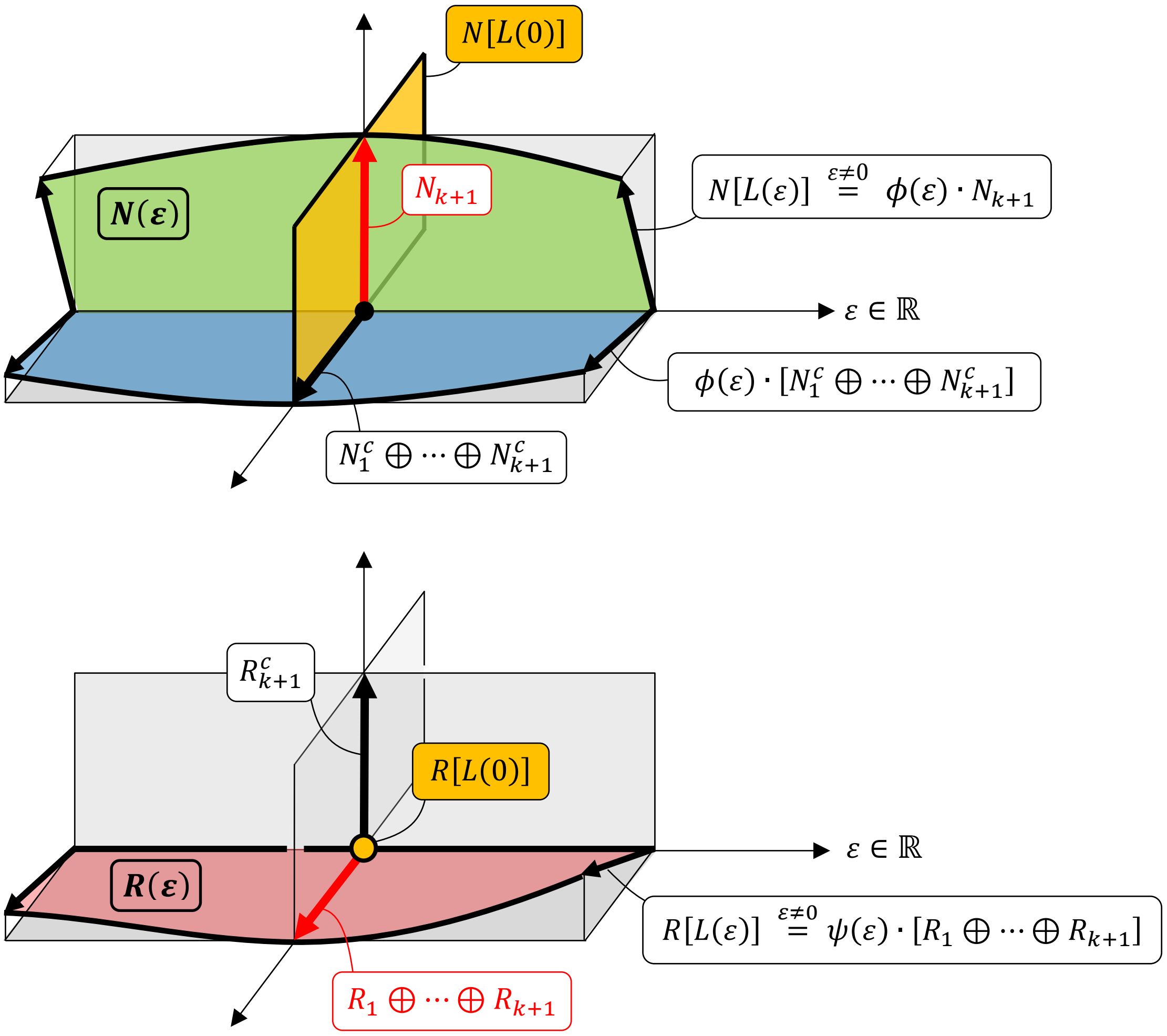

Figure 1 : Breakdown of $N[L(0)]$ and blow up of $R[L(0)]$ when passing to $\varepsilon \neq 0$.

In figure 1 upper part, the green surface of kernels $N[L(\varepsilon)] \subset B$ is smoothly continued to $\varepsilon = 0$, in this way approaching the limit space $N_{k+1}$. If a singularity occurs at $\varepsilon = 0$, then the red marked subspace $N_{k+1}$ typically is a proper subspace of $N[L(0)]$ (orange).

In the lower part of figure 1, the red surface of ranges $R[L(\varepsilon)] \subset \bar{B}$ is continued analytically to $\varepsilon = 0$ with corresponding limit space $R_1 \oplus \cdots \oplus R_{k+1}$, allowed to contain the subspace $R[L(0)]$.

Reversely, during the passage from $\varepsilon = 0$ to $\varepsilon \neq 0$ the kernel $N[L(0)]$ typically collapses to $N[L(\varepsilon)]$, whereas the range $R[L(0)]$ is blown up to $R[L(\varepsilon)]$.

Now, the embedding of $N[L(\varepsilon)]$ and $R[L(\varepsilon)]$ into analytic families $N(\varepsilon)$ and $R(\varepsilon)$ offers the possibility to define smooth generalized inverses $L^{-1}(\varepsilon)$ for $\varepsilon \neq 0$ with a pole of order $k \geq 0$ at $\varepsilon = 0$. First, a smooth generalized inverse of the diagonal operator polynomial $\Delta(\varepsilon)$ from (1.2) is given by

$$\Delta^{-1}(\varepsilon) = \varepsilon^{-k} \cdot S_{k+1}^{-1} \mathcal{P}_{k+1} + \cdots + S_1^{-1} \mathcal{P}_1 \tag{1.9}$$

with $\mathcal{P}_1, \dots, \mathcal{P}_{k+1}$ denoting bounded projections to $R_1, \dots, R_{k+1}$ respectively. Then, by (1.1)

$$\begin{aligned} L^{-1}(\varepsilon) &= \phi(\varepsilon) \cdot \Delta^{-1}(\varepsilon) \cdot \psi^{-1}(\varepsilon) \qquad (1.10) \\ &= \left[\, I_B + \textstyle\sum_{i=1}^{\infty} \varepsilon^i \cdot \phi_i \,\right] \cdot \left[\, \varepsilon^{-k} \cdot S_{k+1}^{-1} \mathcal{P}_{k+1} + \cdots + S_1^{-1} \mathcal{P}_1 \,\right] \cdot \left[\, I_{\bar{B}} + \textstyle\sum_{i=1}^{\infty} \varepsilon^i \cdot \psi_i \,\right]^{-1} \\ &= \varepsilon^{-k} \cdot S_{k+1}^{-1} \mathcal{P}_{k+1} + \varepsilon^{-k+1} \cdot r(\varepsilon) \end{aligned}$$

with an analytic remainder function $r(\varepsilon)$, implying a generalized inverse $L^{-1}(\varepsilon)$, which is analytic in a punctured neighbourhood $\varepsilon \in U \setminus \{0\}$ with pole of order $k$ at $\varepsilon = 0$. The coefficients of $L^{-1}(\varepsilon)$ can be calculated in dependance of given coefficients $L_i$.

Finally, for $\varepsilon \in U$ analytic families of projections to the subspaces $\phi(\varepsilon) \cdot [N_1^c \oplus \cdots \oplus N_{k+1}^c] \subset B$ (blue surface in figure 1) and $R(\varepsilon) = \psi(\varepsilon) \cdot [R_1 \oplus \cdots \oplus R_{k+1}] \subset \bar{B}$ (red surface) are given by

$$\begin{aligned} L^{-1}(\varepsilon) \cdot L(\varepsilon) &= \phi(\varepsilon) \cdot (P_1 + \cdots + P_{k+1}) \cdot \phi^{-1}(\varepsilon) \qquad (1.11) \\ L(\varepsilon) \cdot L^{-1}(\varepsilon) &= \psi(\varepsilon) \cdot (\mathcal{P}_1 + \cdots + \mathcal{P}_{k+1}) \cdot \psi^{-1}(\varepsilon) \end{aligned}$$

respectively. Geometrically speaking, the analytic family of generalized inverses $L^{-1}(\varepsilon),\ \varepsilon \neq 0$ is constructed with respect to blue and red isomorphic subspaces $\phi(\varepsilon) \cdot [N_1^c \oplus \cdots \oplus N_{k+1}^c]$ and $\psi(\varepsilon) \cdot [R_1 \oplus \cdots \oplus R_{k+1}]$ in figure 1.

Diagonalization of analytic matrix functions can be found in [Kaashoek]. Concerning diagonalization of analytic operator functions with $L_0$ a Fredholm operator, see [López-Gómez], [Gohberg] and [Mennicken]. In [Gohberg], the family $L(\varepsilon)$ is allowed to be meromorphic.

In [Bart] and [Kaballo], the restriction to finite dimensions is removed and the existence of analytic families $N(\varepsilon)$ and $R(\varepsilon)$ is shown with corresponding family of smooth meromorphic generalized inverses and pole of order $k \geq 0$ at $\varepsilon = 0$. On the other hand, diagonalization of $L(\varepsilon)$ in the sense of (1.1) is not performed in [Bart] and [Kaballo].

In the paper at hand, we perform diagonalization in infinite dimensions. As in [Bart] and [Kaballo], we presuppose stabilization of the Jordan chains at $k \geq 0$, as well as closedness of corresponding subspaces. Stabilization of Jordan chains prevents an essential singularity to appear at $\varepsilon = 0$, i.e. we remain on the level of poles with respect to generalized inverses. Closedness of subspaces is needed for showing convergence of power series during the construction of the transformations. Concerning essential singularities, we refer to [Albrecht].

We did not succeed to ensure the local existence of $\phi(\varepsilon)$ and $\psi(\varepsilon)$ by implicit function theorem or an appropriate contraction mapping principle, but instead we had to set up formal power series $\phi(\varepsilon) := I_B + \sum_{i=1}^{\infty} \varepsilon^i \cdot \phi_i$ and $\psi(\varepsilon) := I_{\bar{B}} + \sum_{i=1}^{\infty} \varepsilon^i \cdot \psi_i$, satisfying diagonalization (1.1) first as a power series relation for arbitrary power series $L(\varepsilon) = \sum_{i=0}^{\infty} \varepsilon^i \cdot L_i$.

In the sense of [Töplitz], $\phi(\varepsilon)$ represents the power series associated to an operator Töplitz matrix with an infinity of rows and columns and, as soon as convergence of $L(\varepsilon)$ is assumed, convergence of $\phi(\varepsilon)$ can be shown by use of an analytic defining equation with solution given by $\phi(\varepsilon)$. Then, the transformation $\phi(\varepsilon)$ of $B$ allows to define the analytic family $S(\varepsilon) = L(\varepsilon) \cdot \phi(\varepsilon)$, characterized by a triangular structure when represented with respect to the subspaces in (1.4).

Finally, it is a straightforward calculation to split of the diagonal operator polynomial $\Delta(\varepsilon)$ from the triangular family $S(\varepsilon)$ to end up with $S(\varepsilon) = \psi(\varepsilon) \cdot \Delta(\varepsilon)$ and $\psi(\varepsilon)$ to be an analytic family of transformations of $\bar{B}$. Geometrically, the reverse transformation $\psi^{-1}(\varepsilon)$ merely is turning the shear structure of $S(\varepsilon)$ into a diagonal structure of $\Delta(\varepsilon) = \psi^{-1}(\varepsilon) \cdot S(\varepsilon) = \psi^{-1}(\varepsilon) \cdot L(\varepsilon) \cdot \phi(\varepsilon)$.

Real and complex Banach spaces are treated without difference, based on elementary analysis of formal power series, designed to solve the system of undetermined coefficients

$$\sum_{i+j=k} L_i \cdot b_j = 0\,, \qquad k = 0, \cdots, \infty \tag{1.12}$$

corresponding to the power series equation

$$L(\varepsilon) \cdot b(\varepsilon) = \left( \sum_{i=0}^{\infty} \varepsilon^i \cdot L_i \right) \cdot \left( \sum_{j=0}^{\infty} \varepsilon^j \cdot b_j \right) = 0\,. \tag{1.13}$$

In other words, high order approximations, i.e. Jordan chains of (1.13) are analyzed and used to define $\phi(\varepsilon)$. The construction of $\phi(\varepsilon) = I_B + \sum_{i=1}^{\infty} \varepsilon^i \cdot \phi_i$ closely follows [S1] where high order approximations are derived in a nonlinear context.

The construction may also be interpreted as an extended version of tranversalization introduced in [Esquinas] and optimized in [López-Gómez] with respect to Fredholm operators of Index 0. Therein, the family $L(\varepsilon)$ is ensured to be a bounded bijection for $\varepsilon \neq 0$ by assuming an algebraic eigenvalue of order $k$ to occur at $\varepsilon = 0$. Then, kernels and ranges simplify to $N[L(\varepsilon)] \equiv \{0\}$ and $R[L(\varepsilon)] \equiv \bar{B}$ for $\varepsilon \neq 0$. In [López-Gómez] real and complex Banach spaces are also treated in a parallel way.

Section 2 aims to give some motivation concerning the correlation between stabilization of Jordan chains in $B$ and stabilization of leading coefficients in $\bar{B}$. Additionally, a more detailled outline of the proof to derive the transformation $\phi(\varepsilon)$ is given. Here, we also stress the significance of transformation steps from $k$ to $2k$. Section 2 may be skipped without losing the main aspects of the proofs.

In section 3, the recursion concerning the direct sums in (1.4) and the coefficients of the transformation $\phi(\varepsilon) := I_B + \sum_{i=1}^{\infty} \varepsilon^i \cdot \phi_i$ is stated. The results concerning Jordan chains and leading coefficients are summarized in Lemma 1.

In section 4, triangularization of $k + 1$ leading coefficients of $L(\varepsilon)$ is performed using polynomial pre-transformations $p_i(\varepsilon)$ of degree $i$ in $\varepsilon$ with $i$ varying from 0 to $k$. The results concerning this kind of partial triangularization are summarized in Theorem 1 and may be interpreted as transversalization without assuming an isolated singularity to occur at $\varepsilon = 0$.

Section 5 completes the procedure of triangularization of $L(\varepsilon)$ up to infinity, thereby using the formal power series $\phi(\varepsilon)$ from section 3 with associated Töplitz structure. The next step from formal triangularization to formal diagonalization of $L(\varepsilon)$ is simple. The results are stated in Theorem 2.

Finally, in sections 6 and 7 it is shown that the power series $\phi(\varepsilon)$ represents a convergent operator function, if analyticity of $L(\varepsilon)$ and closedness of subspaces is supposed, hence turning formal into analytic results, as summarized in Theorem 3.

We close the paper with some relations to commutative algebra in the context of linear Artin approximation.

**Remarks: 1)** From (1.2) we see that $\Delta(\varepsilon)$ may further be factorized according to

$$\begin{aligned}\Delta(\varepsilon) &= S_1P_1 + \varepsilon \cdot S_2P_2 + \cdots + \varepsilon^k \cdot S_{k+1}P_{k+1} \\ &= (S_1P_1 + \cdots + S_{k+1}P_{k+1}) \cdot \left(P_1 + \cdots + \varepsilon^k \cdot P_{k+1}\right) \qquad (1.14) \\ &=: S_P \cdot P(\varepsilon)\end{aligned}$$

yielding from $\psi^{-1}(\varepsilon) \cdot L(\varepsilon) \cdot \phi(\varepsilon) = \Delta(\varepsilon) = S_P \cdot P(\varepsilon)$ factorization to the constant operator $S_P$ by

$$\psi^{-1}(\varepsilon) \cdot L(\varepsilon) \cdot \phi(\varepsilon) \cdot P^{-1}(\varepsilon) = S_P \qquad (1.15)$$

for $\varepsilon \neq 0$. Here $P^{-1}(\varepsilon) \in GL[N_1^c \oplus \cdots \oplus N_{k+1}^c]$, $\varepsilon \neq 0$ with

$$P^{-1}(\varepsilon) = \varepsilon^{-k} \cdot P_{k+1} + \cdots + P_1 \,. \qquad (1.16)$$

Linearization by (1.15) to a constant operator $S_P$ is used in [S2] and [S3], where focus is put on deriving a regular system at $\varepsilon = 0$ by an appropriate blow up scaling procedure.

**2)** Alternatively, we obtain from $\psi^{-1}(\varepsilon) \cdot L(\varepsilon) \cdot \phi(\varepsilon) = S_P \cdot P(\varepsilon)$ the factorization

$$\begin{aligned}L(\varepsilon) &= \psi(\varepsilon) \cdot S_P \cdot P(\varepsilon) \cdot \phi^{-1}(\varepsilon) \\ &= \psi(\varepsilon) \cdot (S_1P_1 + \cdots + S_{k+1}P_{k+1}) \cdot \left(P_1 + \cdots + \varepsilon^k \cdot P_{k+1}\right) \cdot \phi^{-1}(\varepsilon) \\ &=: Q(\varepsilon) \cdot \left(P_1 + \cdots + \varepsilon^k \cdot P_{k+1}\right) \cdot \phi^{-1}(\varepsilon) \qquad (1.17)\end{aligned}$$

with $Q(\varepsilon)$ an analytic family from $B$ to $\bar{B}$.

Now, if $R_{k+1}^c = \{0\}$ and $N_{k+1} = \{0\}$ are assumed in (1.4), then $S_P = S_1P_1 + \cdots + S_{k+1}P_{k+1}$ as well as $Q(\varepsilon)$ are turning into bijections from $B$ to $\bar{B}$ with $P_1 + \cdots + P_{k+1} = I_B$ and we obtain the Smith factorization of $L(\varepsilon)$ by (1.17), as formulated in [López-Gómez]. Note also that in this situation the analytic families of subspaces from (1.6) simplify to $N(\varepsilon) \equiv \{0\}$ and $R(\varepsilon) \equiv \bar{B}$ and the generalized inverse $L^{-1}(\varepsilon)$ from (1.10) turns into the resolvent with pole of order $k$ at $\varepsilon = 0$.

A sufficient condition for $R_{k+1}^c = \{0\}$ and $N_{k+1} = \{0\}$ is given by $L(0) = L_0$ to be a Fredholm operator of index 0. Under these assumptions, the Smith form (1.17) and the resolvent is also derived in [Mennicken].

Then it is also possible to represent the resolvent $L^{-1}(\varepsilon)$ in an optimal way by use of the Keldysh theorem and the Jordan chains can effectively be calculated by contour integrals, even in a global context, as demonstrated in [Beyn] and [Latushkin].

Our main tool for analyzing the passage from $\varepsilon = 0$ to $\varepsilon \neq 0$ is given by curves $b(\varepsilon) = \sum_{i=0}^{\infty} \varepsilon^i \cdot b_i$ in $B$ that are mapped to $\bar{B}$ by $L(\varepsilon)$ according to

$$L(\varepsilon) \cdot b(\varepsilon) = \underbrace{L_0 b_0}_{0-th\ order} + \varepsilon \cdot \underbrace{(L_0\, L_1)\begin{pmatrix} b_1 \\ b_0 \end{pmatrix}}_{1-th\ order} + \varepsilon^2 \cdot \underbrace{(L_0\, L_1\, L_2)\begin{pmatrix} b_2 \\ b_1 \\ b_0 \end{pmatrix}}_{2-th\ order} + \cdots$$

$$= \sum_{l=0}^{\infty} \varepsilon^l \cdot \underbrace{\sum_{i+j=l} L_i b_j}_{k-th\ order} . \tag{2.1}$$

If $L_0$ is surjective, then the 0-th order coefficient $L_0 b_0$ takes every element in $\bar{B}$ implying $R[L(0)] = R[L_0] = \bar{B}$, as well as $R[L(\varepsilon)] = \bar{B}$ by open mapping theorem. Hence, in this regular constellation, $R[L(0)]$ is continued to $R[L(\varepsilon)] = \bar{B}$ in a trivial way.

If $L_0$ is not surjective by $R[L(0)] \subsetneq \bar{B}$, then we may add to $\bar{R}_0 \coloneqq R[L(0)] = \{L_0 b_0 \mid b_0 \in B\}$ leading coefficients of 1-th order $\bar{R}_1 \coloneqq \{L_0 b_1 + L_1 b_0 \mid b_1 \in B,\ L_0 b_0 = 0\}$, obviously comprising leading coefficients of 0-th order according to $\bar{R}_0 = \{L_0 b_1 + L_1 b_0 \mid b_1 \in B, b_0 = 0\} \subset \bar{R}_1$. Now, it can be shown that in case of $\bar{R}_1 = \bar{B}$, we obtain again continuation by $R[L(\varepsilon)] = \bar{R}_1 = \bar{B},\ \varepsilon \neq 0$, or conversely, $\bar{R}_1$ represents the limit space of $R[L(\varepsilon)]$ during the passage $\varepsilon \to 0$.

Now, the general rule of this passage reads as follows. If the leading coefficients stabilize at $k \geq 0$ in the sense that leading coefficients of order $k + l,\ l \geq 1$ do not anymore increase the subspace of leading coefficients according to

$$\bar{R}_{k-1} \subsetneq \bar{R}_k = \bar{R}_{k+l} \quad for \quad l \geq 1 \quad and \quad \bar{R}_{-1} \coloneqq \{0\}, \tag{2.2}$$

then $\bar{R}_k$ will represent the limit space of $R[L(\varepsilon)]$ as $\varepsilon \to 0$.

If an element $\bar{b} \in \bar{B}$ is a leading coefficient, we use the abbreviation $lc(\bar{b}) = i$ in case of $\bar{b} \in \bar{R}_i$ and $\bar{b} \notin \bar{R}_{i-1}$.

Now, the behaviour of the kernels $N[L(\varepsilon)] \subset B$ is closely related to the behaviour of $R[L(\varepsilon)] \subset \bar{B}$. First, note that by (2.1), a leading coefficient of order $k \geq 1$ results from a curve $b(\varepsilon)$, which is an approximation of order $k$ with respect to the equation $L(\varepsilon) \cdot b = 0$, i.e. a $k$-th order leading coefficient arises from a constellation of the following form

$$L(\varepsilon) \cdot b(\varepsilon) = \underbrace{\sum_{l=0}^{k-1} \varepsilon^l \cdot \sum_{i+j=l} L_i b_j}_{=0} + \varepsilon^k \cdot \underbrace{(L_0\, L_1 \cdots\, L_k)\begin{pmatrix} b_k \\ b_{k-1} \\ \vdots \\ b_0 \end{pmatrix}}_{\substack{k-th\ order \\ leading\ coefficients}} + \sum_{l=k+1}^{\infty} \varepsilon^l \cdot \sum_{i+j=l} L_i b_j$$

$$=: \varepsilon^k \cdot \bar{b}(\varepsilon) \tag{2.3}$$

and we see that first $k$ summands have to vanish for a $k$-th order leading coefficient to appear.

The coefficients $(b_0 \cdots b_{k-1}),\ b_0 \neq 0$ of an approximation $b(\varepsilon) = b_0 + \cdots + \varepsilon^{k-1} \cdot b_{k-1}$ of order $k \geq 1$ is called a Jordan chain of length $k$ (or chain of generalized eigenvectors of the eigenvalue $\varepsilon = 0$). The base element $b_0 \in N[L(0)],\ b_0 \neq 0$ is called the root element of the Jordan chain.

Finally, the maximal order of an approximation that can be constructed out of $b \in B$ is called the rank of $b$ with abbreviation $rk(b)$. Note that $rk(b) = 0$, if $b \notin N[L(0)]$ and $rk(0) = \infty$.

In this sense, leading coefficients up to the order of $k \geq 1$ are essentially found in $\bar{B}$, if all root elements with rank equal to $k$ (and corresponding Jordan chains) are calculated in $B$, and we will see that stabilization at $k$ in (2.2) is equivalent to the nonexistence of root elements with finite rank above $k$.

But then, the subspace $N_{k+1} \subset N[L(0)]$ of root elements with infinite rank, allowing approximations of arbitrary high order, represents a first candidate concerning the limit space of kernels $N[L(\varepsilon)]$ as $\varepsilon \to 0$.

In fact, from [Kaballo] it is well known that $N_{k+1}$ can analytically be continued to $N[L(\varepsilon)]$ , $\varepsilon \neq 0$ and one of our aims merely consists in deriving a Taylor expansion of $N[L(\varepsilon)]$ parametrized by $N_{k+1}$ according to

$$L(\varepsilon) \cdot \underbrace{\phi(\varepsilon) \cdot N_{k+1}}_{= N[L(\varepsilon)]} = 0 \qquad with \qquad \phi(\varepsilon) = I_B + \textstyle\sum_{i=1}^{\infty} \varepsilon^i \cdot \phi_i \tag{2.4}$$

and $\phi \in C^{\omega}(U, L[B,B])$. The mapping $\phi(\varepsilon) \in L[B,B]$ may be interpreted as an $\varepsilon$-dependent, isomorphic transformation of the Banach space $B$ derived from approximations $b(\varepsilon)$ of increasing order. In some more detail, for $k \geq 0$ the analytic transformation $\phi(\varepsilon)$ implies a direct sum decomposition of $B$

$$B = \underbrace{N_1^c}_{rk=0} \oplus \underbrace{N_2^c}_{rk=1} \oplus \cdots \oplus \underbrace{N_{k+1}^c}_{rk=k} \oplus \underbrace{N_{k+1}}_{rk=\infty} \tag{2.5}$$

with $rk(n_{i+1}^c) = i$ and $\phi(\varepsilon) \cdot n_{i+1}^c$ defining approximations of order $i$ by

$$L(\varepsilon) \cdot \phi(\varepsilon) \cdot n_{i+1}^c = \varepsilon^i \cdot \bar{b}(\varepsilon) \quad and \quad \bar{b}(0) \neq 0 \tag{2.6}$$

with $i = 0, \dots, k$, $n_{i+1}^c \in N_{i+1}^c$. Now by (2.5) and (2.6), the analytic transformation $\phi(\varepsilon)$ delivers precise information about the $\varepsilon$-expansion rates occuring in $\bar{B}$ by application of $L(\varepsilon)$ to the subspaces $\phi(\varepsilon) \cdot N_{i+1}^c$. More precisely, the coefficients $S_i$ of the transformed family

$$S(\varepsilon) \coloneqq L(\varepsilon) \cdot \phi(\varepsilon) = \sum_{i=0}^{\infty} \varepsilon^i \cdot S_{i+1} \tag{2.7}$$

are buildung up a direct sum of $\bar{B}$ according to the lower part of the following diagram.

$$\begin{array}{ccccccccccc}
B & = & \overbrace{N_1^c}^{rk=0} & \oplus & \overbrace{N_2^c}^{rk=1} & \oplus \cdots \oplus & \overbrace{N_{k+1}^c}^{rk=k} & \oplus & \overbrace{N_{k+1}}^{rk=\infty} \\
 & & \uparrow & & \uparrow & & \uparrow & & \\
S(\varepsilon) & = & \boxed{S_1} & + & \varepsilon^1 \cdot \boxed{S_2} & + \cdots + & \varepsilon^k \cdot \boxed{S_{k+1}} & + & \boxed{\sum_{i=k+1}^{\infty} \varepsilon^i \cdot S_{i+1}} \\
 & & \downarrow & & \downarrow & & \downarrow & & \\
\bar{B} & = & \underbrace{R_1}_{lc=0} & \oplus & \underbrace{R_2}_{lc=1} & \oplus \cdots \oplus & \underbrace{R_{k+1}}_{lc=k} & \oplus & R_{k+1}^c
\end{array} \tag{2.8}$$

In some more detail, the coefficient $S_i$ , $i = 1, \dots, k+1$ is mapping elements from $N_i^c$ of rank $i-1$ isomorphically to a subspace $R_i$ in $\bar{B}$ composed of leading coefficients of order $i-1$. In addition, the subspaces $R_1, \ldots, R_{k+1}$ are in direct sum, i.e. the elements in $R_i$ represent no leading coefficients of order below $i-1$. The subspace $R_{k+1}^c$ represents an arbitrary complement of $R_1 \oplus \cdots \oplus R_{k+1}$.

Finally, $S_i$ is mapping root elements with rank above $i-1$ to zero, ending up with the following operator matrix representation of the normal form $S(\varepsilon)$ with respect to the decompositions of $B$ and $\bar{B}$ from (2.8).

(2.9)

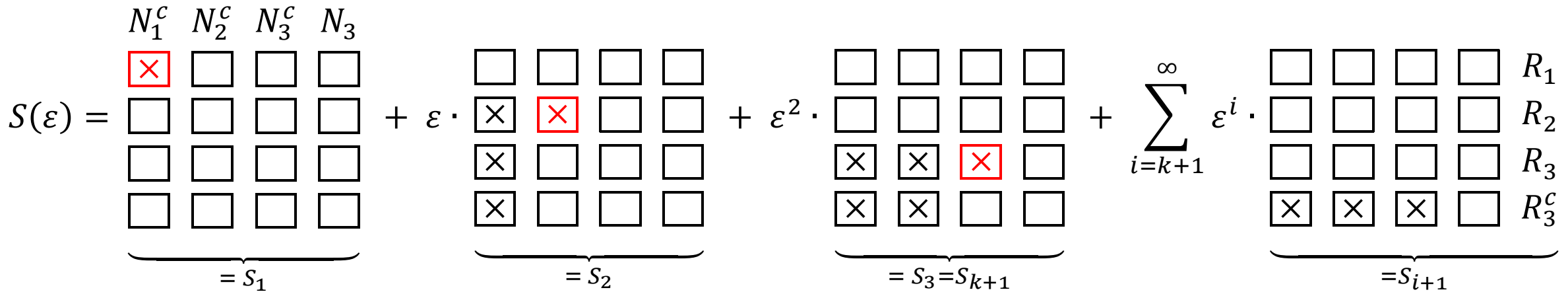

Formula (2.9) shows the case $k = 2$. Squares without entry denote the zero operator, red marked squares indicate bijection according to $S_i \in GL[N_i^c, R_i], i = 1, \ldots, k+1$. Black crosses denote possible entry into the associated subspaces of $\bar{B}$ (indicated on the right).

For $i \geq k+1$, the range of the operator $S_{i+1}$ satisfies $R[S_{i+1}] \subset R_{k+1}^c$, and if $R_{k+1}^c = \{0\}$ occurs, then the normal form $S(\varepsilon)$ reduces to a polynomial with respect to $\varepsilon$ of degree $k$, i.e. the operator family $L(\varepsilon)$ is transformed to the polynomial

$$S(\varepsilon) = L(\varepsilon) \cdot \phi(\varepsilon) = S_1 + \cdots + \varepsilon^k \cdot S_{k+1} \tag{2.10}$$

by pure, analytic right transformation $\phi(\varepsilon)$. The remainder $\sum_{i=k+1}^{\infty} \varepsilon^i \cdot S_{i+1}$ vanishes completely, where this constellation occurs, whenever the family $L(\varepsilon)$ loses surjectivity at most at $\varepsilon = 0$, then yielding an isolated singularity at $\varepsilon = 0$.

The normal form (2.8) shows in detail, the speed of expansion associated to a subspace of $B$ upon mapped by $S(\varepsilon)$ to $\bar{B}$. Exemplarily, the subspace $N_2^c \subset N[S(0)]$ is mapped by $S(\varepsilon)$ to $\bar{B}$ according to

$$S(\varepsilon) \cdot n_2^c = \underbrace{\varepsilon^1 \cdot S_2 \cdot n_2^c}_{\in R_2} + \varepsilon^2 \cdot \underbrace{r(\varepsilon) \cdot n_2^c}_{\in R_3 \oplus \cdots \oplus R_{k+1}^c} \in R_2 \oplus R_3 \oplus \cdots \oplus R_{k+1}^c , \tag{2.11}$$

i.e. at $\varepsilon = 0$, the image of $N_2^c$ equals $\{0\}$ and is subsequently blown up to a family of subspaces with speed of expansion given by $\varepsilon^1$. This means geometrically, if $N_2^c$ is restricted to the unit sphere $\|n_2^c\| \leq 1$, then the image of the sphere in $\bar{B}$ is growing with order of $\varepsilon^1$, as qualitatively indicated in figure 2 below. Moreover, the family of subspaces created from $N_2^c$ is completely contained in the subspace $R_2 \oplus R_3 \oplus \cdots \oplus R_{k+1}^c$ of $\bar{B}$.

Figure 2 visualizes the effect of the mappings $S(\varepsilon)$ within a finite dimensional setting.

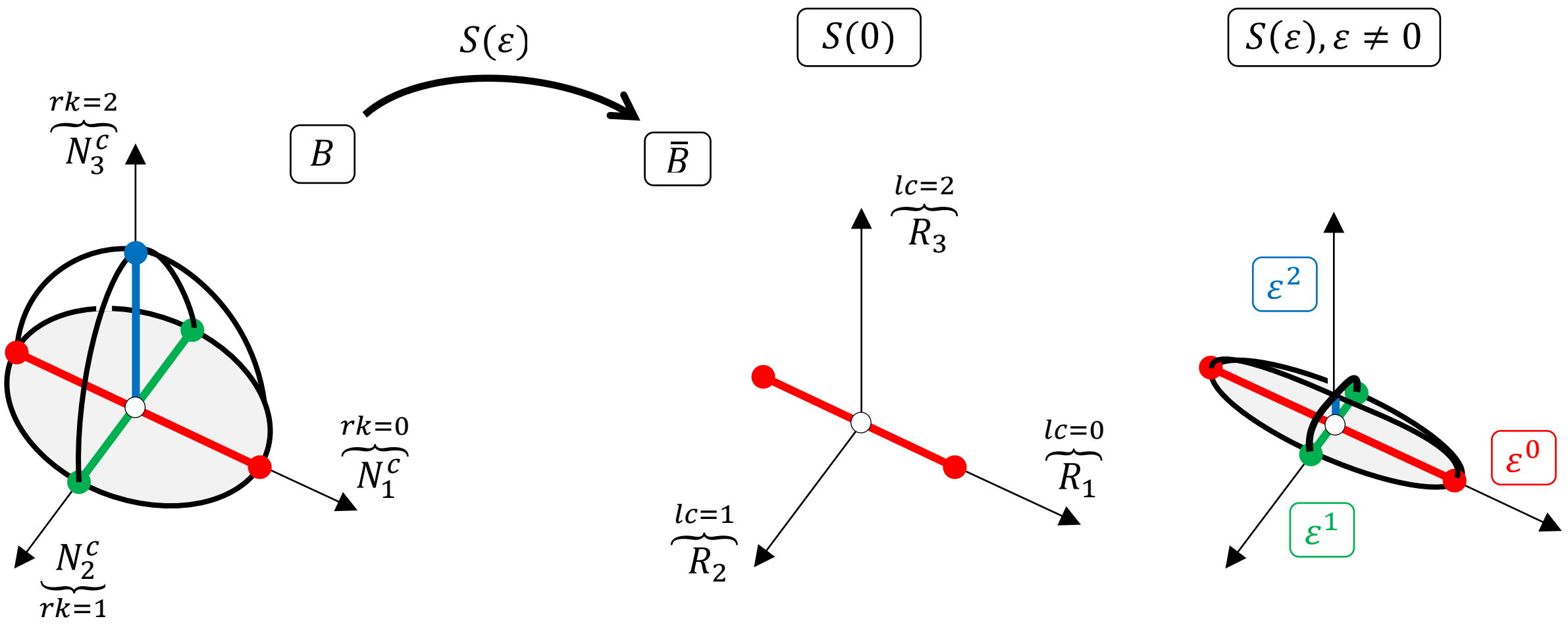

Figure 2 : Geometrical interpretation of the normal form mapping $S(\varepsilon)$.

On the left hand side, the unit sphere in the domain $B$ is indicated, as well as first three subspaces $N_1^c, N_2^c$ and $N_3^c$ with corresponding intervals in the unit sphere along the coordinate axes.

Now, when mapped by $S(\varepsilon)$ to the target space $\bar{B}$, the image of the unit sphere under $S(0)$ is simply given by the red marked line along the first subspace $R_1$ of leading coefficients of order 0 (middle diagram). On the right, we see the constellation during passage to $\varepsilon \neq 0$. Then, some kind of ellipsoid is created out of the red line with speed of expansion by order of $\varepsilon^1$ and $\varepsilon^2$ along the subspaces $R_2$ (green line) and $R_3$ (blue line) respectively.

The derivation of the analytic transformation $\phi(\varepsilon) = I_B + \sum_{i=1}^{\infty} \varepsilon^i \cdot \phi_i$ is performed in several steps. First, Jordan chains of increasing length are calculated and used to define the spaces $N_1^c, \ldots, N_{k+1}^c$ and $R_1^c, \ldots, R_{k+1}^c$ in (2.8) until stabilization occurs (as assumed with $k \geq 0$). Up to this point, a polynomial pre-transformation of the form

$$p_k(\varepsilon) = I_B + \varepsilon^1 \cdot \varphi_1 + \cdots + \varepsilon^k \cdot \varphi_k \in L[B, B] \tag{2.12}$$

has been constructed, implying a transformed system $S_k(\varepsilon)$ given by

$$S_k(\varepsilon) \coloneqq L(\varepsilon) \cdot p_k(\varepsilon) = S_1 + \varepsilon^1 \cdot S_2 + \cdots + \varepsilon^k \cdot S_{k+1} + \sum_{i=k+1}^{\infty} \varepsilon^i \cdot Q_{i+1} \tag{2.13}$$

with $S_1, \ldots, S_{k+1}$ already showing the desired pattern from (2.9). On contrary, the remainders $Q_{i+1} \in L[B, \bar{B}]$, $i \geq k+1$ are still fully occupied and show no special structure.

If $R_{k+1}^c = \{0\}$ and $N_{k+1} = \{0\}$ is assumed, then the transformed system $S_k(\varepsilon) = L(\varepsilon) \cdot p_k(\varepsilon)$ is a transversalization in the sense of [Esquinas] and [López-Gómez]. We note that the transformation used in [Esquinas], [López-Gómez] for obtaining a transversalization is given by a polynomial in $\varepsilon$ of degree $\frac{1}{2}(k + k^2)$ and a simplification is achieved by using $p_k(\varepsilon)$ instead, a polynomial of degree $k$. The polynomial $p_k(\varepsilon)$ is extracted from Jordan chains of length $k+1$.

Next, the polynomial $p_k(\varepsilon) \in L[B, B]$ is refined and eventually extended to a formal power series $\phi(\varepsilon) = I_B + \sum_{i=1}^{\infty} \varepsilon^i \cdot \phi_i$ in such a way that the leading terms $S_1 + \varepsilon^1 \cdot S_2 + \cdots + \varepsilon^k \cdot S_{k+1}$ of $S_k(\varepsilon)$ are not destroyed and all remainder terms $Q_{i+1}$, $i \geq k+1$ turn into operators $S_{i+1}$ of the transformed system $S(\varepsilon) = L(\varepsilon) \cdot \phi(\varepsilon) = S_1 + \varepsilon^1 \cdot S_2 + \cdots + \varepsilon^k \cdot S_{k+1} + \sum_{i=k+1}^{\infty} \varepsilon^i \cdot S_{i+1}$, which are mapping $B$ exclusively into the subspace $R_{k+1}^c \subset \bar{B}$. Additionally, the operators $S_{i+1}$, $i \geq k+1$ are mapping $N_{k+1}$ to zero.

In summary, by construction of $\phi(\varepsilon)$ the transformed remainder terms $S_{i+1}$ are forced to satisfy

$$P_{R_1 \oplus \cdots \oplus R_{k+1}} \cdot S_{i+1} = 0 \qquad and \qquad S_{i+1} \cdot N_{k+1} = 0 \tag{2.14}$$

for $i \geq k+1$, implying the sparse operator matrix pattern from (2.9). Here $P_{R_1 \oplus \cdots \oplus R_{k+1}}$ denotes the projection to $R_1 \oplus \cdots \oplus R_{k+1}$ related to the decomposition of $\bar{B}$ in (2.8).

In the special situation of $N_{k+1} = \{0\}$, we obtain injectivity of $L(\varepsilon)$ for $\varepsilon \neq 0$ and the second condition in (2.14) is trivially satisfied. If $R_{k+1}^c = \{0\}$, surjectivity of $L(\varepsilon)$ occurs for $\varepsilon \neq 0$ and the first condition in (2.14) turns into $S_{i+1} = 0$, i.e. the formal power series $\phi(\varepsilon)$ is constructed in such a way that the remainder terms $Q_{i+1}$ of $S_k(\varepsilon)$ in (2.13) are completely cancelled.

In some more detail, the process to achieve (2.14) is performed along the following lines. For turning $Q_{k+2}$ into $S_{k+2}$, last $k$ summands of the polynomial $p_k(\varepsilon) = I_B + \varepsilon^1 \cdot \varphi_1 + \cdots + \varepsilon^k \cdot \varphi_k$ from (2.12) are first refined and one further monomial of degree $k+1$ has to be added, yielding a polynomial $p_{k+1}(\varepsilon)$ of the form

$$p_{k+1}(\varepsilon) = I_B + \overbrace{\varepsilon^1 \cdot \underbrace{\bar{\varphi}_1}_{=:\phi_1} + \cdots + \varepsilon^k \cdot \bar{\varphi}_k}^{k\ coefficients,\ modified\ from\ p_k(\varepsilon)} + \overbrace{\varepsilon^{k+1} \cdot \varphi_{k+1}}^{new}$$

$$= \underbrace{I_B + \varepsilon^1 \cdot \phi_1}_{2\ final\ coefficients} + \varepsilon^2 \cdot \bar{\varphi}_2 + \cdots + \varepsilon^k \cdot \bar{\varphi}_k + \varepsilon^{k+1} \cdot \varphi_{k+1}$$

of degree $k+1$ with $\bar{\varphi}_1, \ldots, \bar{\varphi}_k$ being modified versions of $\varphi_1, \ldots, \varphi_k$. The coefficient $\bar{\varphi}_1 =: \phi_1$ has already reached its final form, whereas $\varphi_{k+1}$ is the new coefficient constructed from the decompositions of $B$ and $\bar{B}$ in (2.8). Then, the transformed system reads

$$L(\varepsilon) \cdot p_{k+1}(\varepsilon) = \overbrace{S_1 + \varepsilon^1 \cdot S_2 + \cdots + \varepsilon^k \cdot S_{k+1} + \varepsilon^{k+1} \cdot \underbrace{S_{k+2}}_{new}}^{k+2\ final\ coefficients} + \overbrace{\textstyle\sum_{i=k+2}^{\infty} \varepsilon^i \cdot Q_{i+1}}^{remainder}$$

with new operator $S_{k+2}$ satisfying (2.14) as desired, and in total yielding $k+2$ final coefficients.

This process can be continued up to infinity in the sense that for turning $Q_{k+2+l}$ into $S_{k+2+l}$, $l \geq 0$, last $k$ summands of the previous polynomial $p_{k+l}(\varepsilon)$ have to be refined and one further monomial of degree $k+1+l$ has to be added, implying a transformation $p_{k+1+l}(\varepsilon)$ of the form

$$p_{k+1+l}(\varepsilon) = I_B + \varepsilon^1 \cdot \phi_1 + \cdots + \varepsilon^l \cdot \phi_l + \overbrace{\varepsilon^{1+l} \cdot \underbrace{\bar{\varphi}_{1+l}}_{=:\phi_{1+l}} + \cdots + \varepsilon^{k+l} \cdot \bar{\varphi}_{k+l}}^{k\ coefficients,\ modified\ from\ p_{k+l}(\varepsilon)} + \overbrace{\varepsilon^{k+1+l} \cdot \varphi_{k+1+l}}^{new}$$

$$= \underbrace{I_B + \cdots + \varepsilon^{1+l} \cdot \phi_{1+l}}_{2+l\ final\ coefficients} + \varepsilon^{2+l} \cdot \bar{\varphi}_{2+l} + \cdots + \varepsilon^{k+l} \cdot \bar{\varphi}_{k+l} + \varepsilon^{k+1+l} \cdot \varphi_{k+1+l} \tag{2.15}$$

with corresponding transformed system

$$L(\varepsilon)\cdot p_{k+1+l}(\varepsilon) = \overbrace{S_1 + \varepsilon^1\cdot S_2 + \cdots + \varepsilon^{1+l}\cdot S_{2+l} + \varepsilon^{2+l}\cdot S_{3+l} + \cdots + \varepsilon^{k+1+l}\cdot \underbrace{S_{k+2+l}}_{new}}^{k+2+l\ final\ coefficients}$$

$$+ \overbrace{\sum_{i=k+2+l}^{\infty} \varepsilon^i \cdot Q_{i+1}}^{remainder} . \tag{2.16}$$

Note that, after each transformation step, there remains always a gap of $k$ between the number of final coefficients in the transformation (2.15) and the transformed system (2.16).

Note also that for calculation of first $k$ coefficients $\phi_1, \dots, \phi_k$, the basic pre-transformation $p_k(\varepsilon)$ has to be extended up to the polynomial $p_{2k}(\varepsilon)$ of degree $2k$ (use $l = k - 1$ in (2.15)), where this extension is still influenced by the coefficients $\varphi_1, \dots, \varphi_k$ of the pre-transformation $p_k(\varepsilon)$, dependent from the initial step by step build up of the direct sums in (2.8), i.e. no uniform formulas can be expected concerning the calculation of $\phi_1, \dots, \phi_k$.

On contrary, the calculation of higher order coefficients $\phi_{k+1}$, $\phi_{k+2}, \dots$ merely depends of the stabilized decompositions in (2.8) and can iteratively be described by a rather simple and uniform recursive scheme of the form

$$\phi_{k+i} = F_i(\phi_1, \dots, \phi_{k+i-1}) \,,\ i \geq 1 \,, \tag{2.17}$$

ending up with a well defined formal power series $\phi(\varepsilon) = I_B + \sum_{i=1}^{\infty} \varepsilon^i \cdot \phi_i$ satisfying by Cauchy product $L(\varepsilon)\cdot\phi(\varepsilon) = S(\varepsilon)$ with $S(\varepsilon)$ showing the pattern from (2.9).

It remains to show the convergence of $\phi(\varepsilon)$. For this purpose, we ensure that the power series $\phi(\varepsilon)$ satisfies an equation of the form

$$[\, I_B - \varepsilon\cdot f(\varepsilon)\,]\cdot\phi(\varepsilon) = g(\varepsilon)\,, \tag{2.18}$$

with analytic mappings $f \in C^{\omega}(U, L[B,B])$ and $g \in C^{\omega}(U,B)$. Analyticity of $f, g$ results from analyticity of $L$ and the assumption that the projections associated to decompositions (2.8) are continous. Now, equation (2.18) has a unique analytic solution, which can be represented by Neumann series according to

$$\phi(\varepsilon) = [\, I_B - \varepsilon\cdot f(\varepsilon)\,]^{-1}\cdot g(\varepsilon) \tag{2.19}$$

$$= [\, I_B + \varepsilon\cdot f(\varepsilon) + \varepsilon^2\cdot f(\varepsilon)^2 + \cdots ]\cdot g(\varepsilon)\,.$$

Using $f(\varepsilon) = \sum_{i=0}^{\infty} \varepsilon^i \cdot f_i$ and $g(\varepsilon) = \sum_{i=0}^{\infty} \varepsilon^i \cdot g_i$, every coefficient $\phi_i$ of $\phi(\varepsilon) = I_B + \sum_{i=1}^{\infty} \varepsilon^i \cdot \phi_i$ can be formulated in an explicit way by given coefficients $L_i$.

If analyticity of the power series $L(\varepsilon) = \sum_{i=0}^{\infty} \varepsilon^i \cdot L_i$ is not assumed, then the complete process is working along the same lines of reasoning with formal power series $f(\varepsilon)$ and $g(\varepsilon)$ in (2.18) and corresponding formal power series $\phi(\varepsilon)$ in (2.19). In addition, the transformed power series $S(\varepsilon) = L(\varepsilon)\cdot\phi(\varepsilon)$ causes again the normal form pattern from (2.9) and in this sense, the complete construction shows up to be, in essential parts, mainly a topic concerning formal power series, with convergent power series being a subset.

In the context of formal versus analytic transformations, we also refer to [Walcher] with resonant eigenvalues allowing formal, but preventing analytic transformation of vector fields to normal form.

Given a formal power series of linear operators $L(\varepsilon) = \sum_{i=0}^{\infty} \varepsilon^i \cdot L_i\,,\ L_i \in L(B, \bar{B})$ with $B, \bar{B}$ real or complex vector spaces, $\varepsilon \in \mathbb{K} = \mathbb{R}, \mathbb{C}$ and $L(B, \bar{B})$ denoting the vector space of linear mappings, not necessarily bounded, between $B$ and $\bar{B}$. Then, a formal power series $b(\varepsilon) = \sum_{i=0}^{\infty} \varepsilon^i \cdot b_i\,,\ b_i \in B$ is an approximation of order $k \geq 1$ with respect to the equation $L(\varepsilon) \cdot b = 0$, if $b(\varepsilon)$ satisfies

$$L(\varepsilon) \cdot b(\varepsilon) = \underbrace{L_0 b_0 + \cdots + \varepsilon^{k-1} \cdot (L_0 \cdots L_{k-1}) \cdot \begin{pmatrix} b_{k-1} \\ \vdots \\ b_0 \end{pmatrix}}_{=0} + \sum_{l=k}^{\infty} \varepsilon^l \cdot \sum_{i+j=l} L_i b_j = \varepsilon^k \cdot \bar{b}(\varepsilon) \qquad (3.1)$$

or equivalently

$$\underbrace{\begin{pmatrix} L_0 & \cdots & L_{k-1} \\ & \ddots & \vdots \\ & & L_0 \end{pmatrix}}_{=:\,\Delta^k} \cdot \begin{pmatrix} b_{k-1} \\ \vdots \\ b_0 \end{pmatrix} = 0 \quad \Leftrightarrow \quad \begin{pmatrix} b_{k-1} \\ \vdots \\ b_0 \end{pmatrix} \in N[\,\Delta^k\,] \qquad (3.2)$$

with $\Delta^k \in L(B^k, \bar{B}^k)$. Hence, to determine approximations $b(\varepsilon)$ of increasing order, the kernels of $\Delta^k$ have to be calculated for $k = 1,2, \dots$ . For this purpose define

$$\begin{array}{lllll} N_0 := B & N_0^c := \{0\} & \boxed{S_1 := L_0} & N_1 := N[S_1] & N_1^c := N_0/N_1 \\ R_0 := \{0\} & R_0^c := \bar{B} & & R_1 := R[S_1] & R_1^c := R_0^c/R_1 \end{array} \qquad (3.3)$$

with quotient spaces $N_1^c = N_0/N_1$ and $R_1^c = R_0^c/R_1$, implying algebraic direct sum decompositions of $B = N_0$ and $\bar{B} = R_0^c$ according to

$$\begin{array}{ccccccc} B = & \overbrace{N_0}^{rk \geq 0} & = & \overbrace{N_1^c}^{rk=0} & \oplus & \overbrace{N_1}^{rk \geq 1} & \\ & & & \uparrow & & & \\ & & & \boxed{S_1 = L_0} & & & \\ & & & \downarrow & & & \\ \bar{B} = & \underbrace{R_0^c}_{lc \geq 0} & = & \underbrace{R_1}_{lc=0} & \oplus & \underbrace{R_1^c}_{lc \geq 1} & . \end{array} \qquad (3.4)$$

Note that by (3.4) a split of the vector spaces $B$ and $\bar{B}$ is defined, guided by the geometric behaviour with respect to the rank in $B$ and leading coefficients in $\bar{B}$, where the subspaces $N_1^c$ with $rk(n_1^c) = 0$ and $R_1$ with $lc(r_1) = 0$ are in bijection by $S_1$, as indicated by arrows in (3.4).

Then, the approximations of order $k = 1$, characterized by $N[\Delta^1]$, may be written (in a rather complicated way) by

$$N[\,\Delta^1\,] = N[\,L_0\,] = N_1 = R\left[\,\boxed{M}^{\,1}_{\,|\,N_1}\,\right] \qquad (3.5)$$

under consideration of the setting $\boxed{M}^{\,1} := I_B \in L(B^1, B^1)$. Now, for arbitrary $k \geq 1$, we will see that approximations of order $k+1$ are given by

$$N[\,\Delta^{k+1}\,] = R\left[\,\boxed{M}^{\,k+1}_{\,|\,N_1 \times \cdots \times N_{k+1}}\,\right], \qquad (3.6)$$

with a $(k+1) \times (k+1)$ upper triangular operator matrix $\boxed{M}^{k+1} \in L(B^{k+1}, B^{k+1})$ with diagonal composed of the identity map in $B$ according to

$$\boxed{M}^{k+1} = \begin{pmatrix} I_B & * & * \\ & \ddots & * \\ & & I_B \end{pmatrix}. \tag{3.7}$$

The operator matrix $\boxed{M}^{k+1}$ is iteratively defined for $k = 1,2,\ldots$ in the following way. First define the operator

$$\bar{S}_{k+1} \coloneqq [\, L_1, \ldots, L_k \,] \cdot M_k \;\in\; L(B, \bar{B}) \tag{3.8}$$

with $M_k$ denoting the last column $k$ of the previous operator $\boxed{M}^{k}$, i.e. for $k = 1$ we simply obtain $\bar{S}_2 \coloneqq [L_1] \cdot M_1 = L_1 \cdot I_B = L_1$. Next, the definitions from (3.3), related to $S_1$, are transferred to $S_{k+1}$ by

$$\boxed{S_{k+1} \coloneqq P_{R_k^c}\, \bar{S}_{k+1} \in\; L(B, R_k^c)} \qquad \begin{matrix} N_{k+1} \coloneqq N\big[S_{k+1}|_{N_k}\big] & N_{k+1}^c \coloneqq N_k / N_{k+1} \\ R_{k+1} \coloneqq R\big[S_{k+1}|_{N_k}\big] & R_{k+1}^c \coloneqq R_k^c / R_{k+1} \end{matrix} \tag{3.9}$$

with corresponding algebraic direct sum decompositions of $N_k$ and $R_k^c$ according to

$$\begin{matrix} \overbrace{N_k}^{rk\geq k} & = & \overbrace{N_{k+1}^c}^{rk=k} & \oplus & \overbrace{N_{k+1}}^{rk\geq k+1} \\ & & \uparrow & & \\ & & \boxed{S_{k+1}} & & \\ & & \downarrow & & \\ \underbrace{R_k^c}_{lc\geq k} & = & \underbrace{R_{k+1}}_{lc=k} & \oplus & \underbrace{R_{k+1}^c}_{lc\geq k+1} \end{matrix} \quad . \tag{3.10}$$

Again, the split of the subspaces $N_k$ and $R_k^c$ is geometrically guided by subspaces $N_{k+1}^c$ and $R_{k+1}$ with rank $rk(n_{k+1}^c) = k$ and leading coefficient number $lc(r_{k+1}) = k$, which are in bijection by the operator $S_{k+1}$. For consistency, we still add $\bar{S}_1 \coloneqq L_0$ and note $S_1 = P_{R_0^c}\, \bar{S}_1 = P_{\bar{B}}\, L_0 = L_0$.

At this stage, the following decompositions of $B$ and $\bar{B}$ are achieved.

$$
\begin{array}{l}
B = \overbrace{N_1^c \oplus \overbrace{N_2^c \oplus \cdots \oplus \overbrace{N_{k+1}^c \oplus {\color{red}N_{k+1}}}^{=N_k}}^{=N_1}}^{=N_0} \\
\qquad \uparrow \qquad \uparrow \qquad\qquad \uparrow \\
\quad\ \boxed{S_1} \quad \boxed{S_2} \qquad \boxed{S_{k+1}} \\
\qquad \downarrow \qquad \downarrow \qquad\qquad \downarrow \\
\bar{B} = \underbrace{R_1 \oplus \underbrace{R_2 \oplus \cdots \oplus \underbrace{R_{k+1} \oplus {\color{red}R_{k+1}^c}}_{=R_k^c}}_{=R_1^c}}_{=R_0^c}
\end{array}
\tag{3.11}
$$

Now, for closing the recursion, it remains to define the operator matrix $\boxed{M}^{k+1} \in L\left(B^{k+1}, B^{k+1}\right)$ from (3.7) for $k \geq 1$. First, set $E_1 \coloneqq E_{1,1} \coloneqq I_B$ and for $k \geq 1$ define an operator column vector $E_{k+1} \in L\left(B, B^{k+1}\right)$ with $k+1$ components from bottom to top by

$$E_{k+1,k+1} \coloneqq I_B \in L(B,B)$$

$$E_{i,k+1} \coloneqq -S_i^{-1}\,\mathcal{P}_i \cdot \sum_{\nu=i+1}^{k+1} \bar{S}_\nu \cdot E_{\nu,k+1} \in L(B,B)\,, \quad i = k, \dots, 1\,. \tag{3.12}$$

Here, $\mathcal{P}_i \in L(\bar{B}, \bar{B})$ denotes the projection to $R_i$ with respect to the direct sum of $\bar{B}$ in (3.11) and $S_i^{-1} \in L(R_i, N_i^c)$ represents the inverse of $S_i$ with respect to the subspaces $R_i$ and $N_i^c$ satisfying $lc(r_i) = i - 1$ and $rk(n_i^c) = i - 1$ respectively.

Exemplarily, first 3 column vectors $E_1, E_2, E_3$, collected within a common matrix $\boxed{E}^3$, read

$$
\boxed{E}^3 \coloneqq \overbrace{\begin{pmatrix} E_{1,1} & E_{1,2} & E_{1,3} \\ & E_{2,2} & E_{2,3} \\ & & E_{3,3} \end{pmatrix}}^{=:\, trng[\,E_1, E_2, E_3\,]} \tag{3.13}
$$

$$
= \begin{pmatrix} I_B & -S_1^{-1}\,\mathcal{P}_1 \cdot \bar{S}_2 & -S_1^{-1}\,\mathcal{P}_1 \cdot [\,I_B - \bar{S}_2 \cdot S_2^{-1}\,\mathcal{P}_2\,] \cdot \bar{S}_3 \\ & I_B & -S_2^{-1}\,\mathcal{P}_2 \cdot \bar{S}_3 \\ & & I_B \end{pmatrix}.
$$

In general, we use the abbreviation $\boxed{E}^{k+1} \coloneqq trng[E_1, \dots, E_{k+1}]$ for $k \geq 0$.

Now, the last column of the operator matrix $\boxed{M}^{k+1}$ is well defined by

$$M_{k+1} := \begin{pmatrix} I_B & | & \\ - & - & - \\ & | & \boxed{M}^k \end{pmatrix} \cdot E_{k+1} \in L\left(B, B^{k+1}\right), \tag{3.14}$$

whereas first $k$ leading columns of $\boxed{M}^{k+1}$ are simply given by $\boxed{M}^k$, i.e. by last columns $M_1, \dots, M_k$ of previously calculated matrix operators, finally yielding

$$\boxed{M}^{k+1} := \begin{pmatrix} \underbrace{M_{1,1}}_{I_B} & M_{1,2} & \cdots & M_{1,k+1} \\ & \underbrace{M_{2,2}}_{I_B} & \cdots & M_{2,k+1} \\ & & \ddots & \vdots \\ & & & \underbrace{M_{k+1,k+1}}_{I_B} \end{pmatrix} \in L\left(B^{k+1}, B^{k+1}\right). \tag{3.15}$$

Again, we will use the abbrevation $\boxed{M}^{k+1} = trng[M_1, \dots, M_{k+1}]$ and note that by direct inspection $M_{i,i} = I_B, i = 1, \dots, k+1$.

Next, by induction we easily show $N[\Delta^{k+1}] = R[\boxed{M}^{k+1}{}_{| N_1 \times \cdots \times N_{k+1}}], k \geq 0$, as stated in (3.6). For $k = 0$, the settings in (3.5) yield the assertion. Then, assume $N[\Delta^k] = R[\boxed{M}^k{}_{| N_1 \times \cdots \times N_k}], k \geq 0$. For calculation of $N[\Delta^{k+1}]$, we have to plug $N[\Delta^k]$ into the $k$-th leading coefficient of the Cauchy product in (3.1) according to

$$\begin{aligned} (L_0 \cdots L_k) \cdot \underbrace{\begin{pmatrix} b_k \\ \vdots \\ b_0 \end{pmatrix}}_{\in B \times N[\Delta^k]} &= L_0 \cdot \underbrace{b_k}_{=: n_0 \in B} + (L_1 \cdots L_k) \cdot \boxed{M}^k \cdot \underbrace{\begin{pmatrix} n_1 \\ \vdots \\ n_k \end{pmatrix}}_{\in N_1 \times \cdots \times N_k} \\ &= L_0 \cdot n_0 + (L_1 \cdots L_k) \cdot trng[M_1, \dots, M_k] \cdot \begin{pmatrix} n_1 \\ \vdots \\ n_k \end{pmatrix} \\ &\overset{(3.8)}{=} \bar{S}_1 \cdot n_0 + (\bar{S}_2 \cdots \bar{S}_{k+1}) \cdot \begin{pmatrix} n_1 \\ \vdots \\ n_k \end{pmatrix} \\ &= (\bar{S}_1 \cdots \bar{S}_{k+1}) \cdot \begin{pmatrix} n_0 \\ \vdots \\ n_k \end{pmatrix} \end{aligned} \tag{3.16}$$

and equate it to zero. By decomposition (3.11), we obtain the equivalences

$$(\bar{S}_1 \cdots \bar{S}_{k+1}) \cdot \begin{pmatrix} n_0 \\ \vdots \\ n_k \end{pmatrix} = 0 \in \bar{B}$$

$$\Leftrightarrow \quad \mathcal{P}_i \cdot (\bar{S}_1 \cdots \bar{S}_{k+1}) \cdot \begin{pmatrix} n_0 \\ \vdots \\ n_k \end{pmatrix} = 0\,, i = 1, \dots, k+1 \quad \wedge \quad P_{R_{k+1}^c} \cdot (\bar{S}_1 \cdots \bar{S}_{k+1}) \cdot \begin{pmatrix} n_0 \\ \vdots \\ n_k \end{pmatrix} = 0$$

$$\Leftrightarrow \begin{pmatrix} \boxed{S_1} & \mathcal{P}_1\bar{S}_2 & \cdots & \mathcal{P}_1\bar{S}_{k+1} \\ & \ddots & \ddots & \vdots \\ & & \boxed{S_k} & \mathcal{P}_k\bar{S}_{k+1} \\ & & & \boxed{S_{k+1}} \end{pmatrix} \cdot \begin{pmatrix} n_0 \\ \vdots \\ n_{k-1} \\ n_k \end{pmatrix} = 0 \tag{3.17}$$

$$\wedge \quad \bar{S}_1 n_0 + \cdots + \bar{S}_{k+1} n_k \in R_1 \oplus \cdots \oplus R_{k+1}\,,$$

where the settings in (3.9) imply $\bar{S}_1 n_0 \in R_1$ , … , $\bar{S}_{k+1} n_k \in R_1 \oplus \cdots \oplus R_{k+1}$ and we can restrict to the operator matrix equation (3.17). Then, by the definitions of the operators $E_{i,j}$ in (3.12), the following equivalences result from bottom up solution of the triangular system (3.17)

$$\Leftrightarrow \begin{cases} n_k & = \bar{n}_{k+1} & = E_{k+1,k+1} \cdot \bar{n}_{k+1}\,, & \bar{n}_{k+1} \in N_{k+1} \\ n_{k-1} & = \bar{n}_k - S_k^{-1}\mathcal{P}_k\bar{S}_{k+1} \cdot \bar{n}_{k+1} & = [\, E_{k,k} \quad E_{k,k+1} \,] \cdot \begin{pmatrix} \bar{n}_k \\ \bar{n}_{k+1} \end{pmatrix}, & \bar{n}_k \in N_k \\ \vdots & \vdots & \vdots & \vdots \end{cases}$$

$$\Leftrightarrow \begin{pmatrix} n_0 \\ \vdots \\ n_k \end{pmatrix} = \begin{pmatrix} E_{1,1} & \cdots & E_{1,k+1} \\ & \ddots & \vdots \\ & & E_{k+1,k+1} \end{pmatrix} \cdot \begin{pmatrix} \bar{n}_1 \\ \vdots \\ \bar{n}_{k+1} \end{pmatrix}, \quad \begin{pmatrix} \bar{n}_1 \\ \vdots \\ \bar{n}_{k+1} \end{pmatrix} \in N_1 \times \cdots \times N_{k+1}\,. \tag{3.18}$$

Thus, from the first equality in (3.16), we obtain the equivalences

$$\begin{pmatrix} b_k \\ \vdots \\ b_0 \end{pmatrix} \in N[\,\Delta^{k+1}\,] \tag{3.19}$$

$$\Leftrightarrow \begin{pmatrix} b_k \\ \vdots \\ b_0 \end{pmatrix} = \begin{pmatrix} I_B & | & \\ - & - & - \\ & | & \boxed{M}^k \end{pmatrix} \cdot \begin{pmatrix} n_0 \\ \vdots \\ n_k \end{pmatrix} = \begin{pmatrix} I_B & | & \\ - & - & - \\ & | & \boxed{M}^k \end{pmatrix} \cdot \begin{pmatrix} E_{1,1} & \cdots & E_{1,k+1} \\ & \ddots & \vdots \\ & & E_{k+1,k+1} \end{pmatrix} \cdot \begin{pmatrix} \bar{n}_1 \\ \vdots \\ \bar{n}_{k+1} \end{pmatrix}$$

$$= \begin{pmatrix} I_B & & & \\ & M_{1,1} & \cdots & M_{1,k} \\ & & \ddots & \vdots \\ & & & M_{k,k} \end{pmatrix} \cdot \begin{pmatrix} E_{1,1} & E_{1,2} & \cdots & E_{1,k+1} \\ & E_{2,2} & \cdots & E_{2,k+1} \\ & & \ddots & \vdots \\ & & & E_{k+1,k+1} \end{pmatrix} \cdot \begin{pmatrix} \bar{n}_1 \\ \vdots \\ \bar{n}_{k+1} \end{pmatrix}$$

$$= trng\left[ I_B \cdot E_1,\ \begin{pmatrix} I_B & \\ & \boxed{M}^1 \end{pmatrix} \cdot E_2\ , \dots ,\ \begin{pmatrix} I_B & \\ & \boxed{M}^k \end{pmatrix} \cdot E_{k+1} \right] \cdot \begin{pmatrix} \bar{n}_1 \\ \vdots \\ \bar{n}_{k+1} \end{pmatrix}$$

$$\underset{(3.15)}{\overset{(3.14)}{=}} trng[\, M_1\,, \dots, M_{k+1}\,] \cdot \begin{pmatrix} \bar{n}_1 \\ \vdots \\ \bar{n}_{k+1} \end{pmatrix}$$

$$= \boxed{M}^{k+1} \cdot \begin{pmatrix} \bar{n}_1 \\ \vdots \\ \bar{n}_{k+1} \end{pmatrix},$$

implying the following result with respect to Jordan chains of of increasing length.

**Lemma 1:** Given a power series of linear operators $L(\varepsilon) = \sum_{i=0}^{\infty} \varepsilon^i \cdot L_i$, $L_i \in L(B, \bar{B})$ with $B, \bar{B}$ real or complex vector spaces. Then, the iteration (3.3)-(3.15) is well defined for $k \geq 0$ with following properties.

(i) The kernel of $\Delta^{k+1} \in L\left(B^{k+1}, \bar{B}^{k+1}\right)$ can be represented by the range of $\boxed{M}^{k+1} \in L\left(B^{k+1}, B^{k+1}\right)$ restricted to $N_1 \times \cdots \times N_{k+1}$, i.e. we obtain for $k \geq 0$

$$N[\,\Delta^{k+1}\,] = N[\begin{pmatrix} L_0 & \cdots & L_k \\ & \ddots & \vdots \\ & & L_0 \end{pmatrix}] = R\left[\boxed{M}^{k+1}{}_{|\,N_1 \times \cdots \times N_{k+1}}\right]$$

$$= \begin{pmatrix} I_B \\ 0 \\ 0 \\ \vdots \\ 0 \end{pmatrix} \cdot N_1 \oplus \begin{pmatrix} M_{1,2} \\ I_B \\ 0 \\ \vdots \\ 0 \end{pmatrix} \cdot N_2 \oplus \cdots \oplus \begin{pmatrix} M_{1,k+1} \\ M_{2,k+1} \\ \vdots \\ M_{k,k+1} \\ I_B \end{pmatrix} \cdot N_{k+1}$$

(ii) For $k \geq 0$ :

$$rk(b) \geq k \quad \Leftrightarrow \quad b \in N_k$$

$$rk(b) = k \quad \Leftrightarrow \quad b \in N_k \backslash N_{k+1}$$

$$rk(b) = k \quad for \quad b \in N_{k+1}^c \subset N_k\,,\ \ b \neq 0$$

(iii) For $k \geq 0$ :

$$lc(\bar{b}) \leq k \quad \Leftrightarrow \quad \bar{b} \in \bar{R}_k = R_1 \oplus \cdots \oplus R_{k+1}$$

$$lc(\bar{b}) = k \quad \Leftrightarrow \quad \bar{b} \in \bar{R}_k \backslash \bar{R}_{k-1}\,,\ \ \bar{R}_{-1} = \{0\}$$

$$lc(\bar{b}) = k \quad for \quad \bar{b} \in R_{k+1} \subset \bar{R}_k\,,\ \ b \neq 0$$

(i) follows from (3.19). Concerning (ii), $rk(b) \geq k$ iff $b \neq 0$ can be extended to a $k$-tupel that lies in $N[\Delta^k]$. By triangularity of $\boxed{M}^k$ this is possible iff $b \in N_k$. The remaining assertions in (ii) follow immediately.

Concerning (iii), the leading coefficients $\bar{R}_k \subset \bar{B}$ of order $k$ are obtained by inserting $B \times N[\Delta^k], k \geq 1$ into the $k$-th leading coefficient of the Cauchy product in (3.1) implying by (3.16)

$$\bar{R}_k = (L_0 \cdots L_k) \cdot \begin{pmatrix} B \\ N[\Delta^k] \end{pmatrix} = (\,\bar{S}_1 \cdots \bar{S}_{k+1}\,) \cdot \begin{pmatrix} N_0 \\ \vdots \\ N_k \end{pmatrix} \tag{3.20}$$

$$= \underbrace{\underbrace{\underbrace{\bar{S}_1 \cdot N_0}_{=\bar{R}_0} + \bar{S}_2 \cdot N_1}_{=\bar{R}_1} + \cdots + \bar{S}_k \cdot N_{k-1}}_{=\bar{R}_{k-1}} + \bar{S}_{k+1} \cdot N_k$$

and all leading coefficients of order below $k$ are contained in $\bar{R}_k$. But then, $lc(\bar{b}) \le k$ iff $\bar{b} \in \bar{R}_k$. Further, we have

$$\bar{R}_k = \bar{S}_1 \cdot N_0 + \bar{S}_2 \cdot N_1 + \cdots + \bar{S}_k \cdot N_{k-1} + \bar{S}_{k+1} \cdot N_k \tag{3.21}$$

$$= \overbrace{\left(P_{R_0^c}\right)}^{=I_{\bar{B}}} \cdot \bar{S}_1 \cdot N_0 + \overbrace{\left(\mathcal{P}_1 + P_{R_1^c}\right)}^{=I_{\bar{B}}} \cdot \bar{S}_2 \cdot N_1 + \cdots + \overbrace{\left(\mathcal{P}_1 + \cdots + \mathcal{P}_k + P_{R_k^c}\right)}^{=I_{\bar{B}}} \cdot \bar{S}_{k+1} \cdot N_k$$

$$\overset{(3.9)}{=} \underbrace{S_1 \cdot N_0 + S_2 \cdot N_1 + \cdots + S_{k+1} \cdot N_k}_{=R_1 \oplus \cdots \oplus R_{k+1}} + \underbrace{(\mathcal{P}_1) \cdot \bar{S}_2 \cdot N_1 + \cdots + (\mathcal{P}_1 + \cdots + \mathcal{P}_k) \cdot \bar{S}_{k+1} \cdot N_k}_{\subset R_1 \oplus \cdots \oplus R_k}$$

$$= R_1 \oplus \cdots \oplus R_{k+1}$$

and the first statement in (iii) is shown. Again, the remaining statements in (iii) follow immediately.

## *4. Pre-Transformation and Transversalization*

In this section, the pre-transformation of $L(\varepsilon)$ into the form

$$S_k(\varepsilon) \coloneqq L(\varepsilon) \cdot p_k(\varepsilon) = S_1 + \varepsilon^1 \cdot S_2 + \cdots + \varepsilon^k \cdot S_{k+1} + \sum_{i=k+1}^{\infty} \varepsilon^i \cdot Q_{i+1} , \tag{4.1}$$

with operators $S_i, i = 1, \dots, k+1$ given by (3.9), is derived. We assume that the iteration is performed up to decomposition (3.15) with $k \ge 0$. Stabilization is not presupposed.

Due to $S_i \cdot N_{k+1} = 0, i = 1, \dots, k+1$, a necessary condition for $p_k(\varepsilon)$ to obtain (4.1) reads

$$L(\varepsilon) \cdot p_k(\varepsilon) \cdot n_{k+1} = \sum_{i=k+1}^{\infty} \varepsilon^i \cdot Q_{i+1} \cdot n_{k+1} \tag{4.2}$$

for $n_{k+1} \in N_{k+1}$, i.e. $b(\varepsilon) = p_k(\varepsilon) \cdot n_{k+1}$ has to define an approximation of $L(\varepsilon) \cdot b = 0$ of order $k+1$. But then, it is plausible that $p_k(\varepsilon)$ might be derived from Jordan chains of length $k+1$ given by

$$N[\,\Delta^{k+1}\,] = trng[\,M_1 , \dots, M_{k+1}\,] \cdot \begin{pmatrix} N_1 \\ \vdots \\ N_{k+1} \end{pmatrix} = \boxed{M}^{k+1} \cdot \begin{pmatrix} N_1 \\ \vdots \\ N_{k+1} \end{pmatrix}. \tag{4.3}$$

In particular, by Lemma 1 (i) the approximations with root elements different from zero are simply constructed by the last column $M_{k+1}$ of $\boxed{M}^{k+1}$ according to

$$b(\varepsilon) = (\varepsilon^k, \cdots, \varepsilon^1, 1) \cdot M_{k+1} \cdot n_{k+1} = (\varepsilon^k, \cdots, \varepsilon^1, 1) \cdot \begin{pmatrix} M_{1,k+1} \\ \vdots \\ M_{k,k+1} \\ I_B \end{pmatrix} \cdot n_{k+1}$$

$$= (I_B + \varepsilon \cdot M_{k,k+1} + \cdots + \varepsilon^k \cdot M_{1,k+1}) \cdot n_{k+1} \tag{4.4}$$

and setting

$$p_k(\varepsilon) := I_B + \varepsilon \cdot M_{k,k+1} + \cdots + \varepsilon^k \cdot M_{1,k+1} \tag{4.5}$$

implies, at least, (4.2) by construction. Next, we obtain

$$L(\varepsilon) \cdot p_k(\varepsilon)$$

$$= (L_0 + \varepsilon \cdot L_1 + \cdots + \varepsilon^k \cdot L_k + \textstyle\sum_{i=k+1}^{\infty} \varepsilon^i \cdot L_i) \cdot (I_B + \varepsilon \cdot M_{k,k+1} + \cdots + \varepsilon^k \cdot M_{1,k+1})$$

$$= \underbrace{L_0 I_B}_{=S_1} + \underbrace{\varepsilon \cdot (L_0\ L_1) \begin{pmatrix} M_{k,k+1} \\ I_B \end{pmatrix}}_{=S_2} + \cdots + \underbrace{\varepsilon^k \cdot (L_0 \cdots L_k) \overbrace{\begin{pmatrix} M_{1,k+1} \\ \vdots \\ M_{k,k+1} \\ I_B \end{pmatrix}}^{=M_{k+1}}}_{=S_{k+1}} + \textstyle\sum_{i=k+1}^{\infty} \varepsilon^i \cdot Q_{i+1} \tag{4.6}$$

with remainder coefficients $Q_{i+1} \in L(B, \bar{B}), i \geq k+1$. Now for proving (4.1), it remains to verify the red marked identities in (4.6) with the first identity $L_0 I_B = S_1$ being true by (3.3). Note also that the column vectors in (4.6) are given by an increasing number of last components of $M_{k+1}$ until the complete vector $M_{k+1}$ is build up within the coefficient of $\varepsilon^k$.

For proving the red marked identities in (4.6) we use the following lemma.

**Lemma 2:** For $k \geq 0$

$$(L_0 \cdots L_k) \cdot \boxed{M}^{k+1} = (\bar{S}_1 \cdots \bar{S}_{k+1}) \cdot \boxed{E}^{k+1} = (S_1 \cdots S_{k+1}). \tag{4.7}$$

Assuming Lemma 2, the operator $(L_0 \cdots L_k) \cdot M_{k+1}$ equals $S_{k+1}$ and the identity concerning $S_{k+1}$ in (4.6) is shown too. Next, concerning $S_k$, we obtain, under consideration of definition (3.14) and Lemma 2

$$(L_0 \cdots L_{k-1}) \cdot \begin{pmatrix} M_{2,k+1} \\ \vdots \\ M_{k,k+1} \\ I_B \end{pmatrix} \overset{(3.14)}{=} (L_0 \cdots L_{k-1}) \cdot \boxed{M}^k \cdot \begin{pmatrix} E_{2,k+1} \\ \vdots \\ E_{k,k+1} \\ I_B \end{pmatrix}$$

$$\overset{Lemma\ 2}{=} (S_1 \cdots S_k) \cdot \begin{pmatrix} E_{2,k+1} \\ \vdots \\ E_{k,k+1} \\ I_B \end{pmatrix}. \tag{4.8}$$

Further, by the definition of $E_{i,k+1}, i = k, \ldots, 2$ in (3.12), we see $R[E_{i,k+1}] \subset N_i^c \subset N_{i-1}$ implying

$$\underbrace{S_1 \cdot \overbrace{E_{2,k+1}}^{\to N_2^c \subset N_1}}_{=0} + \cdots + \underbrace{S_{k-1} \cdot \overbrace{E_{k,k+1}}^{\to N_k^c \subset N_{k-1}}}_{=0} + S_k \cdot I_B = S_k \tag{4.9}$$

and the identity concerning $S_k$ in (4.6) is also shown. The remaining identities with respect to $S_{k-1}, \ldots, S_2$ follow in the same way. We only have to employ the equality

$$\boxed{M}_l^{k+1} = \boxed{M}^l \cdot \boxed{E}_l^{l+1} \cdot \ldots \cdot \boxed{E}_l^{k+1} \quad for \quad k \geq 1 \quad and \quad 1 \leq l \leq k \tag{4.10}$$

with $\boxed{M}_l^*$ and $\boxed{E}_l^*$ denoting the $(l \times l)$ triangular matrix composed of last $l$ rows and last $l$ columns of the matrices $\boxed{M}^*$ and $\boxed{E}^*$ respectively. (4.10) follows by direct calculation from the definition of $\boxed{M}^{k+1}$ in (3.14), (3.15).

**Proof of Lemma 2:** By (3.15) we obtain

$$\begin{aligned}(L_0 \cdots L_k) \cdot \boxed{M}^{k+1} &= (L_0 \cdots L_k) \cdot trng[M_1, \ldots, M_{k+1}] \\ &= (L_0 \cdot M_1 \mid \cdots \mid (L_0 \cdots L_k) \cdot M_{k+1}) \end{aligned} \tag{4.11}$$

and the first component of the row vector satisfies Lemma 2 by $L_0 \cdot M_1 = L_0 \cdot I_B = \bar{S}_1 \cdot E_{1,1} = S_1$. Concerning the remaining components in (4.11), we see for $1 \leq i \leq k$, under consideration of definition (3.8),

$$(L_0 \cdots L_i) \cdot M_{i+1} = (L_0 \cdots L_i) \cdot \begin{pmatrix} I_B & | & \\ - & - & - \\ & | & \boxed{M}^i \end{pmatrix} \cdot E_{i+1}$$

$$= (L_0 \cdot I_B \mid L_1 \cdot M_1 \mid \cdots \mid (L_1 \cdots L_i) \cdot M_i) \cdot E_{i+1} \tag{4.12}$$

$$\stackrel{(3.8)}{=} (\bar{S}_1 \mid \bar{S}_2 \mid \cdots \mid \bar{S}_{i+1}) \cdot E_{i+1}$$

and summarizing (4.11), (4.12)

$$\begin{aligned}(L_0 \cdots L_k) \cdot \boxed{M}^{k+1} &= (\bar{S}_1 \cdot E_1 \mid (\bar{S}_1 \mid \bar{S}_2) \cdot E_2 \mid \cdots \mid (\bar{S}_1 \mid \bar{S}_2 \mid \cdots \mid \bar{S}_{k+1}) \cdot E_{k+1}) \\ &= (\bar{S}_1 \cdots \bar{S}_{k+1}) \cdot \boxed{E}^{k+1} \end{aligned} \tag{4.13}$$

yielding the first equality of Lemma 2.

Next, we note that the iteratively defined mappings $E_{i,k+1} = -S_i^{-1}\mathcal{P}_i \cdot \sum_{\nu=i+1}^{k+1} \bar{S}_\nu E_{\nu,k+1}$ in (3.12) of the solution operator $\boxed{E}^{k+1}$ of the triangular system (3.17) may also be written for $k \geq 1$ in an explicit way according to

$$B_i \coloneqq \bar{S}_i\, S_i^{-1}\mathcal{P}_i \in L(\bar{B},\bar{B}),\;\; i = 1,\dots,k \tag{4.14}$$

and

$$\begin{aligned}
E_{1,k+1} &= -S_1^{-1}\,\mathcal{P}_1 \cdot (I_{\bar{B}} - B_2)\cdot \;\dots\; \cdot (I_{\bar{B}} - B_k)\cdot I_{\bar{B}} \cdot \bar{S}_{k+1}\\
&\;\vdots\\
E_{k-1,k+1} &= -S_{k-1}^{-1}\,\mathcal{P}_{k-1} \cdot (I_{\bar{B}} - B_k)\cdot I_{\bar{B}} \cdot \bar{S}_{k+1}\\
E_{k,k+1} &= -S_k^{-1}\,\mathcal{P}_k \cdot I_{\bar{B}} \cdot \bar{S}_{k+1}\\
E_{k+1,k+1} &= I_B\,.
\end{aligned} \tag{4.15}$$

The components $E_{1,k+1},\cdots,E_{k-1,k+1}$ only appear in case of $k \geq 2$. The formulas in (4.15) follow by direct inspection of the bottom up solution process of the triangular system (3.17). Note also that all of the components $E_{1,k+1}$ up to $E_{k,k+1}$ are multiplied by $\bar{S}_{k+1}$ from the right, motivating for $k \geq 2$ the abbreviation

$$E_{i,k+1} = \overbrace{-S_i^{-1}\,\mathcal{P}_i \cdot (I_{\bar{B}} - B_{i+1})\cdot \;\dots\; \cdot (I_{\bar{B}} - B_k)\cdot I_{\bar{B}}}^{=:\, e_{i,k+1}} \cdot \bar{S}_{k+1}\;, \qquad i = 1,\dots,k-1$$

$$E_{k,k+1} = \overbrace{-S_k^{-1}\,\mathcal{P}_k \cdot I_{\bar{B}}}^{=:\, e_{k,k+1}} \cdot \bar{S}_{k+1} \tag{4.16}$$

with $e_{i,k+1} \in L(\bar{B},B), i = 1,\dots,k$. For completeness, we add $e_{1,2} \coloneqq -S_1^{-1}\,\mathcal{P}_1$ and obviously, we obtain for $k \geq 2$ the relation

$$\begin{pmatrix} e_{1,k+1} \\ \vdots \\ e_{k-1,k+1} \end{pmatrix} = \begin{pmatrix} e_{1,k} \\ \vdots \\ e_{k-1,k} \end{pmatrix} \cdot (I_{\bar{B}} - B_k)\,. \tag{4.17}$$

To finish the proof of Lemma 2, we need the following statement concerning $e_{i,k+1},\;\; i = 1,\dots,k$.

**Lemma 3:** For $k \geq 1$

$$(\,\bar{S}_1 \cdots \bar{S}_k\,)\cdot \begin{pmatrix} e_{1,k+1} \\ \vdots \\ e_{k,k+1} \end{pmatrix} = -(\,\mathcal{P}_1 + \cdots + \mathcal{P}_k\,)\,. \tag{4.18}$$

Using Lemma 3, the second equality of Lemma 2 results along the following lines for $k \geq 1$ (for $k = 0$, Lemma 2 is obviously true).

$$(\,\bar{S}_1 \cdots \bar{S}_{k+1}\,)\cdot \boxed{E}^{\,k+1} = (\,\bar{S}_1 \cdots \bar{S}_{k+1}\,)\cdot trng[\,E_1\,,\dots,E_{k+1}\,]$$

$$= (\,\bar{S}_1 \cdot E_1 \;\mid\; (\bar{S}_1\;\bar{S}_2\,)\cdot E_2 \;\mid\; \cdots \;\mid\; (\bar{S}_1 \cdots \bar{S}_{k+1}\,)\cdot E_{k+1}\;) \tag{4.19}$$

$$= \left( \bar{S}_1 \cdot I_B \;\middle|\; \underbrace{\bar{S}_1 \cdot e_{1,2}}_{\underset{(4.18)}{=} -\mathcal{P}_1} \cdot \bar{S}_2 + \bar{S}_2 \cdot I_B \;\middle|\; \cdots \;\middle|\; \underbrace{(\bar{S}_1 \cdots \bar{S}_k) \cdot \begin{pmatrix} e_{1,k+1} \\ \vdots \\ e_{k,k+1} \end{pmatrix}}_{\underset{(4.18)}{=} -\mathcal{P}_1 - \cdots - \mathcal{P}_k} \cdot \bar{S}_{k+1} + \bar{S}_{k+1} \cdot I_B \right)$$

$$= \left( \bar{S}_1 \cdot I_B \;\middle|\; (-\mathcal{P}_1 + I_B) \cdot \bar{S}_2 \;\middle|\; \cdots \;\middle|\; (-\mathcal{P}_1 - \cdots - \mathcal{P}_k + I_B) \cdot \bar{S}_{k+1} \right)$$

$$= \left( \bar{S}_1 \cdot I_B \;\middle|\; P_{R_1^c} \bar{S}_2 \;\middle|\; \cdots \;\middle|\; P_{R_k^c} \bar{S}_{k+1} \right)$$

$$\overset{(3.9)}{=} ( S_1 \; S_2 \; \cdots \; S_{k+1} )$$

**Proof of Lemma 3:** The proof goes by induction for $k \geq 1$. For $k = 1$, Lemma 3 is true by

$$\bar{S}_1 \cdot e_{1,2} = S_1 \cdot (-S_1^{-1} \mathcal{P}_1) = -\mathcal{P}_1 \,. \tag{4.20}$$

Now, assume

$$(\bar{S}_1 \cdots \bar{S}_{k-1}) \cdot \begin{pmatrix} e_{1,k} \\ \vdots \\ e_{k-1,k} \end{pmatrix} = -(\mathcal{P}_1 + \cdots + \mathcal{P}_{k-1}) \tag{4.21}$$

for $k \geq 2$. Then, by (4.16), (4.17) and (4.21)

$$(\bar{S}_1 \cdots \bar{S}_k) \cdot \begin{pmatrix} e_{1,k+1} \\ \vdots \\ e_{k,k+1} \end{pmatrix} = (\bar{S}_1 \cdots \bar{S}_{k-1}) \cdot \begin{pmatrix} e_{1,k+1} \\ \vdots \\ e_{k-1,k+1} \end{pmatrix} + \bar{S}_k \cdot e_{k,k+1}$$

$$\underset{(4.17)}{\overset{(4.16)}{=}} (\bar{S}_1 \cdots \bar{S}_{k-1}) \cdot \begin{pmatrix} e_{1,k} \\ \vdots \\ e_{k-1,k} \end{pmatrix} \cdot (I_{\bar{B}} - B_k) + \bar{S}_k \cdot \left(-S_k^{-1} \mathcal{P}_k\right)$$

$$\overset{(4.21)}{=} -(\mathcal{P}_1 + \cdots + \mathcal{P}_{k-1}) \cdot (I_{\bar{B}} - B_k) - \underbrace{\bar{S}_k S_k^{-1} \mathcal{P}_k}_{= B_k} \tag{4.22}$$

$$= -(\mathcal{P}_1 + \cdots + \mathcal{P}_{k-1}) + (\mathcal{P}_1 + \cdots + \mathcal{P}_{k-1}) \cdot B_k - B_k$$

$$= -(\mathcal{P}_1 + \cdots + \mathcal{P}_{k-1}) + \underbrace{(\mathcal{P}_1 + \cdots + \mathcal{P}_{k-1} - I_{\bar{B}})}_{= -P_{R_{k-1}^c}} \cdot B_k$$

$$= -(\mathcal{P}_1 + \cdots + \mathcal{P}_{k-1}) + \underbrace{P_{R_{k-1}^c} \cdot \bar{S}_k}_{= S_k} S_k^{-1} \mathcal{P}_k$$

$$= -(\mathcal{P}_1 + \cdots + \mathcal{P}_{k-1}) + S_k S_k^{-1} \mathcal{P}_k$$

$$= -(\mathcal{P}_1 + \cdots + \mathcal{P}_{k-1} + \mathcal{P}_k)$$

and the induction is finished. We collect the results concerning the pre-transformation $p_k(\varepsilon)$ and associated triangularization up to the order of $k$ in the following theorem.

**Theorem 1:** Given the iteration (3.3)-(3.15) up to $k \geq 0$ with associated direct sums of the vector spaces $B$ and $\bar{B}$.

$$\begin{array}{ccccccccccccccccc}
B & = & \overbrace{N_0^c & \oplus & N_1^c & \oplus & \cdots & \oplus}^{= N_0^c + \cdots + N_{i-1}^c} & N_i^c & \overbrace{\oplus & \cdots & \oplus & N_{k+1}^c & \oplus & N_{k+1}}^{= N_i} \\
 & & \uparrow & & \uparrow & & & & \uparrow & & & & \uparrow & & \\
 & & \boxed{S_0 := 0} & & \boxed{S_1} & & & & \boxed{S_i} & & & & \boxed{S_{k+1}} & & \\
 & & \downarrow & & \downarrow & & & & \downarrow & & & & \downarrow & & \\
\bar{B} & = & R_0 & \oplus & R_1 & \oplus & \cdots & \oplus & R_i & \underbrace{\oplus & \cdots & \oplus & R_{k+1} & \oplus & R_{k+1}^c}_{= R_{i-1}^c}
\end{array} \tag{4.23}$$

Then, the polynomial of degree $k$ with respect to $\varepsilon$

$$p_k(\varepsilon) = I_B + \varepsilon \cdot M_{k,k+1} + \cdots + \varepsilon^k \cdot M_{1,k+1} = \left( \varepsilon^k \ \cdots \ \varepsilon^1 \ 1 \right) \cdot M_{k+1} \tag{4.24}$$

transforms the power series $L(\varepsilon) = \sum_{i=0}^{\infty} \varepsilon^i \cdot L_i \,,\ L_i \in L(B,\bar{B})$ into

$$S_k(\varepsilon) = L(\varepsilon) \cdot p_k(\varepsilon) = S_1 + \varepsilon^1 \cdot S_2 + \cdots + \varepsilon^k \cdot S_{k+1} + \sum_{i=k+1}^{\infty} \varepsilon^i \cdot Q_{i+1} \tag{4.25}$$

with $k+1$ leading coefficients $S_1, \dots, S_{k+1} \in L(B,\bar{B})$ satisfying

$$S_i \cdot (N_0^c + \cdots + N_{i-1}^c) \subset R_{i-1}^c \qquad \wedge \qquad S_i \cdot N_i^c = R_i \qquad \wedge \qquad S_i \cdot N_i = 0\,, \tag{4.26}$$

i.e. the operators $S_1, \dots, S_{k+1}$ are showing a matrix representation with respect to (4.23) of the form

$S_i =$ 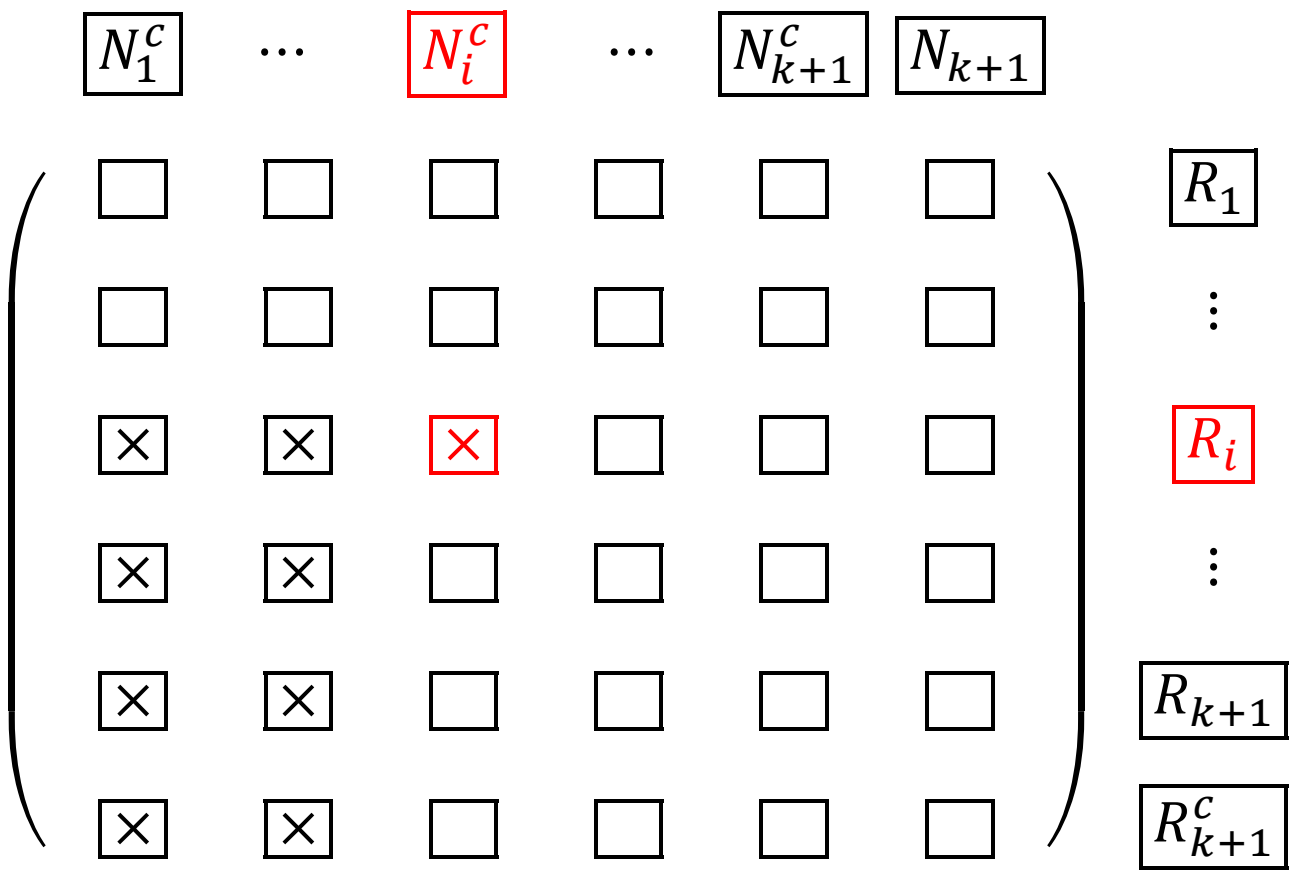 (4.27)

implying lower triangularity of the operator polynomial $P_k(\varepsilon) := S_1 + \varepsilon^1 \cdot S_2 + \cdots + \varepsilon^k \cdot S_{k+1}$ by

$$\begin{array}{cc} & \begin{array}{cccc} \boxed{N_1^c} & \cdots & \boxed{N_{k+1}^c} & \boxed{N_{k+1}} \end{array} \\ P_k(\varepsilon) = & \begin{pmatrix} \boxed{S_1} & \square & \square & \square \\ \boxtimes & \boxed{\ddots} & \square & \square \\ \boxtimes & \boxtimes & \boxed{\varepsilon^k S_{k+1}} & \square \\ \boxtimes & \boxtimes & \square & \square \end{pmatrix} \begin{array}{c} \boxed{R_1} \\ \vdots \\ \boxed{R_{k+1}} \\ \boxed{R_{k+1}^c} \end{array} . \end{array} \tag{4.28}$$

**Remarks: 1**) We have $S_k(\varepsilon) = P_k(\varepsilon) + \sum_{i=k+1}^{\infty} \varepsilon^i \cdot Q_{i+1}$ and up to this point, the remainders $Q_{i+1}$ show no special structure in the sense that in general the corresponding matrix representations are fully occupied, i.e. only partial triangularization is achieved. In the next section we will see that by extending and refining the transformations $p_k(\varepsilon)$ to the power series $\phi(\varepsilon)$, the remainders $Q_{i+1}$ are forced to map into $R_{k+1}^c$ and complete triangularization is obtained. However, this refinement only works, if the Jordan chains (or the leading coefficients) stabilize at some $k \geq 0$. In contrast, closedness of the subspaces within the direct sums of $B$ and $\bar{B}$ is not required for triangularization.

**2)** If $N_{k+1} = \{0\}$ and $R_{k+1}^c = \{0\}$, then a matrix representation of the form

$$\begin{array}{cc} & \begin{array}{ccc} \boxed{N_1^c} & \cdots & \boxed{N_{k+1}^c} \end{array} \\ P_k(\varepsilon) = & \begin{pmatrix} \boxed{S_1} & \square & \square \\ \boxtimes & \boxed{\ddots} & \square \\ \boxtimes & \boxtimes & \boxed{\varepsilon^k S_{k+1}} \end{pmatrix} \begin{array}{c} \boxed{R_1} \\ \vdots \\ \boxed{R_{k+1}} \end{array} \end{array} \tag{4.29}$$

is obtained from (4.28) and $P_k(\varepsilon)$ turns into a bijection for $\varepsilon \neq 0$. This case is treated in [Esquinas] and [López-Gómez] investigating a family of Fredholm operators of index 0. There, the transformations used are polynomials $p_k(\varepsilon)$ of degree $\frac{1}{2}(k + k^2)$ . The transformed family $S_k(\varepsilon)$ is called a transversalization, as soon as bijectivity of $P_k(\varepsilon)$ is achieved for $\varepsilon \neq 0$. In [S1] an alternative approach for transversalization was performed using the polynomial $p_k(\varepsilon)$ of degree $k$ from (4.24). Besides $R_{k+1}^c = \{0\}$, no other assumptions concerning finite dimensionality were needed.

**3)** The partial triangularization $S_k(\varepsilon)$ from (4.25) is characterized by rather simple formulas concerning approximations of order $i = 1, \ldots, k+1$, which are given by

$$b_i(\varepsilon) = n_i + \cdots + \varepsilon^{i-1} \cdot n_1 \qquad \textit{or equivalently} \qquad N[\begin{pmatrix} S_1 & \cdots & S_i \\ & \ddots & \vdots \\ & & S_1 \end{pmatrix}] = \begin{pmatrix} N_1 \\ \vdots \\ N_i \end{pmatrix}. \tag{4.30}$$

In other words, when performing the iteration (3.3)-(3.15) with $L(\varepsilon)$ replaced by $S_k(\varepsilon)$, then the triangular matrix $\boxed{M}^i$ in (3.15) simplifies to a diagonal matrix $\boxed{M}^i = Diag[I_B \cdots I_B]$ and the $i$-tuples $(n_1 \cdots n_i) \in N_1 \times \cdots \times N_i$ in (4.30) define Jordan chains of length $i$ for $n_i \neq 0$.

**4)** By construction, the transformation $p_k(\varepsilon)$ is mapping root elements from $N_{k+1}$ to corresponding approximations of order $k+1$ by

$$p_k(\varepsilon)\cdot n_{k+1} = \left( I_B + \varepsilon \cdot M_{k,k+1} + \cdots + \varepsilon^k \cdot M_{1,k+1} \right) \cdot n_{k+1} \tag{4.31}$$

and we can be sure that $p_k(\varepsilon)$ transforms elements from $N_{k+1}$ in the correct way for reaching (4.25). Now, when choosing $b \notin N_{k+1}$ it is not obvious that $p_k(\varepsilon)$ is working as well. From a geometrical point of view, the possibility to use $p_k(\varepsilon)$ results from the fact that $p_k(\varepsilon)$ is also mapping elements from $N_{i+1}^c, i = 1, \dots, k$ to a corresponding approximation of order $i$, i.e. the steps from $p_i(\varepsilon) \cdot n_{i+1}^c$ up to $p_k(\varepsilon) \cdot n_{i+1}^c$ do not destroy the approximation property of order $i$ with respect to $N_{i+1}^c$.

In some more detail, we have

$$p_k(\varepsilon)\cdot n_{i+1}^c = \left( I_B + \varepsilon \cdot M_{k,k+1} + \cdots + \varepsilon^k \cdot M_{1,k+1} \right) \cdot n_{i+1}^c$$

$$= \left( I_B + \varepsilon \cdot M_{k,k+1} + \cdots + \varepsilon^{i-1} \cdot M_{k-i+2,k+1} + \cdots + \varepsilon^k \cdot M_{1,k+1} \right) \cdot n_{i+1}^c \tag{4.32}$$

$$= \underbrace{\left( I_B + \cdots + \varepsilon^{i-1} \cdot M_{k-i+2,k+1} \right) \cdot n_{i+1}^c}_{Approximation\ of\ order\ i} + \varepsilon^i \cdot \left( M_{k-i+1,k+1} + \cdots + \varepsilon^{k-i} \cdot M_{1,k+1} \right) \cdot n_{i+1}^c$$

with coefficients of $\left( I_B + \cdots + \varepsilon^{i-1} \cdot M_{k-i+2,k+1} \right) \cdot n_{i+1}^c$ satisfying

$$\begin{pmatrix} M_{k-i+2,k+1} \\ \vdots \\ I_B \end{pmatrix} \cdot n_{i+1}^c \in N\left[ \Delta^i \right] = \boxed{M}^{\,i} \cdot \begin{pmatrix} N_1 \\ \vdots \\ N_i \end{pmatrix}, \tag{4.33}$$

i.e. by (4.33) the first summand $\left( I_B + \cdots + \varepsilon^{i-1} \cdot M_{k-i+2,k+1} \right) \cdot n_{i+1}^c$ in (4.32) is an approximation of order $i$ and all in all, we obtain $p_k(\varepsilon) \cdot n_{i+1}^c = \varepsilon^i \cdot r(\varepsilon) \cdot n_{i+1}^c$ with a remainder operator polynomial $r(\varepsilon)$.

To show (4.33) by direct calculation, (4.10) may be used.

## *5. Triangularization and Diagonalization of Operator Power Series*

In this section, the pre-transformation $p_k(\varepsilon)$ is refined and eventually extended to a formal power series $\phi(\varepsilon) = I_B + \sum_{i=1}^{\infty} \varepsilon^i \cdot \phi_i$ in such a way that all remainder terms $Q_{i+1}, i \geq k+1$ in (4.25) turn into operators $S_{i+1}$ of the transformed system $S(\varepsilon) \coloneqq L(\varepsilon) \cdot \phi(\varepsilon) = S_1 + \varepsilon^1 \cdot S_2 + \cdots + \varepsilon^k \cdot S_{k+1} + \sum_{i=k+1}^{\infty} \varepsilon^i \cdot S_{i+1}$ satisfying for $i \geq k+1$

$$R\left[ S_{i+1} \right] \subset R_{k+1}^c \qquad \wedge \qquad S_{i+1} \cdot N_{k+1} = 0 \tag{5.1}$$

with corresponding matrix representation

$$S_{i+1} = \begin{array}{c} \begin{array}{cccc} \boxed{N_1^c} & \cdots & \boxed{N_{k+1}^c} & \boxed{N_{k+1}} \end{array} \\ \begin{pmatrix} \square & \square & \square & \square \\ \square & \square & \square & \square \\ \square & \square & \square & \square \\ \boxtimes & \boxtimes & \boxtimes & \square \end{pmatrix} \end{array} \begin{array}{c} \boxed{R_1} \\ \vdots \\ \boxed{R_{k+1}} \\ \boxed{R_{k+1}^c} \end{array} . \tag{5.2}$$

To achieve (5.1), we have to assume stabilization of the iteration at a certain $k \geq 0$ in the sense that the subspace $\bar{R}_k = R_1 \oplus \cdots \oplus R_{k+1}$ from (3.21) of leading coefficients up to order $k$ cannot be extended anymore according to

$$\underbrace{\bar{R}_{k-1}}_{= R_1 \oplus \cdots \oplus R_k} \subsetneq \underbrace{\bar{R}_k}_{= R_1 \oplus \cdots \oplus R_{k+1}} = \bar{R}_{k+l} \quad for \quad l \geq 1 \quad and \quad \bar{R}_{-1} = \{0\}, \tag{5.3}$$

or equivalently under consideration of (3.20)

$$N_k \supsetneq N_{k+1} = N_{k+l} \qquad for \quad l \geq 1 \tag{5.4}$$

$$\Leftrightarrow \quad R_{k+1} \neq \{0\}, \quad R_{k+l} = \{0\} \quad for \quad l \geq 2 \tag{5.5}$$

$$\Leftrightarrow \quad N_{k+1}^c \neq \{0\}, \quad N_{k+l}^c = \{0\} \quad for \quad l \geq 2 \,. \tag{5.6}$$

Note that by (5.4), root elements in $N_{k+1}$are satisfying $rk(n_{k+1}) = \infty$, i.e. stabilization at $k$ requires the existence of root elements $n_k$ with $rk(n_k) = k$, combined with the nonexistence of root elements with finite rank above $k$.

Note also that, from a geometrical point of view, approximations of arbitrary high order can still be defined via $N[\Delta^{k+l}], l \geq 1$, but when evaluating associated leading coefficients of order $k + l$, the coefficients are not pointing out of $\bar{R}_k = R_1 \oplus \cdots \oplus R_{k+1}$ and hence $\bar{R}_k$ is not increased.

If stabilization occurs at $k \geq 0$, then by (5.3)-(5.6), the following decomposition arises after $l \geq 1$ further steps.

$$\begin{array}{cccccccccccccc}
B & = & N_1^c & \oplus & \cdots & \oplus & N_{k+1}^c & \oplus & \overbrace{\{0\}}^{N_{k+2}^c} & \oplus & \cdots & \oplus & \overbrace{\{0\}}^{N_{k+1+l}^c} & \oplus \; \overbrace{N_{k+1+l}}^{= N_{k+1}} \\
 & & \uparrow & & & & \uparrow & & \uparrow & & & & \uparrow & \\
 & & \boxed{S_1} & & & & \boxed{S_{k+1}} & & \boxed{S_{k+2}} & & & & \boxed{S_{k+1+l}} & \\
 & & \downarrow & & & & \downarrow & & \downarrow & & & & \downarrow & \\
\bar{B} & = & R_1 & \oplus & \cdots & \oplus & R_{k+1} & \oplus & \underbrace{\{0\}}_{R_{k+2}} & \oplus & \cdots & \oplus & \underbrace{\{0\}}_{R_{k+1+l}} & \oplus \; \underbrace{R_{k+1+l}^c}_{= R_{k+1}^c}
\end{array} \tag{5.7}$$

When stabilization does not occur, then the subspace of leading coefficients increases indefinitely in $\bar{B}$ according to

$$\bar{B} \supsetneq R_1 \oplus \cdots \oplus R_{k+1} \oplus \cdots \oplus \underbrace{R_{k+1+l_1}}_{\neq\{0\}} \oplus \cdots \oplus \underbrace{R_{k+1+l_2}}_{\neq\{0\}} \oplus \cdots \tag{5.8}$$

with subspaces $R_{k+1+l_i}$ characterized by arbitrary slow expansion rates $\varepsilon^{k+l_i}, l_i \to \infty$. But then, every kind of generalized inverse family $L^{-1}(\varepsilon)$ will be forced to compensate this slow behaviour by multiplication with $\varepsilon^{-(k+l_i)}$ and some kind of essential singularity will occur at $\varepsilon = 0$. In this sense, the requirement of stabilization at $k$, makes it possible for a pole of order $k$ to arise with respect to $L^{-1}(\varepsilon)$. Concerning essential singularities, see [Albrecht].

Note also that stabilization is ensured, as soon as a first operator $S_{i+1}$ appears in the sequence of operators defined in (3.9) satisfying

$$dim(N_{i+1}) < \infty \qquad or \qquad dim(R_{i+1}^c) < \infty \,, \tag{5.9}$$

i.e. in Banach spaces a semi-Fredholm operator $S_{i+1}$ establishes stabilization at some $k \geq i$. The first possibility for an operator satisfying (5.9) to occur, is given by $S_1 = L_0$ with $dim(N_1) < \infty$ or $dim(R_1^c) < \infty$. In general, finite dimensionality is not necessary for stabilization at $k$.

Theorem 1 from the last section is valid for every $k \geq 0$, but for each $k$ a new transformation $p_k(\varepsilon) = (\,\varepsilon^k \;\cdots\; \varepsilon^1 \;\, 1\,) \cdot M_{k+1}$ has to be build up from the previous transformation $p_{k-1}(\varepsilon) = (\,\varepsilon^{k-1} \;\cdots\; \varepsilon^1 \;\, 1\,) \cdot M_k$. The coefficients of consecutive polynomials $p_0(\varepsilon) = I_B$ , $p_1(\varepsilon) = I_B + \varepsilon \cdot M_{1,2}$ , ... are successively defined by the columns of the $\boxed{M}$ —Matrix, i.e. by $M_1$ , $M_2$ , ... .

Schematically, the $\boxed{M}$ —Matrix shows the following structure.

$$\begin{array}{ccccccccc} M_1 & M_2 & \cdots & \cdots & M_{k+1} & M_{k+2} & M_{k+3} & M_{k+4} & \end{array} \qquad (5.10)$$

$$\boxed{M} = \begin{pmatrix} \boxed{I_B} & M_{1,2} & \cdots & \cdots & M_{1,k+1} & M_{1,k+2} & M_{1,k+3} & M_{1,k+4} & \cdots \\ & \boxed{I_B} & \cdots & \cdots & M_{2,k+1} & M_{2,k+2} & M_{2,k+3} & M_{2,k+4} & \\ & \uparrow & \ddots & & \vdots & \vdots & \vdots & \vdots & \\ & \boxed{p_1(\varepsilon)} & & \boxed{I_B} & M_{k,k+1} & M_{k,k+2} & M_{k,k+3} & M_{k,k+4} & \\ & & & & \boxed{I_B} & M_{k+1,k+2} & M_{k+1,k+3} & M_{k+1,k+4} & \cdots \\ & & & & \uparrow & \boxed{I_B} & M_{k+2,k+3} & M_{k+2,k+4} & \\ & & & & \boxed{p_k(\varepsilon)} & \uparrow & \boxed{I_B} & M_{k+3,k+4} & \\ & & & & & \boxed{p_{k+1}(\varepsilon)} & \uparrow & \boxed{I_B} & \cdots \\ & & & & & & \boxed{p_{k+2}(\varepsilon)} & & \ddots \end{pmatrix}$$

Now, when going from $p_k(\varepsilon)$ to $p_{k+1}(\varepsilon)$ or from $p_{k+1}(\varepsilon)$ to $p_{k+2}(\varepsilon)$, then in general coefficients of common order $\geq 1$ of the polynomials (lying on the same diagonal above the main diagonal) differ completely, as indicated in (5.11).

$$\begin{array}{ccccccc} \varepsilon^k & \rightarrow & M_{1,k+1} & \neq & M_{2,k+2} & \neq & M_{3,k+3} \\ \varepsilon^{k-1} & \rightarrow & M_{2,k+1} & \neq & M_{3,k+2} & \neq & M_{4,k+3} \\ \vdots & \vdots & \vdots & \vdots & \vdots & \vdots & \vdots \\ \varepsilon^1 & \rightarrow & M_{k,k+1} & \neq & M_{k+1,k+2} & \neq & M_{k+2,k+3} \\ & & \uparrow & & \uparrow & & \uparrow \\ & & \boxed{p_k(\varepsilon)} & & \boxed{p_{k+1}(\varepsilon)} & & \boxed{p_{k+2}(\varepsilon)} \end{array} \qquad (5.11)$$

On contrary, we will see below, if stabilization occurs at $k$, then the first order coefficients of $p_{k+1}(\varepsilon)$ and $p_{k+2}(\varepsilon)$ coincide by $M_{k+1,k+2} = M_{k+2,k+3}$, whereas higher order coefficients may remain different. The fact that within matrix $\boxed{M}$ , the element $M_{k+2,k+3}$ is simply obtained from the transfer of the element top left of $M_{k+2,k+3}$, i.e. from $M_{k+1,k+2}$, is depicted in (5.10) by a small arrow according to $M^{\searrow}_{k+2,k+3}$.

Alternatively, we can say that with blue marked element $M_{k+1,k+2}$ of $p_{k+1}(\varepsilon)$, the coefficient of $\varepsilon$ has already reached its final configuration for all subsequent polynomials $p_{k+2}(\varepsilon)$, $p_{k+3}(\varepsilon)$, … . This behaviour continues in the sense that with $M_{k+1,k+3}$, the coefficient of $\varepsilon^2$ attains its final configuration for the polynomials $p_{k+2}(\varepsilon)$, $p_{k+3}(\varepsilon)$, … .

In summary, if stabilization occurs at $k \geq 0$, then the blue marked row $k+1$ of matrix $\boxed{M}$ establishes some sort of final configuration of coefficients up to infinity and, in particular, all rows below row $k+1$ are simply given by copies of row $k+1$ shifted to the right and down, i.e. starting with blue row $k+1$ an operator Töplitz matrix arises. Note also that in case of $k=0$, row $k+1$ turns into the first row and the complete matrix $\boxed{M}$ represents a Töplitz matrix.

Now, as expected, blue row $k+1$ delivers the formal power series $\phi(\varepsilon)$ looked for according to

$$\phi(\varepsilon) := I_B + \sum_{i=1}^{\infty} \varepsilon^i \cdot \overbrace{M_{k+1,k+1+i}}^{=:\,\phi_i} \tag{5.12}$$

and we see that $\phi(\varepsilon)$ corresponds to the Laurent power series that is usually associated to a Töplitz matrix with infinite dimensions, as introduced in [Töplitz]. Due to our assumption $L(\varepsilon) = \sum_{i=0}^{\infty} \varepsilon^i \cdot L_i$, the principal part of the Laurent power series vanishes in (5.12). In case of meromorphic operator families $L(\varepsilon)$, the power series of $\phi(\varepsilon)$ in (5.12) would contain negative exponents as well.

Note that the pre-transformations $p_i(\varepsilon)$ are defined successively by the columns of the matrix $\boxed{M}$, whereas the transformation $\phi(\varepsilon)$ is defined by row $k+1$. Note also that for calculation of the coefficients $\phi_1$, $\phi_2$, … , the last pre-transformation $p_k(\varepsilon)$ before stabilization is not used, but instead polynomials of higher degree $p_{k+1}(\varepsilon)$, $p_{k+2}(\varepsilon)$, … have to be applied.

Schematically, we end up with the following matrix.

$$\boxed{M} = \begin{pmatrix} \boxed{I_B} & M_{1,2} & \cdots & M_{1,k+1} & M_{1,k+2} & \cdots & M_{1,k+1+l} & \cdots \\ & \ddots & & \vdots & \vdots & \vdots & \vdots & \\ & & \boxed{I_B} & M_{k,k+1} & M_{k,k+2} & \cdots & M_{k,k+1+l} & \cdots \\ & \boxed{\phi(\varepsilon)} & \longrightarrow & \boxed{I_B} & \boxed{\phi_1} & \cdots & \boxed{\phi_l} & \cdots \\ & & & \uparrow & \ddots & & \vdots & \\ & & & \boxed{p_k(\varepsilon)} & \uparrow & \ddots & \boxed{\phi_1} & \\ & & & & \boxed{p_{k+1}(\varepsilon)} & & \boxed{I_B} & \cdots \\ & & & & & & \uparrow & \ddots \\ & & & & & & \boxed{p_{k+l}(\varepsilon)} & \end{pmatrix} \tag{5.13}$$

Note that, due to the Töplitz structure, the transformations $p_{k+l}(\varepsilon)$ agree for $l \geq 1$ with $\phi(\varepsilon)$ in first $l+1$ coefficients according to

$$\begin{pmatrix} \phi_l \\ \vdots \\ \phi_1 \\ I_B \end{pmatrix} = \begin{pmatrix} M_{k+1,k+1+l} \\ \vdots \\ M_{k+l,k+1+l} \\ I_B \end{pmatrix} . \tag{5.14}$$

Now, concerning the transformed system, we obtain by Cauchy product

$$\begin{aligned} L(\varepsilon) \cdot \phi(\varepsilon) &= \sum_{l=0}^{\infty} \varepsilon^l \cdot \sum_{i+j=l} L_i \phi_j = \sum_{l=0}^{\infty} \varepsilon^l \cdot ( L_0 \cdots L_l ) \cdot \begin{pmatrix} \phi_l \\ \vdots \\ I_B \end{pmatrix} \\ &= \sum_{l=0}^{\infty} \varepsilon^l \cdot ( L_0 \cdots L_l ) \cdot \begin{pmatrix} M_{k+1,k+1+l} \\ \vdots \\ M_{k+l,k+1+l} \\ I_B \end{pmatrix} \\ &= \sum_{l=0}^{\infty} \varepsilon^l \cdot S_{l+1} \end{aligned} \tag{5.15}$$

with the last identity following from Theorem 1 according to

$$\begin{aligned} L(\varepsilon) \cdot p_{k+l}(\varepsilon) = ( L_0 + \cdots + \varepsilon^l \cdot L_l + \sum_{i=l+1}^{\infty} \varepsilon^i \cdot L_i ) \\ \cdot ( I_B + \varepsilon \cdot M_{k+l,k+1+l} + \cdots + \varepsilon^l \cdot M_{k+1,k+1+l} + \cdots + \varepsilon^{k+l} \cdot M_{1,k+1+l} ) \\ \overset{Th.1}{=} S_1 + \cdots + \varepsilon^l \cdot S_{l+1} + \cdots + \varepsilon^{k+l} \cdot S_{k+1+l} + \sum_{i=k+1+l}^{\infty} \varepsilon^i \cdot Q_{i+1} \end{aligned} \tag{5.16}$$

and $S_{l+1}$ given by

$$S_{l+1} = ( L_0 \cdots L_l ) \cdot \begin{pmatrix} M_{k+1,k+1+l} \\ \vdots \\ M_{k+l,k+1+l} \\ I_B \end{pmatrix} \tag{5.17}$$

for $l \geq 1$. Thus, by (5.15) and (4.7) in Lemma 2, the coefficients of the transformed power series $L(\varepsilon) \cdot \phi(\varepsilon)$ are reproducing the mappings $S_i$ , $i \geq 1$, defined recursively by (3.3)-(3.15), and it remains to show that $S_{l+1}$ satisfies the conditions (5.1) for $l \geq k+1$.

Now, by (4.26) of Theorem 1 and choosing $l \geq k+1$, we see

$$S_{l+1} \cdot ( N_0^c + \cdots + N_l^c ) \subset R_l^c \overset{(5.7)}{=} R_{k+1}^c \quad \wedge \quad S_{l+1} \cdot N_{l+1} \overset{(5.7)}{=} S_{l+1} \cdot N_{k+1} = 0 \tag{5.18}$$

implying (5.1), as desired.

As already mentioned, (5.13) shows in detail that for calculation of the coefficients $\phi_1, \dots, \phi_k$, we cannot use the basic pre-transformation $p_k(\varepsilon)$, but we have to go up to the refined coefficients of the polynomial $p_{2k}(\varepsilon)$ of degree $2k$ (set $l = k$ in (5.13)). In this sense, the step from $k$ to $2k$ seems to be crucial with respect to the behaviour of $L(\varepsilon)$ under stabilization with $k$.

It remains to show that pattern (5.10) of $\boxed{M}$, in fact turns into the Töplitz pattern (5.13) in case of stabilization at $k$. Remember, matrix $\boxed{M}$ is defined column by column using (3.14), (3.15) as well as evaluating matrix $\boxed{E}$ of solution operators of the triangular system (3.17).

Now, in case of stabilization at $k$, matrix $\boxed{E}$ adopts the form

$$
\begin{array}{cccccccccc}
E_1 & E_2 & \cdots & E_{k+1} & E_{k+2} & E_{k+3} & \cdots & E_{k+1+l} & \\
\hline
\end{array}
\tag{5.19}
$$

$$
\boxed{E} = \begin{pmatrix}
\boxed{I_B} & E_{1,2} & \cdots & E_{1,k+1} & E_{1,k+2} & E_{1,k+3} & \cdots & E_{1,k+1+l} & \cdots \\
 & \ddots & & \vdots & \vdots & \vdots & \vdots & \vdots & \\
 & & \boxed{I_B} & E_{k,k+1} & E_{k,k+2} & E_{k,k+3} & \cdots & E_{k,k+1+l} & \cdots \\
 & & & \boxed{I_B} & \boxed{E_{k+1,k+2}} & \boxed{E_{k+1,k+3}} & \cdots & \boxed{E_{k+1,k+1+l}} & \cdots \\
 & & & & \boxed{I_B} & \boxed{0} & \cdots & \boxed{0} & \cdots \\
 & & & & & \boxed{I_B} & \ddots & \vdots & \\
 & & & & & & \ddots & \boxed{0} & \cdots \\
 & & & & & & & \boxed{I_B} & \ddots \\
 & & & & & & & & \ddots
\end{pmatrix}
$$

with green marked zero operators following from (5.5), due to $R_{k+l} = \{0\}$ for $l \geq 2$ and

$$
\begin{aligned}
E_{k+2,i} &= -S_{k+2}^{-1} \overbrace{\mathcal{P}_{k+2}}^{=0} \cdot [\cdots] \cdot \bar{S}_{k+i} = 0, \qquad i \geq k+3 \\
E_{k+3,i} &= -S_{k+3}^{-1} \overbrace{\mathcal{P}_{k+3}}^{=0} \cdot [\cdots] \cdot \bar{S}_{k+i} = 0, \qquad i \geq k+4 \\
&\ \ \vdots \qquad\qquad\qquad\qquad\qquad \vdots \qquad\qquad \vdots
\end{aligned}
\tag{5.20}
$$

by (4.15). Then, from the definition of the columns of $\boxed{M}$ in (3.14), we obtain for $l \geq 2$

$$
M_{k+1+l} = \begin{pmatrix} I_B & | & \\ - & - & - - - \\ & | & \boxed{M}^{k+l} \end{pmatrix} \cdot E_{k+1+l}
\tag{5.21}
$$

and last $l$ components of $M_{k+1+l}$ read

$$\begin{pmatrix} M^{\searrow}_{k+2,k+1+l} \\ \vdots \\ M^{\searrow}_{k+l,k+1+l} \\ M^{\searrow}_{k+1+l,k+1+l} \end{pmatrix} = \begin{pmatrix} \boxed{I_B} & M_{k+1,k+2} & \cdots & M_{k+1,k+l} \\ & \ddots & & \vdots \\ & & \ddots & M_{k+l-1,k+l} \\ & & & \boxed{I_B} \end{pmatrix} \cdot \begin{pmatrix} E_{k+2,k+1+l} \\ \vdots \\ E_{k+l,k+1+l} \\ E_{k+1+l,k+1+l} \end{pmatrix}$$

$$\overset{(5.19)}{=} \begin{pmatrix} \boxed{I_B} & M_{k+1,k+2} & \cdots & M_{k+1,k+l} \\ & \ddots & & \vdots \\ & & \ddots & M_{k+l-1,k+l} \\ & & & \boxed{I_B} \end{pmatrix} \cdot \begin{pmatrix} 0 \\ \vdots \\ 0 \\ \boxed{I_B} \end{pmatrix} = \begin{pmatrix} M_{k+1,k+l} \\ \vdots \\ M_{k+l-1,k+l} \\ \boxed{I_B} \end{pmatrix} . \tag{5.22}$$

Hence, last $l$ components of $M_{k+1+l}$ in (5.10) are simply obtained by shifting last $l$ components of $M_{k+l}$ to the right and down and the repetition structure of $\boxed{M}$ in (5.10), (5.13) is shown.

Concerning stabilization, triangularization and diagonalization of operator power series $L(\varepsilon)$ between real or complex vector spaces, we summarize the results in the following theorem.

**Theorem 2:** Assume stabilization of the iteration (3.3)-(3.15) at $k \geq 0$ and define a power series $\phi(\varepsilon) \in L(B, B)$ by row $k+1$ of matrix $\boxed{M}$ according to

$$\phi(\varepsilon) = \sum_{i=0}^{\infty} \varepsilon^i \cdot M_{k+1,k+1+i} \, . \tag{5.23}$$

(i) Then, $\phi(\varepsilon)$ transforms $L(\varepsilon) = \sum_{i=0}^{\infty} \varepsilon^i \cdot L_i$ by Cauchy product into the triangularization

$$S(\varepsilon) := L(\varepsilon) \cdot \phi(\varepsilon) = S_1 + \varepsilon^1 \cdot S_2 + \cdots + \varepsilon^k \cdot S_{k+1} + \sum_{l=2}^{\infty} \varepsilon^{k+l-1} \cdot S_{k+l} \tag{5.24}$$

with $S_1, \ldots, S_{k+1}$ satisfying (4.26) from Theorem 1 and $S_{k+l}$ satisfying for $l \geq 2$

$$R[\, S_{k+l} \,] \subset R^c_{k+1} \qquad \wedge \qquad S_{k+l} \cdot N_{k+1} = 0 \tag{5.25}$$

with associated matrix representation

$$S_{k+l} = \begin{array}{c} \begin{matrix} \boxed{N_1^c} & \cdots & \boxed{N_{k+1}^c} & \boxed{N_{k+1}} \end{matrix} \\ \begin{pmatrix} \square & \square & \square & \square \\ \square & \square & \square & \square \\ \square & \square & \square & \square \\ \boxtimes & \boxtimes & \boxtimes & \square \end{pmatrix} \end{array} \begin{matrix} \boxed{R_1} \\ \vdots \\ \boxed{R_{k+1}} \\ \boxed{R^c_{k+1}} \end{matrix} . \tag{5.26}$$

In particular, for all $i \geq 0$, we have $S_{i+1} \cdot N_{k+1} = 0$.

(ii) For every power series $b(\varepsilon) = \sum_{i=0}^{\infty} \varepsilon^i \cdot b_i \,, b_i \in B$ with

$$L(\varepsilon) \cdot b(\varepsilon) = 0\,, \tag{5.27}$$

there exists a unique power series $n_{k+1}(\varepsilon) = \sum_{i=0}^{\infty} \varepsilon^i \cdot n_{k+1}^i \,, n_{k+1}^i \in N_{k+1}$ with

$$b(\varepsilon) = \phi(\varepsilon) \cdot n_{k+1}(\varepsilon) \tag{5.28}$$

and coefficients recursively determined by

$$\begin{aligned} n_{k+1}^0 &= b_0 \\ n_{k+1}^i &= b_i - (\, \phi_1 \ \cdots \ \phi_i \,) \cdot \begin{pmatrix} n_{k+1}^{i-1} \\ \vdots \\ n_{k+1}^0 \end{pmatrix}, \ \ i \geq 1\,. \end{aligned} \tag{5.29}$$

(iii) The set of all power series solutions of $L(\varepsilon) \cdot b(\varepsilon) = 0$ is given by

$$X_\infty := \{\, (b_i)_{i \in \mathbb{N}_0} \ | \ b_i = (I_B \ \cdots \ \phi_i) \cdot \begin{pmatrix} n_{k+1}^i \\ \vdots \\ n_{k+1}^0 \end{pmatrix}, \ \ n_{k+1}^i \in N_{k+1} \,\}\,. \tag{5.30}$$

(iv) The power series $S(\varepsilon) = \sum_{i=0}^{\infty} \varepsilon^i \cdot S_{i+1}\,, \ S_{i+1} \in L(B, \bar{B})$ from (i) can be factorized by a power series $\psi(\varepsilon)$ and a diagonal operator polynomial $\Delta(\varepsilon)$ according to

$$\begin{aligned} S(\varepsilon) &= \psi(\varepsilon) \cdot \Delta(\varepsilon) \\ &:= [\, I_{\bar{B}} \ + \sum_{i=1}^{\infty} \varepsilon^i \cdot \psi_i \,] \cdot [\, S_1 P_1 + \varepsilon \cdot S_2 P_2 + \cdots + \varepsilon^k \cdot S_{k+1} P_{k+1} \,] \end{aligned} \tag{5.31}$$

with $\psi_i \in L(\bar{B}, \bar{B})$ and $P_1, \dots, P_{k+1}$ denoting projections to $N_1^c, \dots, N_{k+1}^c$ respectively.

Concerning the proof of (ii), we first obtain from $L(\varepsilon) \cdot b(\varepsilon) = 0$ for $k \geq 0$ and $l \geq 1$ an approximation of order $k + l$ according to

$$L(\varepsilon) \cdot (\, b_0 + \cdots + \varepsilon^{k+1-l} \cdot b_{k+1-l} \,) = \varepsilon^{k+l} \cdot \bar{b}(\varepsilon) \tag{5.32}$$

with remainder power series $\bar{b}(\varepsilon)$ in $\bar{B}$, and from Lemma 1 (i) we see

$$\begin{pmatrix} b_{k+1-l} \\ \vdots \\ b_0 \end{pmatrix} \in N[\, \Delta^{k+l} \,] = \boxed{M}^{\,k+l} \cdot \begin{pmatrix} N_1 \\ \vdots \\ N_{k+l} \end{pmatrix}. \tag{5.33}$$

Then, starting at $l = 1$, there exist $(n_1, \dots, n_{k+1})$ with

$$\boxed{M}^{\,k+1} \cdot \begin{pmatrix} n_1 \\ \vdots \\ n_{k+1} \end{pmatrix} = \begin{pmatrix} b_k \\ \vdots \\ b_0 \end{pmatrix} \tag{5.34}$$

and the last component yields $I_B \cdot n_{k+1} = b_0$. Next, with $l = 2$, there exist $(n_1, \dots, n_{k+1}, n_{k+2})$ satisfying

$$\boxed{M}^{k+2} \cdot \begin{pmatrix} n_1 \\ \vdots \\ n_{k+1} \\ n_{k+2} \end{pmatrix} = \begin{pmatrix} b_{k+1} \\ \vdots \\ b_1 \\ b_0 \end{pmatrix} \tag{5.35}$$

and, due to stabilization at $k$, we have $N_{k+2} = N_{k+1}$ and the last two components in (5.35) read

$$\begin{pmatrix} I_B & \phi_1 \\ & I_B \end{pmatrix} \cdot \begin{pmatrix} n_{k+1}^1 \\ n_{k+1}^0 \end{pmatrix} \overset{(5.24)}{=} \begin{pmatrix} b_1 \\ b_0 \end{pmatrix}, \tag{5.36}$$

where we used an upper index to distinguish the two elements from $N_{k+1}$. We reproduce $n_{k+1}^0 = b_0$ from the second equation in (5.36), whereas the first equation has a unique solution given by

$$n_{k+1}^1 = b_1 - \phi_1 \cdot n_{k+1}^0 \tag{5.37}$$

and (5.29) is shown for $i = 1$. Now, the general induction step from $i$ to $i + 1$ works along the same lines of reasoning.

Concerning (iii), we see from (5.29) that every power series solution is contained in $X_\infty$. Reversly, choosing an element $(b_i)_{i\in\mathbb{N}_0}$ from $X_\infty$ we obtain

$$\begin{aligned} L(\varepsilon) \cdot \sum_{i=0}^{\infty} \varepsilon^i \cdot b_i &= \sum_{l=0}^{\infty} \varepsilon^l \cdot (L_0 \ \cdots \ L_l) \cdot \begin{pmatrix} b_l \\ \vdots \\ b_0 \end{pmatrix} \\ &= \sum_{l=0}^{\infty} \varepsilon^l \cdot (L_0 \ \cdots \ L_l) \cdot \begin{pmatrix} I_B & \cdots & \phi_l \\ & \ddots & \vdots \\ & & I_B \end{pmatrix} \cdot \begin{pmatrix} n_{k+1}^l \\ \vdots \\ n_{k+1}^0 \end{pmatrix} \end{aligned} \tag{5.38}$$

with the blue matrix indicating last $l$ columns and rows of the matrix $\boxed{M}^{k+1+l}$ according to

$$\boxed{M}^{k+1+l} \overset{(5.13)}{=} \begin{pmatrix} \boxed{I_B} & M_{1,2} & \cdots & M_{1,k+1} & M_{1,k+2} & \cdots & M_{1,k+1+l} \\ & \ddots & & \vdots & \vdots & \vdots & \vdots \\ & & \boxed{I_B} & M_{k,k+1} & M_{k,k+2} & \cdots & M_{k,k+1+l} \\ & & & \boxed{I_B} & \boxed{\phi_1} & \cdots & \boxed{\phi_l} \\ & & & & \ddots & & \vdots \\ & & & & & \ddots & \boxed{\phi_1} \\ & & & & & & \boxed{I_B} \end{pmatrix}. \tag{5.39}$$

Then, from Lemma 1 we conclude

$$N[\,\Delta^{k+1+l}\,] = N[\begin{pmatrix} L_0 & \cdots & \cdots & \cdots & L_{k+l} \\ & \ddots & & & \vdots \\ & & L_0 & \cdots & L_l \\ & & & \ddots & \vdots \\ & & & & L_0 \end{pmatrix}] = R[\,M^{k+1+l}{}_{|\,N_1\times\cdots\times N_k\times N_{k+1}^l\cdots N_{k+1}^0}\,] \tag{5.40}$$

implying all coefficients in (5.38) to be zero by direct inspection.

Concerning (iv), define

$$\psi(\varepsilon) := \underbrace{I_{\bar{B}}}_{=:\psi_0} + \sum_{i=1}^{\infty}\varepsilon^i\cdot\psi_i\,,\qquad \psi_i := S_{i+1}\cdot S_1^{-1}\mathcal{P}_1 + \cdots + S_{k+i+1}\cdot S_{k+1}^{-1}\mathcal{P}_{k+1} \tag{5.41}$$

and remember $P_1,\dots,P_{k+1}$ as well as $\mathcal{P}_1,\dots,\mathcal{P}_{k+1}$ denoting projections with respect to the subspaces $N_1^c,\dots,N_{k+1}^c$ and $R_1,\dots,R_{k+1}$ within the direct sum decompositions

$$\begin{aligned} B &= \overbrace{N_1^c}^{\downarrow P_1} \oplus \cdots \oplus \overbrace{N_{k+1}^c}^{\downarrow P_{k+1}} \oplus N_{k+1} \\ \bar{B} &= \underbrace{R_1}_{\uparrow \mathcal{P}_1} \oplus \cdots \oplus \underbrace{R_{k+1}}_{\uparrow \mathcal{P}_{k+1}} \oplus R_{k+1}^c \end{aligned} \tag{5.42}$$

from (5.7). Then, rewrite $\psi(\varepsilon)$ by

$$\begin{aligned} \psi(\varepsilon) &= I_{\bar{B}} + \sum_{i=1}^{\infty}\varepsilon^i\cdot\left(S_{i+1}\cdot S_1^{-1}\mathcal{P}_1 + \cdots + S_{k+i+1}\cdot S_{k+1}^{-1}\mathcal{P}_{k+1}\right) \\ &= I_{\bar{B}} + \sum_{i=1}^{\infty}\varepsilon^i\cdot S_{i+1}\cdot S_1^{-1}\mathcal{P}_1 + \cdots + \sum_{i=1}^{\infty}\varepsilon^i\cdot S_{k+i+1}\cdot S_{k+1}^{-1}\mathcal{P}_{k+1} \\ &= I_{\bar{B}} + \varepsilon\cdot\left(\sum_{i=1}^{\infty}\varepsilon^{i-1}\cdot S_{i+1}\right)\cdot S_1^{-1}\mathcal{P}_1 + \cdots + \varepsilon\cdot\left(\sum_{i=1}^{\infty}\varepsilon^{i-1}\cdot S_{k+i+1}\right)\cdot S_{k+1}^{-1}\mathcal{P}_{k+1} \\ &=: I_{\bar{B}} + \varepsilon\cdot r_1^c(\varepsilon)\cdot S_1^{-1}\mathcal{P}_1 + \cdots + \varepsilon\cdot r_{k+1}^c(\varepsilon)\cdot S_{k+1}^{-1}\mathcal{P}_{k+1} \end{aligned} \tag{5.43}$$

with power series $r_i^c(\varepsilon), i = 1,\dots,k+1$. Now, the factorization $S(\varepsilon) = \psi(\varepsilon)\cdot\Delta(\varepsilon)$ can be simply shown by Cauchy product under consideration of splittings (5.42). We obtain

$$\begin{aligned} &\psi(\varepsilon)\cdot\Delta(\varepsilon) \\ &= \left[\,I_{\bar{B}} + \varepsilon\cdot r_1^c(\varepsilon)\cdot S_1^{-1}\mathcal{P}_1 + \cdots + \varepsilon\cdot r_{k+1}^c(\varepsilon)\cdot S_{k+1}^{-1}\mathcal{P}_{k+1}\,\right]\cdot\left[\,S_1P_1 + \cdots + \varepsilon^k\cdot S_{k+1}P_{k+1}\,\right] \\ &= I_{\bar{B}}\cdot\left[\,S_1P_1 + \cdots + \varepsilon^k\cdot S_{k+1}P_{k+1}\,\right] + \varepsilon\cdot r_1^c(\varepsilon)\cdot P_1 + \cdots + \varepsilon^{k+1}\cdot r_{k+1}^c(\varepsilon)\cdot P_{k+1} \end{aligned}$$

$$= [\, S_1 + \varepsilon \cdot r_1^c(\varepsilon)\,] \cdot P_1 \; + \cdots + \; [\, \varepsilon^k \cdot S_{k+1} + \varepsilon^{k+1} \cdot r_{k+1}^c(\varepsilon)\,] \cdot P_{k+1} \tag{5.44}$$

$$= [\, S(\varepsilon)\,] \cdot P_1 \; + \cdots + \; \left[ S(\varepsilon) - \left( S_1 + \cdots + \varepsilon^{k-1} \cdot S_k \right) \right] \cdot P_{k+1}$$

$$= S(\varepsilon) \cdot P_1 + \cdots + \; S(\varepsilon) \cdot P_{k+1}$$

$$= S(\varepsilon) \cdot \underbrace{\left( P_1 + \cdots + P_{k+1} + P_{N_{k+1}} \right)}_{= I_B} \; ,$$

due to $S(\varepsilon) \cdot P_{N_{k+1}} = 0$ by Theorem 2 (i) and $P_{N_{k+1}}$ denoting projection to $N_{k+1}$ within (5.42). This finishes the proof of Theorem 2.

## *6. Defining Equation for the Transformation*

In this section, a defining equation for $\phi(\varepsilon) = \sum_{i=0}^{\infty} \varepsilon^i \cdot \phi_i$ is derived, finally allowing to show analyticity of $\phi(\varepsilon)$, if stabilization at $k \geq 0$, analyticity of the operator family $L(\varepsilon) = \sum_{i=0}^{\infty} \varepsilon^i \cdot L_i$ and continuity of the projections is assumed.

Up to now, stabilization at $k$ implies the blue marked Töplitz repetition structure of $\boxed{M}$ in (5.13), as well as the green marked zero operators of $\boxed{E}$ in (5.19). These patterns are sufficient to prove triangularization and diagonalization as stated in Theorem 2. Yet, more structure is present in $\boxed{M}$ and $\boxed{E}$.

First, note that in case of stabilization at $k \geq 0$, we have $R_i = \{0\}$ for $i \geq k+2$ and

$$I_{\bar{B}} - B_i \overset{(4.14)}{=} I_{\bar{B}} - \bar{S}_i \, S_i^{-1} \underbrace{\mathcal{P}_i}_{=0} \; = \; I_{\bar{B}} \; , \tag{6.1}$$

yielding by (4.16) for $l \geq 2$

$$\begin{pmatrix} e_{1,k+1+l} \\ \vdots \\ e_{k+1,k+1+l} \end{pmatrix} = \begin{pmatrix} e_{1,k+2} \\ \vdots \\ e_{k+1,k+2} \end{pmatrix} \cdot \underbrace{(I_{\bar{B}} - B_{k+2})}_{= I_{\bar{B}}} \cdot \cdots \cdot \underbrace{(I_{\bar{B}} - B_{k+l})}_{= I_{\bar{B}}} = \begin{pmatrix} e_{1,k+2} \\ \vdots \\ e_{k+1,k+2} \end{pmatrix} . \tag{6.2}$$

Hence, first $k+1$ components of the columns $E_{k+2}, E_{k+3}, \ldots$ are almost equal in the sense that they differ only by multiplication with $\bar{S}_{k+2}, \bar{S}_{k+3}, \ldots$ according to

$$\begin{pmatrix} E_{1,k+2} \\ \vdots \\ E_{k+1,k+2} \end{pmatrix} = \begin{pmatrix} e_{1,k+2} \\ \vdots \\ e_{k+1,k+2} \end{pmatrix} \cdot \bar{S}_{k+2} \, , \quad \begin{pmatrix} E_{1,k+3} \\ \vdots \\ E_{k+1,k+3} \end{pmatrix} = \begin{pmatrix} e_{1,k+2} \\ \vdots \\ e_{k+1,k+2} \end{pmatrix} \cdot \bar{S}_{k+3} \, , \quad \cdots \, . \tag{6.3}$$

We indicate these elements by brown colour within matrix $\boxed{E}$.

$$
\begin{array}{c}
\begin{array}{cccccccc}
E_1 & E_2 & \cdots & E_{k+1} & E_{k+2} & E_{k+3} & \cdots & E_{k+1+l}
\end{array} \\
\hline
\end{array}
\tag{6.4}
$$

$$
\boxed{E} = \begin{pmatrix}
\boxed{I_B} & E_{1,2} & \cdots & E_{1,k+1} & E_{1,k+2} & E_{1,k+3} & \cdots & E_{1,k+1+l} & \cdots \\
 & \ddots & & \vdots & \vdots & \vdots & \vdots & \vdots & \vdots \\
 & & \boxed{I_B} & E_{k,k+1} & E_{k,k+2} & E_{k,k+3} & \cdots & E_{k,k+1+l} & \cdots \\
 & & & \boxed{I_B} & \boxed{E_{k+1,k+2}} & \boxed{E_{k+1,k+3}} & \cdots & \boxed{E_{k+1,k+1+l}} & \cdots \\
 & & & & \boxed{I_B} & \boxed{0} & \cdots & \boxed{0} & \cdots \\
 & & & & & \boxed{I_B} & \ddots & \vdots & \cdots \\
 & & & & & & \ddots & \boxed{0} & \cdots \\
 & & & & & & & \boxed{I_B} & \cdots
\end{pmatrix}
$$

Now, green zeros and brown elements have the following consequences upon constructing the columns of the matrix $\boxed{M}$. For $k \geq 0$ and $l \geq 1$, we have

$$
M_{k+1+l} = \begin{pmatrix} I_B & | \\ - & - - - - \\ & | \ \boxed{M} \end{pmatrix}^{k+l} \cdot E_{k+1+l} = \begin{pmatrix} I_B & | \\ - & - - - - \\ & | \ \boxed{M} \end{pmatrix}^{k+l} \cdot \begin{pmatrix} E_{1,k+1+l} \\ \vdots \\ E_{k+1,k+1+l} \\ \boxed{0} \\ \vdots \\ \boxed{0} \\ \boxed{I_B} \end{pmatrix}
\tag{6.5}
$$

$$
= \begin{pmatrix} \begin{pmatrix} I_B & | \\ - & - - - \\ & | \ \boxed{M} \end{pmatrix}^{k} \cdot \begin{pmatrix} E_{1,k+1+l} \\ \vdots \\ E_{k+1,k+1+l} \end{pmatrix} \\ \hline 0 \\ \vdots \\ 0 \end{pmatrix} + \begin{pmatrix} 0 \\ M_{k+l} \end{pmatrix} \cdot \boxed{I_B}
$$

$$
\overset{(6.3)}{=} \begin{pmatrix} \begin{pmatrix} I_B & | \\ - & - - - \\ & | \ \boxed{M} \end{pmatrix}^{k} \cdot \begin{pmatrix} e_{1,k+2} \\ \vdots \\ e_{k+1,k+2} \end{pmatrix} \cdot \bar{S}_{k+1+l} \\ \hline 0 \\ \vdots \\ 0 \end{pmatrix} + \begin{pmatrix} 0 \\ M_{k+l} \end{pmatrix},
$$

which motivates the abbreviation

$$\bar{H} = \begin{pmatrix} H_1 \\ \vdots \\ H_{k+1} \end{pmatrix} := \begin{pmatrix} I_B & | & \\ - & - & - \; - \\ & | & \boxed{M}^k \end{pmatrix} \cdot \begin{pmatrix} e_{1,k+2} \\ \vdots \\ e_{k+1,k+2} \end{pmatrix} \tag{6.6}$$

with $H_i \in L(\bar{B}, B)$ and (6.5) turns into

$$M_{k+1+l} = \begin{pmatrix} \bar{H} \\ 0 \\ \vdots \\ 0 \end{pmatrix} \cdot \bar{S}_{k+1+l} + \begin{pmatrix} 0 \\ M_{k+l} \end{pmatrix}, \tag{6.7}$$

characterized by the operator vector $\bar{H}$, which is independent of $l \geq 1$.

In the next step, we restrict columns $M_{k+1+l}$ to first $k+1$ components and we concentrate on the columns to the right of $M_{2k+1}$ to obtain the following result.

**Lemma 4:** Assume stabilization with $k \geq 0$. Then, for $l \geq k+1$

$$\begin{pmatrix} M_{1,k+1+l} \\ \vdots \\ M_{k+1,k+1+l} \end{pmatrix} = \begin{pmatrix} H_1 & & \ddots & \\ \vdots & & & \\ H_{k+1} & \cdots & & H_1 \end{pmatrix} \cdot \begin{pmatrix} \bar{S}_{k+1+l} \\ \vdots \\ \bar{S}_{1+l} \end{pmatrix} =: \boxed{H} \cdot \begin{pmatrix} \bar{S}_{k+1+l} \\ \vdots \\ \bar{S}_{1+l} \end{pmatrix}. \tag{6.8}$$

Lemma 4 means that first $k+1$ components of columns $M_{2k+2}$ , $M_{2k+3}$ , … show a similar structure in the sense that multiplication of a fixed operator matrix $\boxed{H}$ with a vector composed of $k+1$ operators $\bar{S}_{1+l}$ , … , $\bar{S}_{k+1+l}$ is performed. But in contrast to the matrix $\boxed{E}$ , this kind of repetition is not valid starting with column $k+2$, but only valid for columns to the right of $M_{2k+1}$. Hence, we see a delay of exactly $k$ columns to occur until the structure within $\boxed{E}$ , caused by stabilization, is completely transferred to matrix $\boxed{M}$ .

Schematically, the matrix $\boxed{M}$ is composed of regions with different properties, as indicated below by different colours.

$$\begin{array}{ccccccccc} M_1 & M_2 & \cdots & M_{k+1} & M_{k+2} & \cdots & M_{2k+1} & M_{2k+2} & \cdots \\ \hline \end{array} \tag{6.9}$$

$$\boxed{M} = \begin{pmatrix} \boxed{I_B} & M_{1,2} & \cdots & M_{1,k+1} & M_{1,k+2} & \cdots & M_{1,2k+1} & M_{1,2k+2} & M_{1,2k+3} \\ & \ddots & & \vdots & \vdots & \vdots & \vdots & \vdots & \vdots \\ & & \boxed{I_B} & M_{k,k+1} & M_{k,k+2} & \cdots & M_{k,2k+1} & M_{k,2k+2} & M_{k+1,2k+3} \\ & \boxed{\phi(\varepsilon)} & \longrightarrow & \boxed{I_B} & \boxed{\phi_1} & \cdots & \boxed{\phi_k} & \boxed{\phi_{k+1}} & \boxed{\phi_{k+2}} \\ & & & & \boxed{I_B} & \cdots & \cdots & \boxed{\phi_k} & \boxed{\phi_{k+1}} \\ & & & & & \ddots & \vdots & \vdots & \vdots \end{pmatrix}$$

Top left, the black matrix $\boxed{M}^{k+1}$ characterizing the process up to stabilization at $k \geq 0$. Hence, these operators are needed to define the decompositions of $B$ and $\bar{B}$ and thus to construct the subspace $\bar{R}_k$ of leading coefficients in $\bar{B}$ with associated Jordan chains of rank 1 to $k$ in $B$.

Secondly, $k$ red marked columns $M_{k+2}, \dots, M_{2k+1}$ with corresponding coefficients $\phi_1, \dots, \phi_k$ in row $k+1$. These coefficients are still depending from the process up to stabilization at $k$.

Thirdly, the brown columns $M_{2k+2}, M_{2k+3}, \dots$ with coefficients $\phi_{k+1}, \phi_{k+2}, \dots$, characterized by the simple repetition structure in (6.8) that eventually will allow us to derive the defining equation of $\phi(\varepsilon) = \sum_{i=0}^{\infty} \varepsilon^i \cdot \phi_i$ looked for.

Finally, we have blue Töplitz copies of row $k+1$ in the lower part of $\boxed{M}$.

**Proof of Lemma 4:** The abbreviation $\Theta_i$ is used for the operator zero column vector with $i \geq 2$ components. Then, from (6.7) we obtain backwards step-by-step

$$M_{2k+2} \overset{l=k+1}{=} \begin{pmatrix} \bar{H} \\ \Theta_{k+1} \end{pmatrix} \cdot \bar{S}_{2k+2} + \begin{pmatrix} 0 \\ M_{2k+1} \end{pmatrix} \tag{6.10}$$

$$\overset{l=k}{=} \begin{pmatrix} \bar{H} \\ \Theta_{k+1} \end{pmatrix} \cdot \bar{S}_{2k+2} + \begin{pmatrix} 0 \\ \begin{pmatrix} \bar{H} \\ \Theta_k \end{pmatrix} \cdot \bar{S}_{2k+1} + \begin{pmatrix} 0 \\ M_{2k} \end{pmatrix} \end{pmatrix}$$

$$= \begin{pmatrix} \bar{H} \\ \Theta_{k+1} \end{pmatrix} \cdot \bar{S}_{2k+2} + \begin{pmatrix} 0 \\ \bar{H} \\ \Theta_k \end{pmatrix} \cdot \bar{S}_{2k+1} + \begin{pmatrix} \Theta_2 \\ M_{2k} \end{pmatrix}$$

$$\overset{l=k-1}{=} \cdots$$

$$\overset{l=1}{=} \begin{pmatrix} \bar{H} \\ \Theta_{k+1} \end{pmatrix} \cdot \bar{S}_{2k+2} + \begin{pmatrix} 0 \\ \bar{H} \\ \Theta_k \end{pmatrix} \cdot \bar{S}_{2k+1} + \cdots + \begin{pmatrix} \Theta_k \\ \bar{H} \\ 0 \end{pmatrix} \cdot \bar{S}_{k+2} + \begin{pmatrix} \Theta_{k+1} \\ M_{k+1} \end{pmatrix}$$

and restricting to first $k+1$ components

$$\begin{pmatrix} M_{1,2k+2} \\ \vdots \\ M_{k+1,2k+2} \end{pmatrix} = \begin{pmatrix} H_1 \\ \vdots \\ \vdots \\ H_{k+1} \end{pmatrix} \cdot \bar{S}_{2k+2} + \begin{pmatrix} 0 \\ H_1 \\ \vdots \\ H_k \end{pmatrix} \cdot \bar{S}_{2k+1} + \cdots + \begin{pmatrix} 0 \\ \vdots \\ 0 \\ H_1 \end{pmatrix} \cdot \bar{S}_{k+2} \tag{6.11}$$

$$= \begin{pmatrix} H_1 & & \\ \vdots & \ddots & \\ H_k & \cdots & H_1 \end{pmatrix} \cdot \begin{pmatrix} \bar{S}_{2k+2} \\ \vdots \\ \bar{S}_{k+2} \end{pmatrix}$$

we see that Lemma 4 is proved for $l = k+1$. Now, we proceed by induction with respect to $l \geq k+1$. Assume (6.8) with $l$ replaced by $l-1$ according to

$$\begin{pmatrix} M_{1,k+l} \\ \vdots \\ M_{k+1,k+l} \end{pmatrix} = \begin{pmatrix} H_1 & & \\ \vdots & \ddots & \\ H_{k+1} & \cdots & H_1 \end{pmatrix} \cdot \begin{pmatrix} \bar{S}_{k+l} \\ \vdots \\ \bar{S}_l \end{pmatrix}. \tag{6.12}$$

Then, we obtain from (6.7) under consideration of first $k+1$ components

$$\begin{pmatrix} M_{1,k+1+l} \\ \vdots \\ \vdots \\ M_{k+1,k+1+l} \end{pmatrix} = \begin{pmatrix} H_1 \\ \vdots \\ \vdots \\ H_{k+1} \end{pmatrix} \cdot \bar{S}_{k+1+l} + \begin{pmatrix} 0 \\ M_{1,k+l} \\ \vdots \\ M_{k,k+l} \end{pmatrix} \tag{6.13}$$

and using components 1 to $k$ in (6.12), Lemma 4 is proved according to

$$\begin{pmatrix} M_{1,k+1+l} \\ \vdots \\ \vdots \\ M_{k+1,k+1+l} \end{pmatrix} = \begin{pmatrix} H_1 \\ \vdots \\ \vdots \\ H_{k+1} \end{pmatrix} \cdot \bar{S}_{k+1+l} + \begin{pmatrix} 0 \\ \begin{pmatrix} H_1 & & \\ \vdots & \ddots & \\ H_k & \cdots & H_1 \end{pmatrix} \cdot \begin{pmatrix} \bar{S}_{k+l} \\ \vdots \\ \bar{S}_{1+l} \end{pmatrix} \end{pmatrix} \tag{6.14}$$

$$= \begin{pmatrix} H_1 & & \\ \vdots & \ddots & \\ H_{k+1} & \cdots & H_1 \end{pmatrix} \cdot \begin{pmatrix} \bar{S}_{k+1+l} \\ \vdots \\ \bar{S}_{1+l} \end{pmatrix}.$$

Now, the construction of the defining equation is straightforward. The repetition structure in $\boxed{M}$ starts with brown column $2k+2$. We define a vector of operator power series by

$$\bar{d}(\varepsilon) := \begin{pmatrix} M_{1,2k+2} \\ \vdots \\ \underbrace{M_{k+1,2k+2}}_{=\phi_{k+1}} \end{pmatrix} + \varepsilon \cdot \begin{pmatrix} M_{1,2k+3} \\ \vdots \\ \underbrace{M_{k+1,2k+3}}_{=\phi_{k+2}} \end{pmatrix} + \cdots \tag{6.15}$$

$$= \sum_{l=k+1}^{\infty} \varepsilon^{l-k-1} \cdot \begin{pmatrix} M_{1,k+1+l} \\ \vdots \\ M_{k+1,k+1+l} \end{pmatrix}$$

and the transformation $\phi(\varepsilon)$ from (5.23) can be written as

$$\phi(\varepsilon) = \sum_{i=0}^{\infty} \varepsilon^i \cdot \phi_i = I_B + \varepsilon \cdot \phi_1 + \cdots + \varepsilon^k \cdot \phi_k + \varepsilon^{k+1} \cdot d_{k+1}(\varepsilon) \tag{6.16}$$

with $d_{k+1}(\varepsilon)$ denoting component $k+1$ of the operator vector $\bar{d}(\varepsilon)$. Now, by Lemma 4

$$\bar{d}(\varepsilon) = \sum_{l=k+1}^{\infty} \varepsilon^{l-k-1} \cdot \begin{pmatrix} H_1 & & \\ \vdots & \ddots & \\ H_{k+1} & \cdots & H_1 \end{pmatrix} \cdot \begin{pmatrix} \bar{S}_{k+1+l} \\ \vdots \\ \bar{S}_{1+l} \end{pmatrix} \tag{6.17}$$

$$= \boxed{H} \cdot \sum_{l=k+1}^{\infty} \varepsilon^{l-k-1} \cdot \begin{pmatrix} \bar{S}_{k+1+l} \\ \vdots \\ \bar{S}_{1+l} \end{pmatrix}$$

$$= \boxed{H} \cdot \left[ \begin{pmatrix} \bar{S}_{2k+2} \\ \bar{S}_{2k+1} \\ \vdots \\ \bar{S}_{k+3} \\ \bar{S}_{k+2} \end{pmatrix} + \varepsilon \cdot \begin{pmatrix} \bar{S}_{2k+3} \\ \bar{S}_{2k+2} \\ \vdots \\ \bar{S}_{k+4} \\ \bar{S}_{k+3} \end{pmatrix} + \cdots + \varepsilon^{k-1} \cdot \begin{pmatrix} \bar{S}_{2k+1+k} \\ \bar{S}_{2k+k} \\ \vdots \\ \bar{S}_{2k+2} \\ \bar{S}_{2k+1} \end{pmatrix} + \varepsilon^{k} \cdot \begin{pmatrix} \bar{S}_{2k+2+k} \\ \bar{S}_{2k+1+k} \\ \vdots \\ \bar{S}_{2k+3} \\ \bar{S}_{2k+2} \end{pmatrix} + \cdots \right]$$

where each component contains the same power series, starting with $\bar{S}_{2k+2}$, according to

$$c(\varepsilon) \coloneqq \bar{S}_{2k+2} + \varepsilon \cdot \bar{S}_{2k+3} + \cdots = \sum_{i=0}^{\infty} \varepsilon^{i} \cdot \bar{S}_{2k+2+i} \quad . \tag{6.18}$$

Then, the vector of power series $\bar{d}(\varepsilon)$ may be splitted by

$$\bar{d}(\varepsilon) = \boxed{H} \cdot \begin{pmatrix} 0 \\ \bar{S}_{2k+1} \\ \vdots \\ \bar{S}_{k+3} + \varepsilon \cdot \bar{S}_{k+4} + \cdots + \varepsilon^{k-2} \cdot \bar{S}_{2k+1} \\ \bar{S}_{k+2} + \varepsilon \cdot \bar{S}_{k+3} + \cdots + \varepsilon^{k-2} \cdot \bar{S}_{2k} + \varepsilon^{k-1} \cdot \bar{S}_{2k+1} \end{pmatrix} + \boxed{H} \cdot \begin{pmatrix} 1 \\ \varepsilon \\ \vdots \\ \varepsilon^{k-1} \\ \varepsilon^{k} \end{pmatrix} \cdot c(\varepsilon)$$

$$=: \bar{p}(\varepsilon) + \boxed{H} \cdot \begin{pmatrix} 1 \\ \varepsilon \\ \vdots \\ \varepsilon^{k} \end{pmatrix} \cdot c(\varepsilon) \tag{6.19}$$

with a $k+1$ dimensional operator vector of polynomials $\bar{p}(\varepsilon)$ of degree $k-1$ with respect to $\varepsilon$. In the next step, the operator coefficients $\bar{S}_{2k+2+i}$ of $c(\varepsilon)$ in (6.18) are replaced by their definitions $\bar{S}_1 = L_0$ and $\bar{S}_{i+1} = [L_1 \ldots L_i] \cdot M_i\,, i \geq 1$ from (3.8), to obtain

$$\bar{d}(\varepsilon) = \bar{p}(\varepsilon) + \boxed{H} \cdot \begin{pmatrix} 1 \\ \varepsilon \\ \vdots \\ \varepsilon^{k} \end{pmatrix} \cdot \sum_{i=0}^{\infty} \varepsilon^{i} \cdot \overbrace{[L_1 \ldots L_{2k+1+i}] \cdot M_{2k+1+i}}^{= \bar{S}_{2k+2+i}} \tag{6.20}$$

$$= \bar{p}(\varepsilon) + \boxed{H} \cdot \begin{pmatrix} 1 \\ \varepsilon \\ \vdots \\ \varepsilon^{k} \end{pmatrix} \cdot \left( [L_1 \ldots L_{2k+1}] \cdot M_{2k+1} + \varepsilon \cdot [L_1 \ldots L_{2k+2}] \cdot M_{2k+2} + \cdots \right)$$

$$\overset{(6.9)}{=} \bar{p}(\varepsilon) + \boxed{H} \cdot \begin{pmatrix} 1 \\ \varepsilon \\ \vdots \\ \varepsilon^k \end{pmatrix} \cdot \underbrace{\left( [L_1 \dots L_{2k+1}] \cdot \begin{pmatrix} M_{1,2k+1} \\ \vdots \\ M_{k+1,2k+1} \\ \phi_{k-1} \\ \vdots \\ I_B \end{pmatrix} + \varepsilon \cdot [L_1 \dots L_{2k+2}] \cdot \begin{pmatrix} M_{1,2k+2} \\ \vdots \\ M_{k+1,2k+2} \\ \phi_k \\ \vdots \\ I_B \end{pmatrix} + \cdots \right)}_{= c(\varepsilon)} .$$

Then, motivated by different parts within the matrix $\boxed{M}$ in (6.9), the column vectors of $c(\varepsilon)$ can be splitted according to

$$c(\varepsilon) = [L_1 \dots L_{k+1}] \cdot \begin{pmatrix} M_{1,2k+1} \\ \vdots \\ M_{k+1,2k+1} \end{pmatrix} + [L_{k+2} \dots L_{2k+1}] \cdot \begin{pmatrix} \phi_{k-1} \\ \vdots \\ I_B \end{pmatrix} \tag{6.21}$$

$$+ \varepsilon \cdot [L_1 \dots L_{k+1}] \cdot \begin{pmatrix} M_{1,2k+2} \\ \vdots \\ M_{k+1,2k+2} \end{pmatrix} + \varepsilon \cdot [L_{k+2} \dots L_{2k+2}] \cdot \begin{pmatrix} \phi_k \\ \vdots \\ I_B \end{pmatrix}$$

$$+ \varepsilon^2 \cdot [L_1 \dots L_{k+1}] \cdot \begin{pmatrix} M_{1,2k+3} \\ \vdots \\ M_{k+1,2k+3} \end{pmatrix} + \varepsilon^2 \cdot L_{k+2} \cdot \phi_{k+1} + \varepsilon^2 \cdot [L_{k+3} \dots L_{2k+3}] \cdot \begin{pmatrix} \phi_k \\ \vdots \\ I_B \end{pmatrix}$$

$$+ \varepsilon^3 \cdot [L_1 \dots L_{k+1}] \cdot \begin{pmatrix} M_{1,2k+4} \\ \vdots \\ M_{k+1,2k+4} \end{pmatrix} + \varepsilon^3 \cdot [L_{k+2}\ L_{k+3}] \cdot \begin{pmatrix} \phi_{k+2} \\ \phi_{k+1} \end{pmatrix} + \varepsilon^3 \cdot [L_{k+4} \dots L_{2k+3}] \cdot \begin{pmatrix} \phi_k \\ \vdots \\ I_B \end{pmatrix}$$

$$+ \cdots$$

and collected by colours we obtain

$$c(\varepsilon) = [L_1 \dots L_{k+1}] \cdot \begin{pmatrix} M_{1,2k+1} \\ \vdots \\ M_{k+1,2k+1} \end{pmatrix} + \varepsilon \cdot [L_1 \dots L_{k+1}] \cdot \bar{d}(\varepsilon)$$

$$+ \varepsilon^2 \cdot \left( L_{k+2} \cdot \phi_{k+1} + \varepsilon \cdot [L_{k+2}\ L_{k+3}] \cdot \begin{pmatrix} \phi_{k+2} \\ \phi_{k+1} \end{pmatrix} + \cdots \right)$$

$$+ \left( [0\ \ L_{k+2} \dots L_{2k+1}] + \varepsilon \cdot [L_{k+2} \dots L_{2k+2}] + \cdots \right) \cdot \begin{pmatrix} \phi_k \\ \vdots \\ I_B \end{pmatrix}$$

$$= [L_1 \dots L_{k+1}] \cdot \begin{pmatrix} M_{1,2k+1} \\ \vdots \\ M_{k+1,2k+1} \end{pmatrix} + \varepsilon \cdot [L_1 \dots L_{k+1}] \cdot \bar{d}(\varepsilon) \tag{6.22}$$

$$+\ \varepsilon^2 \cdot (L_{k+2} + \varepsilon \cdot L_{k+3} + \cdots) \cdot \underbrace{(\phi_{k+1} + \varepsilon \cdot \phi_{k+2} + \cdots)}_{= d_{k+1}(\varepsilon)}$$

$$+\ (\,[0\ \ L_{k+2} \ldots L_{2k+1}] + \varepsilon \cdot [L_{k+2} \ldots L_{2k+2}] + \cdots) \cdot \begin{pmatrix} \phi_k \\ \vdots \\ I_B \end{pmatrix}$$

under consideration of (6.16) and the definition in (6.15). Hence, using (6.20) we arrive at

$$\bar{d}(\varepsilon) = \bar{p}(\varepsilon)\ + \boxed{H} \cdot \begin{pmatrix} 1 \\ \varepsilon \\ \vdots \\ \varepsilon^k \end{pmatrix} \cdot c(\varepsilon) \tag{6.23}$$

$$= \bar{p}(\varepsilon)\ + \boxed{H} \cdot \begin{pmatrix} 1 \\ \varepsilon \\ \vdots \\ \varepsilon^k \end{pmatrix} \cdot [L_1 \ldots L_{k+1}] \cdot \begin{pmatrix} M_{1,2k+1} \\ \vdots \\ M_{k+1,2k+1} \end{pmatrix}$$

$$+ \boxed{H} \cdot \begin{pmatrix} 1 \\ \varepsilon \\ \vdots \\ \varepsilon^k \end{pmatrix} \cdot (\,[0\ \ L_{k+2} \ldots L_{2k+1}] + \varepsilon \cdot [L_{k+2} \ldots L_{2k+2}] + \cdots) \cdot \begin{pmatrix} \phi_k \\ \vdots \\ I_B \end{pmatrix}$$

$$+ \boxed{H} \cdot \begin{pmatrix} 1 \\ \varepsilon \\ \vdots \\ \varepsilon^k \end{pmatrix} \cdot (\,\varepsilon \cdot [L_1\ \ldots\ L_{k+1}] \cdot \bar{d}(\varepsilon) + \varepsilon^2 \cdot [0\ \ldots\ 0 \mid L_{k+2} + \varepsilon \cdot L_{k+3} + \cdots] \cdot \underbrace{\begin{pmatrix} d_1(\varepsilon) \\ \vdots \\ d_{k+1}(\varepsilon) \end{pmatrix}}_{= \bar{d}(\varepsilon)})$$

$$=:\ \bar{q}(\varepsilon) + \boxed{H} \cdot \begin{pmatrix} 1 \\ \varepsilon \\ \vdots \\ \varepsilon^k \end{pmatrix} \cdot (\,[L_1\ \ldots\ L_{k+1}] + \varepsilon \cdot [0\ \ldots\ 0 \mid L_{k+2} + \varepsilon \cdot L_{k+3} + \cdots]\,) \cdot \varepsilon \cdot \bar{d}(\varepsilon)$$

$$= \bar{q}(\varepsilon) + \boxed{H} \cdot \begin{pmatrix} 1 \\ \varepsilon \\ \vdots \\ \varepsilon^k \end{pmatrix} \cdot [L_1\ \ldots\ L_k \mid L_{k+1} + \varepsilon \cdot L_{k+2} + \varepsilon^2 \cdot L_{k+3} + \cdots] \cdot \varepsilon \cdot \bar{d}(\varepsilon)$$

$$=:\ \bar{q}(\varepsilon) +\ Q(\varepsilon) \cdot \varepsilon \cdot \bar{d}(\varepsilon)$$

with $\bar{q}(\varepsilon)$ an operator power series vector with $k+1$ components and $Q(\varepsilon)$ a $(k+1) \times (k+1)$ operator power series matrix.

Now, by (6.23), the defining equation for the power series vector $\bar{d}(\varepsilon)$ reads

$$\bar{d}(\varepsilon) = \bar{q}(\varepsilon) +\ \varepsilon \cdot Q(\varepsilon) \cdot \bar{d}(\varepsilon) \qquad or \qquad [\,I_{B^{k+1}}\ -\ \varepsilon \cdot Q(\varepsilon)\,] \cdot \bar{d}(\varepsilon) = \bar{q}(\varepsilon) \tag{6.24}$$

with $I_{B^{k+1}} \coloneqq Diag(I_B \cdots I_B)$. Note that the representation of $\bar{d}(\varepsilon)$ by some sort of regular part $\bar{q}(\varepsilon)$ plus a small perturbation, caused by the factor $\varepsilon$, reminds a Nakayama Lemma type representation of $\bar{d}(\varepsilon)$ that allows for solution of the equation.

Note also that by (6.6), (6.19) and (6.23), the perturbation matrix $Q(\varepsilon)$ depends from elements of the black matrix $\boxed{M}^{k+1}$ in (6.9), whereas the inhomogeneity $\bar{q}(\varepsilon)$ depends from elements of first $k+1$ rows within the matrix $\boxed{M}^{2k+1}$(black and red parts in (6.9)).

In addition, $Q(\varepsilon)$ and $\bar{q}(\varepsilon)$ depend of the power series

$$\bar{L}_{k+1}(\varepsilon) := L_{k+1} + \varepsilon \cdot L_{k+2} + \varepsilon^2 \cdot L_{k+3} + \cdots \tag{6.25}$$

and

$$\begin{aligned} \bar{L}_{k+2}(\varepsilon) &:= L_{k+2} + \varepsilon \cdot L_{k+3} + \varepsilon^2 \cdot L_{k+4} + \cdots \\ &\vdots \\ \bar{L}_{2k+1}(\varepsilon) &:= L_{2k+1} + \varepsilon \cdot L_{2k+2} + \varepsilon^2 \cdot L_{2k+3} + \cdots \end{aligned} \tag{6.26}$$

respectively, representing certain sections of the given power series $L(\varepsilon) = \sum_{i=0}^{\infty} \varepsilon^i \cdot L_i$.

Finally by definition, $\bar{d}(\varepsilon)$ is the power series composed of first $k+1$ components of the columns $M_{2k+2}$ , $M_{2k+3}$ , ... (brown part in (6.9)) and by (6.16) the transformation $\phi(\varepsilon) = I_B + \cdots + \varepsilon^k \cdot \phi_k + \varepsilon^{k+1} \cdot d_{k+1}(\varepsilon)$ is obtained by use of the last component $d_{k+1}(\varepsilon)$ of $\bar{d}(\varepsilon)$.

Now, concerning explicit calculation of the power series coefficients, we first note that the unique left and right inverse power series of $I_{B^{k+1}} - \varepsilon \cdot Q(\varepsilon)$ reads by Neumann series

$$\begin{aligned} \Gamma(\varepsilon) := \left[\, I_{B^{k+1}} - \varepsilon \cdot Q(\varepsilon) \,\right]^{-1} &= I_{B^{k+1}} + \varepsilon \cdot Q(\varepsilon) + \varepsilon^2 \cdot Q(\varepsilon)^2 + \cdots \\ &= I_{B^{k+1}} + \sum_{i=1}^{\infty} \varepsilon^i \cdot \Gamma_i \end{aligned} \tag{6.27}$$

with coefficients $\Gamma_i$ explicitly given by the coefficients of $Q(\varepsilon) = \sum_{i=0}^{\infty} \varepsilon^i \cdot Q_{i+1}$ according to

$$\begin{aligned} \Gamma_0 &= I_{B^{k+1}} \\ \Gamma_1 = Q_1 = Q(0) &= \begin{pmatrix} H_1 \\ \vdots \\ H_{k+1} \end{pmatrix} \cdot [L_1 \;\; \ldots \;\; L_k \;\; L_{k+1}] \\ i \geq 2 \;:\; \Gamma_i &= \sum_{\substack{k=1 \\ j_1+\cdots+j_k=i}}^{n} Q_{j_1} \cdot \cdots \cdot Q_{j_k} \end{aligned} \tag{6.28}$$

with $k$-tuple $(j_1, \ldots, j_k)$ and $j_1 \geq 1, \ldots, j_k \geq 1$. Using $\bar{q}(\varepsilon) = \sum_{j=0}^{\infty} \varepsilon^j \cdot \bar{q}_j$, the power series $\bar{d}(\varepsilon)$ may finally be rewritten by Cauchy product in the form

$$\begin{aligned} \bar{d}(\varepsilon) = \Gamma(\varepsilon) \cdot \bar{q}(\varepsilon) &= \left( I_{B^{k+1}} + \sum_{i=1}^{\infty} \varepsilon^i \cdot \Gamma_i \right) \cdot \left( \bar{q}_0 + \sum_{j=1}^{\infty} \varepsilon^j \cdot \bar{q}_j \right) \\ &= \bar{q}_0 + \sum_{l=1}^{\infty} \varepsilon^l \cdot \sum_{i+j=l} \Gamma_i \cdot \bar{q}_j \,. \end{aligned} \tag{6.29}$$

## 7. Triangularization and Diagonalization of analytic Operator Functions

In this section, we restrict to analytic operator families $L(\varepsilon)=\sum_{i=0}^{\infty}\varepsilon^i\cdot L_i\,,\; L\in C^{\omega}(U,L[B,\bar{B}])$ with $B,\bar{B}$ real or complex Banach spaces, $0\in U\subset\mathbb{K}=\mathbb{R},\mathbb{C}$ and $L[B,\bar{B}]$ bounded linear operators from $B$ to $\bar{B}$.

Further, assume stabilization at $k\geq 0$ according to (5.3). Then, the defining equation

$$\bar{d}=\bar{q}(\varepsilon)+\;\varepsilon\cdot Q(\varepsilon)\cdot\bar{d}\;,\;\;\bar{d}\in B^{k+1} \tag{7.1}$$

of $\bar{d}(\varepsilon)$ from (6.24) turns into an analytic equation, as soon as analyticity of $\bar{q}(\varepsilon)$ and $Q(\varepsilon)$ is assured. Now, these operator power series depend from a finite combination of operators selected from $\boxed{M}^{2k+1}$ and the power series (6.25), (6.26) respectively. In some more detail, the elements of $\boxed{M}^{2k+1}$ result from the iteration (3.3)-(3.15) up to $2k$, depending from $2k+1$ operator coefficients $L_0\,,\dots,L_{2k}$ of $L(\varepsilon)$ as well as from the projections to the subspaces of the direct sums

$$\begin{aligned} B\;&=\;N_1^c\;\oplus\;\cdots\;\oplus\;N_{k+1}^c\;\oplus\;N_{k+1}\\ \bar{B}\;&=\;R_1\;\oplus\;\cdots\;\oplus\;R_{k+1}\;\oplus\;R_{k+1}^c \end{aligned}\;. \tag{7.2}$$

Hence, to assure analyticity of $\bar{q}(\varepsilon)$ and $Q(\varepsilon)$ in (7.1), above all we have to assume continuity of these projections or equivalently, we have to assume closedness of the subspaces. This is well known from [Bart], [Gohberg] and [Kaballo], where the existence of complemented subspaces $N_{k+1}$ and $R_{k+1}^c$ is required too.

It remains to consider analyticity of the power series $\bar{L}_{k+1}(\varepsilon)\,,\dots,\bar{L}_{2k+1}(\varepsilon)$ in (6.25) and (6.26). Now, analyticity of $L(\varepsilon)$ implies by Taylor's formula

$$L(\varepsilon)=L_0+\cdots+\varepsilon^k\cdot L_k+\underbrace{\varepsilon^{k+1}\cdot\frac{1}{k!}\cdot\int_0^1(1-\tau)^k\,\overbrace{L^{(k)}(\tau\cdot\varepsilon)}^{analytic}\;d\tau}_{=\,\bar{L}_{k+1}(\varepsilon)} \tag{7.3}$$

yielding analyticity of $\bar{L}_{k+1}(\varepsilon)$, as desired. Analogously, for $\bar{L}_{k+2}(\varepsilon)\,,\dots,\bar{L}_{2k+1}(\varepsilon)$. Thus, equation (7.1) has a unique analytic solution

$$\bar{d}(\varepsilon)=\left[\,I_{B^{k+1}}\;-\;\varepsilon\cdot Q(\varepsilon)\,\right]^{-1}\cdot\bar{q}(\varepsilon) \tag{7.4}$$

with coefficients given by (6.28), (6.29) and the analytic diffeomorphic transformation $\phi(\varepsilon)$ is obtained from component $k+1$ of $\bar{d}(\varepsilon)$ according to

$$\phi(\varepsilon)=I_B+\varepsilon\cdot\phi_1+\cdots+\varepsilon^k\cdot\phi_k+\varepsilon^{k+1}\cdot d_{k+1}(\varepsilon)\;. \tag{7.5}$$

In summary, combining closedness of subspaces in (7.2) with analyticity of $L(\varepsilon)$, we obtain analyticity of the transformation $\phi(\varepsilon)$ and, obviously, analyticity of the normal form $S(\varepsilon)=L(\varepsilon)\cdot\phi(\varepsilon)$ from Theorem 2 (i).

Next, consider the factorization $S(\varepsilon)=\psi(\varepsilon)\cdot\Delta(\varepsilon)$ from Theorem 2 (iv). We will show that analyticity of $S(\varepsilon)$ implies analyticity of the power series $\psi(\varepsilon)$. Again, closedness of the subspaces in (7.2) is assumed. By (5.43) we have

$$\psi(\varepsilon) = I_{\bar{B}} + \varepsilon \cdot \underbrace{\left( \sum_{i=1}^{\infty} \varepsilon^{i-1} \cdot S_{i+1} \right)}_{= r_1^c(\varepsilon)} \cdot S_1^{-1} \mathcal{P}_1 + \cdots + \varepsilon \cdot \underbrace{\left( \sum_{i=1}^{\infty} \varepsilon^{i-1} \cdot S_{k+i+1} \right)}_{= r_{k+1}^c(\varepsilon)} \cdot S_{k+1}^{-1} \mathcal{P}_{k+1} \tag{7.6}$$

with power series $r_i^c(\varepsilon), i = 1, \dots, k+1$ that are analytic by Taylor's formula according to

$$S(\varepsilon) = \sum_{i=0}^{\infty} \varepsilon^i \cdot S_{i+1} = S_1 + \varepsilon \cdot \underbrace{\left( \sum_{i=1}^{\infty} \varepsilon^{i-1} \cdot S_{i+1} \right)}_{= r_1^c(\varepsilon)} = S_1 + \varepsilon \cdot \underbrace{\int_0^1 \overbrace{S(\tau \cdot \varepsilon)}^{analytic} \, d\tau}_{analytic} . \tag{7.7}$$

Analogously for $r_2^c(\varepsilon), \dots, r_{k+1}^c(\varepsilon)$. The diagonal power polynomial $\Delta(\varepsilon) = S_1 P_1 + \varepsilon \cdot S_2 P_2 + \cdots + \varepsilon^k \cdot S_{k+1} P_{k+1}$ satisfies $\Delta \in C^{\omega}(U, L[B, \bar{B}])$ by construction.

Concerning analytic operator functions $L(\varepsilon) = \sum_{i=0}^{\infty} \varepsilon^i \cdot L_i$, $L_i \in L[B, \bar{B}]$, we summarize the results in the following theorem.

**Theorem 3:** Assume stabilization at $k \geq 0$ of the analytic operator family $L(\varepsilon)$ with $B, \bar{B}$ real or complex Banach spaces and subspaces closed in (7.2).

(i) Then, $\phi(\varepsilon)$ from (7.5) transforms $L(\varepsilon)$ into the analytic normal form

$$S(\varepsilon) := L(\varepsilon) \cdot \phi(\varepsilon) = S_1 + \varepsilon^1 \cdot S_2 + \cdots + \varepsilon^k \cdot S_{k+1} + \sum_{l=2}^{\infty} \varepsilon^{k+l-1} \cdot S_{k+l} \tag{7.8}$$

with $S_1, \dots, S_{k+1}$ satisfying

$$S_l \cdot (\, N_0^c + \cdots + N_{l-1}^c \,) \subset R_{l-1}^c \qquad \wedge \qquad S_l \cdot N_l^c = R_l \qquad \wedge \qquad S_l \cdot N_l = 0 \tag{7.9}$$

for $l = 1, \dots, k+1$ and $S_{k+l}$ satisfying

$$R[\, S_{k+l} \,] \subset R_{k+1}^c \qquad \wedge \qquad S_{k+l} \cdot N_{k+1} = 0 \tag{7.10}$$

for $l \geq 2$. In particular, for all $l \geq 0$, we have $S_{l+1} \cdot N_{k+1} = 0$.

(ii) Using the analytic near identity transformation $\psi(\varepsilon)$ of $\bar{B}$ from (5.41), the family $L(\varepsilon)$ is diagonalized according to

$$\psi^{-1}(\varepsilon) \cdot L(\varepsilon) \cdot \phi(\varepsilon) = \Delta(\varepsilon) \tag{7.11}$$

with $\Delta(\varepsilon)$ given by

$$\Delta(\varepsilon) = S_1 P_1 + \varepsilon \cdot S_2 P_2 + \cdots + \varepsilon^k \cdot S_{k+1} P_{k+1} \, , \tag{7.12}$$

i.e. the family $L(\varepsilon)$ and the polynomial $\Delta(\varepsilon)$ of degree $k \geq 0$ are analytically equivalent.

(iii) Within a punctured neighbourhood $\varepsilon \in U \setminus \{0\}$, kernels and ranges of $L(\varepsilon)$ are given by

$$N[\, L(\varepsilon) \,] = \phi(\varepsilon) \cdot N_{k+1} \quad \textit{and} \quad R[\, L(\varepsilon) \,] = \psi(\varepsilon) \cdot [\, R_1 \oplus \cdots \oplus R_{k+1} \,] \, . \tag{7.13}$$

In particular, kernels $N[L(\varepsilon)], \varepsilon \neq 0$ and ranges $R[L(\varepsilon)], \varepsilon \neq 0$ are analytically embedded within the families $N(\varepsilon) = \phi(\varepsilon) \cdot N_{k+1}, \varepsilon \in U$ and $R(\varepsilon) = \psi(\varepsilon) \cdot [R_1 \oplus \cdots \oplus R_{k+1}]$, $\varepsilon \in U$, i.e. smoothing of kernels and ranges occurs.

(iv) A smooth generalized inverse of the diagonal operator polynomial $\Delta(\varepsilon)$ reads for $\varepsilon \neq 0$

$$\Delta^{-1}(\varepsilon) = \varepsilon^{-k} \cdot S_{k+1}^{-1}\mathcal{P}_{k+1} + \cdots + S_1^{-1}\mathcal{P}_1 \tag{7.14}$$

with pole of order $k \geq 0$ at $\varepsilon = 0$. Correspondingly, the family

$$L^{-1}(\varepsilon) = \phi(\varepsilon) \cdot \Delta^{-1}(\varepsilon) \cdot \psi^{-1}(\varepsilon) \tag{7.15}$$

defines a generalized inverse of $L(\varepsilon)$, which is analytic for $\varepsilon \in U \setminus \{0\}$ with pole of order $k$ at $\varepsilon = 0$. For $\varepsilon \in U \setminus \{0\}$, the analytic families

$$L^{-1}(\varepsilon) \cdot L(\varepsilon) = \phi(\varepsilon) \cdot (P_1 + \cdots + P_{k+1}) \cdot \phi^{-1}(\varepsilon) \tag{7.16}$$

$$L(\varepsilon) \cdot L^{-1}(\varepsilon) = \psi(\varepsilon) \cdot (\mathcal{P}_1 + \cdots + \mathcal{P}_{k+1}) \cdot \psi^{-1}(\varepsilon) \ .$$

are representing projections to the subspaces $\phi(\varepsilon) \cdot [N_1^c \oplus \cdots \oplus N_{k+1}^c] \subset B$ and $R[L(\varepsilon)] = \psi(\varepsilon) \cdot [R_1 \oplus \cdots \oplus R_{k+1}] \subset \bar{B}$ respectively. The families can analytically be continued to $\varepsilon = 0$ by

$$\phi(0) \cdot (P_1 + \cdots + P_{k+1}) \cdot \phi^{-1}(0) = P_1 + \cdots + P_{k+1} \tag{7.17}$$

$$\psi(0) \cdot (\mathcal{P}_1 + \cdots + \mathcal{P}_{k+1}) \cdot \psi^{-1}(0) = \mathcal{P}_1 + \cdots + \mathcal{P}_{k+1} \ .$$

**Proof of Theorem 3:** (i) and (ii) are repetitions of Theorem 2 (i) and (iv) under consideration of the analyticity of $\phi(\varepsilon)$ and $\psi(\varepsilon)$ from (7.3) and (7.7).

Concerning (iii), note that for $\varepsilon \neq 0$, we have $N[\Delta(\varepsilon)] \equiv N_{k+1}$, implying $L(\varepsilon) \cdot b = \psi(\varepsilon) \cdot \Delta(\varepsilon) \cdot \phi^{-1}(\varepsilon) \cdot b = 0$ iff $\phi^{-1}(\varepsilon) \cdot b \in N_{k+1}$ or $N[L(\varepsilon)] = \phi(\varepsilon) \cdot N_{k+1}$.

Again for $\varepsilon \neq 0$, the identity $R[\Delta(\varepsilon)] \equiv R_1 \oplus \cdots \oplus R_{k+1}$ is true, yielding $R[L(\varepsilon)] = R[\psi(\varepsilon) \cdot \Delta(\varepsilon) \cdot \phi^{-1}(\varepsilon)] = \psi(\varepsilon) \cdot [R_1 \oplus \cdots \oplus R_{k+1}]$, as desired.

Concerning (iv), note that

$$\begin{aligned}
&\Delta^{-1}(\varepsilon) \cdot \Delta(\varepsilon) \\
&= [\, \varepsilon^{-k} \cdot S_{k+1}^{-1}\mathcal{P}_{k+1} + \cdots + S_1^{-1}\mathcal{P}_1 \,] \cdot [\, S_1P_1 + \varepsilon \cdot S_2P_2 + \cdots + \varepsilon^k \cdot S_{k+1}P_{k+1} \,] \qquad (7.18) \\
&= S_1^{-1}\mathcal{P}_1 \cdot S_1P_1 + \cdots + S_{k+1}^{-1}\mathcal{P}_{k+1} \cdot S_{k+1}P_{k+1} \\
&= I_{N_1^c} \cdot P_1 + \cdots + I_{N_{k+1}^c} \cdot P_{k+1} \\
&= P_1 + \cdots + P_{k+1}
\end{aligned}$$

and

$$\begin{aligned}
&\Delta(\varepsilon) \cdot \Delta^{-1}(\varepsilon) \\
&= [\, S_1P_1 + \varepsilon \cdot S_2P_2 + \cdots + \varepsilon^k \cdot S_{k+1}P_{k+1} \,] \cdot [\, \varepsilon^{-k} \cdot S_{k+1}^{-1}\mathcal{P}_{k+1} + \cdots + S_1^{-1}\mathcal{P}_1 \,] \qquad (7.19) \\
&= S_1P_1 \cdot S_1^{-1}\mathcal{P}_1 + \cdots + S_{k+1}P_{k+1} \cdot S_{k+1}^{-1}\mathcal{P}_{k+1} \\
&= I_{R_1} \cdot \mathcal{P}_1 + \cdots + I_{R_{k+1}} \cdot \mathcal{P}_{k+1}
\end{aligned}$$

$$= \mathcal{P}_1 + \cdots + \mathcal{P}_{k+1}$$

yielding

$$\Delta(\varepsilon) \cdot \Delta^{-1}(\varepsilon) \cdot \Delta(\varepsilon)$$

$$= \underbrace{[\mathcal{P}_1 + \cdots + \mathcal{P}_{k+1}]}_{projection\ to\ R_1 \oplus \cdots \oplus R_{k+1}} \cdot \underbrace{[S_1 P_1 + \varepsilon \cdot S_2 P_2 + \cdots + \varepsilon^k \cdot S_{k+1} P_{k+1}]}_{\in R_1 \oplus \cdots \oplus R_{k+1}} \tag{7.20}$$

$$= S_1 P_1 + \varepsilon \cdot S_2 P_2 + \cdots + \varepsilon^k \cdot S_{k+1} P_{k+1}$$

$$= \Delta(\varepsilon)$$

and

$$\Delta^{-1}(\varepsilon) \cdot \Delta(\varepsilon) \cdot \Delta^{-1}(\varepsilon)$$

$$= \underbrace{[P_1 + \cdots + P_{k+1}]}_{projection\ to\ N_1^c \oplus \cdots \oplus N_{k+1}^c} \cdot \underbrace{[S_1^{-1} \mathcal{P}_1 + \cdots + \varepsilon^{-k} \cdot S_{k+1}^{-1} \mathcal{P}_{k+1}]}_{\in N_1^c \oplus \cdots \oplus N_{k+1}^c} \tag{7.21}$$

$$= S_1^{-1} \mathcal{P}_1 + \cdots + \varepsilon^{-k} \cdot S_{k+1}^{-1} \mathcal{P}_{k+1}$$

$$= \Delta^{-1}(\varepsilon)$$

and (7.14) is shown. But then,

$$L^{-1}(\varepsilon) \cdot L(\varepsilon)$$

$$= \phi(\varepsilon) \cdot \Delta^{-1}(\varepsilon) \cdot \psi^{-1}(\varepsilon) \cdot \psi(\varepsilon) \cdot \Delta(\varepsilon) \cdot \phi^{-1}(\varepsilon) \tag{7.22}$$

$$= \phi(\varepsilon) \cdot \Delta^{-1}(\varepsilon) \cdot \Delta(\varepsilon) \cdot \phi^{-1}(\varepsilon)$$

$$\overset{(7.18)}{=} \phi(\varepsilon) \cdot (P_1 + \cdots + P_{k+1}) \cdot \phi^{-1}(\varepsilon)$$

as well as

$$L(\varepsilon) \cdot L^{-1}(\varepsilon)$$

$$= \psi(\varepsilon) \cdot \Delta(\varepsilon) \cdot \phi^{-1}(\varepsilon) \cdot \phi(\varepsilon) \cdot \Delta^{-1}(\varepsilon) \cdot \psi^{-1}(\varepsilon) \tag{7.23}$$

$$= \psi(\varepsilon) \cdot \Delta(\varepsilon) \cdot \Delta^{-1}(\varepsilon) \cdot \psi^{-1}(\varepsilon)$$

$$\overset{(7.18)}{=} \psi(\varepsilon) \cdot (\mathcal{P}_1 + \cdots + \mathcal{P}_{k+1}) \cdot \psi^{-1}(\varepsilon)$$

implying (7.16). The property of $\phi(\varepsilon) \cdot (P_1 + \cdots + P_{k+1}) \cdot \phi^{-1}(\varepsilon)$ and $\psi(\varepsilon) \cdot (\mathcal{P}_1 + \cdots + \mathcal{P}_{k+1}) \cdot \psi^{-1}(\varepsilon)$ to be projections to corresponding subspaces is obvious. Finally,

$$L(\varepsilon) \cdot L^{-1}(\varepsilon) \cdot L(\varepsilon)$$

$$= \psi(\varepsilon) \cdot (\mathcal{P}_1 + \cdots + \mathcal{P}_{k+1}) \cdot \psi^{-1}(\varepsilon) \cdot \psi(\varepsilon) \cdot \Delta(\varepsilon) \cdot \phi^{-1}(\varepsilon) \tag{7.24}$$

$$= \psi(\varepsilon) \cdot \underbrace{[\,\mathcal{P}_1 + \cdots + \mathcal{P}_{k+1}\,]}_{projection\ to\ R_1 \oplus \cdots \oplus R_{k+1}} \cdot \underbrace{\Delta(\varepsilon)}_{\in R_1 \oplus \cdots \oplus R_{k+1}} \cdot \phi^{-1}(\varepsilon)$$

$$= \psi(\varepsilon) \cdot \Delta(\varepsilon) \cdot \phi^{-1}(\varepsilon)$$

$$= L(\varepsilon)$$

and

$$L^{-1}(\varepsilon) \cdot L(\varepsilon) \cdot L^{-1}(\varepsilon)$$

$$= \phi(\varepsilon) \cdot (P_1 + \cdots + P_{k+1}) \cdot \phi^{-1}(\varepsilon) \cdot \phi(\varepsilon) \cdot \Delta^{-1}(\varepsilon) \cdot \psi^{-1}(\varepsilon) \tag{7.25}$$

$$= \phi(\varepsilon) \cdot \underbrace{[\,P_1 + \cdots + P_{k+1}\,]}_{projection\ to\ N_1^c \oplus \cdots \oplus N_{k+1}^c} \cdot \underbrace{\Delta^{-1}(\varepsilon)}_{\in N_1^c \oplus \cdots \oplus N_{k+1}^c} \cdot \psi^{-1}(\varepsilon)$$

$$= \phi(\varepsilon) \cdot \Delta^{-1}(\varepsilon) \cdot \psi^{-1}(\varepsilon)$$

$$= L^{-1}(\varepsilon)$$

yielding $L^{-1}(\varepsilon)$ to be a generalized inverse of $L(\varepsilon)$.

We close the paper with some remarks concerning global properties of $L(\varepsilon)$ and possible relations to commutative algebra, see [Hauser] and [Rond].

**Remarks: 1)** The results in Theorem 3 are stated with respect to $\varepsilon = 0$ and a sufficiently small neighbourhood $U \subset \mathbb{K} = \mathbb{R}, \mathbb{C}$ of $\varepsilon = 0$, i.e. only local results are formulated. Then, if stabilization occurs at $\varepsilon = 0$ with $k_0 := k \geq 0$ and corresponding subspaces closed, smoothing of kernels and ranges is assured, i.e. root elements $N_{k+1}$ of Jordan chains with infinite rank as well as leading coefficients $R_1 \oplus \cdots \oplus R_{k+1}$ of order 0 to $k$ can smoothly be continued to closed subspaces $N[L(\varepsilon)]$ and $R[L(\varepsilon)]$. Now, by setting

$$\bar{L}(\varepsilon) := L(\bar{\varepsilon} + \varepsilon) = L(\,\bar{\varepsilon}\,) + \varepsilon \cdot L'(\,\bar{\varepsilon}\,) + \cdots \; = \; \bar{L}_0 + \varepsilon \cdot \bar{L}_1 + \cdots$$

for $\bar{\varepsilon} \neq 0$ fixed in $U$, the investigation of Theorem 3 can also be performed at $\bar{\varepsilon} \neq 0$ with trivial limiting behaviour of associated subspaces according to

$$N[\,\bar{L}(\varepsilon)\,] \to N[\,\bar{L}(0)\,] = N[\,\bar{L}_0\,] \qquad and \qquad R[\,\bar{L}(\varepsilon)\,] \to R[\,\bar{L}(0)\,] = R[\,\bar{L}_0\,]\,,$$

implying stabilization with $k_{\bar{\varepsilon}} \equiv 0$ within the punctured neighbourhood $U \setminus \{0\}$. Hence, within $U$ at most for $\varepsilon = 0$ a value of $k \geq 1$ may appear.

Let us now turn to arbitrary $\Omega \subset \mathbb{K}$, $\Omega$ open with an analytic family of Fredholm operators $L \in C^{\omega}(\Omega, L[B, \bar{B}])$ of arbitrary index. Then, for every $\varepsilon$ in $\Omega$, stabilization with $k_\varepsilon \geq 0$ is assured by finite dimensionality of (closed) subspaces and we obtain the disjoint split

$$\Omega = \underbrace{\{\,\varepsilon \in \Omega \mid k_\varepsilon = 0\,\}}_{=:\,\rho\ open} \;\cup\; \underbrace{\{\,\varepsilon \in \Omega \mid k_\varepsilon \geq 1\,\}}_{=:\,\sigma\ discrete}$$

with discrete and closed set $\sigma$. Note that $\sigma$ represents the singularities or jumping points of $k_\varepsilon$ in $\Omega$ with kernel blown up and range collapsed, when compared to the neighbourhood. In case of $L(\varepsilon) \equiv 0$, we have $N[L(\varepsilon)] \equiv B$ and $R[L(\varepsilon)] \equiv \{0\}$ and no jumping point occurs, implying $\sigma$ to be the empty set. In [Kaballo] global aspects of the set of singularities $\sigma$ are treated, even in a more general context, by weakening Fredholm properties to properties concerning closedness of kernel and ranges. In addition meromorphic families $L(\varepsilon)$ are allowed.

At this point we note that a meromorphic family $L(\varepsilon) = \varepsilon^{-n} \cdot L_{-n} + \cdots$ with pole of order $n$ at $\varepsilon = 0$ locally turns into a holomorphic family after multiplication by $\varepsilon^n$. Then, the Jordan chains of the holomorphic family turn into associated chains of $L(\varepsilon)$ and it should be possible to transfer most of the results of the previous sections to the meromorphic case.

**2)** ***Flatness of power series over analytic functions***

From Theorem 2 (ii) we know that every power series $b(\varepsilon) = \sum_{i=0}^{\infty} \varepsilon^i \cdot b_i$ satisfying $L(\varepsilon) \cdot b(\varepsilon) = 0$ can uniquely be represented by $b(\varepsilon) = \phi(\varepsilon) \cdot n_{k+1}(\varepsilon)$ with power series $n_{k+1}(\varepsilon) = \sum_{i=0}^{\infty} \varepsilon^i \cdot n_{k+1}^i$ and $n_{k+1}^i \in N_{k+1}$. Now, if $\phi(\varepsilon)$ is analytic, we obtain

$$\underbrace{b(\varepsilon)}_{\text{power series in } B} = \underbrace{\phi(\varepsilon)}_{\text{analytic in } L[B,B]} \cdot \underbrace{n_{k+1}(\varepsilon)}_{\text{power series in } N_{k+1}} \tag{7.26}$$

and the power series solution $b(\varepsilon)$ is constructed by multiplication of the power series $n_{k+1}(\varepsilon)$ with the analytic operator family $\phi \in C^\omega(U, L[B,B])$. This means that every power series solution $b(\varepsilon)$ of the analytic equation $L(\varepsilon) \cdot b = 0$ can be constructed by use of the same analytic operator $\phi(\varepsilon)$ applied to an appropriate power series from $N_{k+1}$. In finite dimensions, this property is known as the flatness of the power series over the analytic functions.

In some more detail, when restricting the Banach space $B$ to finite dimension, say $B = \mathbb{K}^n, n \geq 1$, then $\phi(\varepsilon)$ turns into an analytic $(n \times n)$-matrix family and $b(\varepsilon)$ as well as $n_{k+1}(\varepsilon)$ turn into $n$-dimensional power series vectors according to

$$b(\varepsilon) = \begin{pmatrix} b_1(\varepsilon) \\ \vdots \\ b_n(\varepsilon) \end{pmatrix} = \overbrace{\begin{pmatrix} \phi_{1,1}(\varepsilon) & \cdots & \phi_{1,n}(\varepsilon) \\ \vdots & \ddots & \vdots \\ \phi_{n,1}(\varepsilon) & \cdots & \phi_{n,n}(\varepsilon) \end{pmatrix}}^{=: \boxed{\phi(\varepsilon)}} \cdot \overbrace{\begin{pmatrix} n_1(\varepsilon) \\ \vdots \\ n_n(\varepsilon) \end{pmatrix}}^{\in N_{k+1}}$$

$$= \boxed{\phi(\varepsilon)} \cdot (\bar{n}_1 \cdot c_1(\varepsilon) + \cdots + \bar{n}_l \cdot c_l(\varepsilon))$$

$$= \boxed{\phi(\varepsilon)} \cdot \bar{n}_1 \cdot c_1(\varepsilon) + \cdots + \boxed{\phi(\varepsilon)} \cdot \bar{n}_l \cdot c_l(\varepsilon)$$

$$=: \bar{a}_1(\varepsilon) \cdot c_1(\varepsilon) + \cdots + \bar{a}_l(\varepsilon) \cdot c_l(\varepsilon) \tag{7.27}$$

with $\bar{n}_1, \ldots, \bar{n}_l$ denoting a basis of $N_{k+1}$ of dimension $1 \leq l \leq n$ and $\bar{a}_1(\varepsilon), \ldots, \bar{a}_l(\varepsilon)$ representing analytic solutions of $L(\varepsilon) \cdot b = 0$ due to

$$L(\varepsilon) \cdot \bar{a}_1(\varepsilon) = L(\varepsilon) \cdot \boxed{\phi(\varepsilon)} \cdot \bar{n}_1 = 0\,, \quad \cdots \quad, \; L(\varepsilon) \cdot \bar{a}_l(\varepsilon) = L(\varepsilon) \cdot \boxed{\phi(\varepsilon)} \cdot \bar{n}_l = 0 \tag{7.28}$$

by Theorem 3 (iii). Hence, every power series solution $b(\varepsilon)$ can be written by (7.27) as a linear combination of $l$ analytic solutions $\bar{a}_i(\varepsilon)$, independent of $b(\varepsilon)$, and $l$ coefficients $c_i(\varepsilon)$, dependent from $b(\varepsilon)$, with $c_i(\varepsilon)$ from the power series ring $\mathbb{K}[[\varepsilon]]$.

In this sense, the set $\{\bar{a}_1(\varepsilon), \dots, \bar{a}_l(\varepsilon)\}$ delivers an analytic basis, composed of analytic solutions of $L(\varepsilon) \cdot b = 0$, for linear combination of all power series solutions $b(\varepsilon)$. As already mentioned, in commutative algebra this means flatness of the power series over the analytic functions [Hauser] and (7.26) gives a corresponding result in the context of Banach spaces by simply replacing the basis by the analytic family $\phi(\varepsilon) \in L[B, B]$.

**3) *Linear analytic Artin approximation***

Consider $L(\varepsilon)$ to be analytic with assumptions from Theorem 3 and possessing a power series solution $b(\varepsilon)$ according to $L(\varepsilon) \cdot b(\varepsilon) = 0$. Then by (7.26), we obtain $b(\varepsilon) = \phi(\varepsilon) \cdot n_{k+1}(\varepsilon)$ with appropriate power series $n_{k+1}(\varepsilon) = \sum_{i=0}^{\infty} \varepsilon^i \cdot n_{k+1}^i$, $n_{k+1}^i \in N_{k+1}$, and truncation of $n_{k+1}(\varepsilon)$ at order $c \geq 0$ yields the polynomial

$$\bar{n}_{k+1}(\varepsilon) := n_{k+1}^0 + \cdots + \varepsilon^c \cdot n_{k+1}^c \,, \tag{7.29}$$

implying the identity

$$L(\varepsilon) \cdot \phi(\varepsilon) \cdot \bar{n}_{k+1}(\varepsilon) = 0 \tag{7.30}$$

due to $L(\varepsilon) \cdot \phi(\varepsilon) \cdot N_{k+1} \equiv 0$ from Theorem 3 (iii). Hence, $b^*(\varepsilon) := \phi(\varepsilon) \cdot \bar{n}_{k+1}(\varepsilon)$ represents an analytic solution agreeing up to order $c \geq 0$ with the power series solution $b(\varepsilon)$ and we see that Theorem 3 delivers a simple version of linear analytic Artin approximation [Hauser] in general Banach spaces.

**4) *Linear Strong Artin approximation, Greenberg function and Artin-Rees Lemma***

Consider a power series $L(\varepsilon) = \sum_{i=0}^{\infty} \varepsilon^i \cdot L_i$ with stabilization of the iteration (3.3)-(3.15) at $k \geq 0$, as given in Theorem 2. Then, a strong Artin approximation result for linear relations in vector spaces reads as follows.

For every approximation $b(\varepsilon) = \sum_{i=0}^{\infty} \varepsilon^i \cdot b_i$ of order $k + l$, $l \geq 1$,

$$L(\varepsilon) \cdot b(\varepsilon) = \varepsilon^{k+l} \cdot \sum_{i=0}^{\infty} \varepsilon^i \cdot \bar{b}_{k+l+i} \,, \tag{7.31}$$

there exists a power series solution $\hat{b}(\varepsilon) = \sum_{i=0}^{\infty} \varepsilon^i \cdot \hat{b}_i$ satisfying

$$L(\varepsilon) \cdot \hat{b}(\varepsilon) = 0 \qquad \textit{and} \qquad b_i = \hat{b}_i \,, \ i = 0, \dots, l-1 \,, \tag{7.32}$$

i.e. the approximation $b(\varepsilon)$ and the power series solution $\hat{b}(\varepsilon)$ agree in first $l$ coefficients.

Note that by the previous remark, if $L(\varepsilon)$ is supposed to be analytic, then the power series solution $\hat{b}(\varepsilon)$ further gives rise to an analytic solution $b^*(\varepsilon) = \phi(\varepsilon) \cdot \bar{n}_{k+1}(\varepsilon)$.

Concerning the proof of (7.32), first note that $k + l$ leading terms of the approximation

$$b(\varepsilon) = b_0 + \cdots + \varepsilon^{k+l-1} \cdot b_{k+l-1} + \sum_{i=k+l}^{\infty} \varepsilon^i \cdot b_i \tag{7.33}$$

define a Jordan chain of length $k + l$ and we obtain from Lemma 1 (i)

$$\begin{pmatrix} b_{k+l-1} \\ \vdots \\ b_l \\ b_{l-1} \\ \vdots \\ b_0 \end{pmatrix} = \boxed{M}^{k+l} \cdot \begin{pmatrix} n_1 \\ \vdots \\ n_k \\ n_{1,k+1} \\ \vdots \\ n_{l,k+1} \end{pmatrix} \tag{7.34}$$

as well as

$$\begin{pmatrix} b_{l-1} \\ \vdots \\ b_0 \end{pmatrix} \overset{(4.10)}{=} \boxed{M}_l^{k+l} \cdot \begin{pmatrix} n_{1,k+1} \\ \vdots \\ n_{l,k+1} \end{pmatrix} \tag{7.35}$$

$$\overset{(5.13)}{=} \begin{pmatrix} \boxed{I_B} & \cdots & \boxed{\phi_{l-1}} \\ & \ddots & \vdots \\ & & \boxed{I_B} \end{pmatrix} \cdot \begin{pmatrix} n_{1,k+1} \\ \vdots \\ n_{l,k+1} \end{pmatrix}$$

with $n_{i,k+1} \in N_{k+1}, \ i = 1, \dots l$ due to stabilization at $k$. But then, the power series solution reads

$$\hat{b}(\varepsilon) = \phi(\varepsilon) \cdot \left( n_{l,k+1} + \cdots + \varepsilon^{l-1} \cdot n_{1,k+1} \right), \tag{7.36}$$

as can be seen along the following lines. First, $\hat{b}(\varepsilon)$ is a solution by $L(\varepsilon) \cdot \phi(\varepsilon) \cdot N_{k+1} = 0$ from Theorem 2 (i). Secondly, this solution agrees up to order $l-1$ with the approximation $b(\varepsilon)$ according to

$$\hat{b}(\varepsilon) = \left( I_B + \cdots + \varepsilon^{l-1} \cdot \phi_{l-1} + \sum_{i=l}^{\infty} \varepsilon^i \cdot \phi_i \right) \cdot \left( n_{l,k+1} + \cdots + \varepsilon^{l-1} \cdot n_{1,k+1} \right)$$

$$= (I_B) \cdot n_{l,k+1} + \cdots + \varepsilon^{l-1} \cdot (I_B \ \cdots \ \phi_{l-1}) \cdot \begin{pmatrix} n_{1,k+1} \\ \vdots \\ n_{l,k+1} \end{pmatrix} + \varepsilon^l \cdot r(\varepsilon)$$

$$= \left( \varepsilon^{l-1} \cdots 1 \right) \cdot \begin{pmatrix} \boxed{I_B} & \cdots & \boxed{\phi_{l-1}} \\ & \ddots & \vdots \\ & & \boxed{I_B} \end{pmatrix} \cdot \begin{pmatrix} n_{1,k+1} \\ \vdots \\ n_{l,k+1} \end{pmatrix} + \varepsilon^l \cdot r(\varepsilon) \tag{7.37}$$

$$\overset{(7.35)}{=} \left( \varepsilon^{l-1} \cdots 1 \right) \cdot \begin{pmatrix} b_{l-1} \\ \vdots \\ b_0 \end{pmatrix} + \varepsilon^l \cdot r(\varepsilon)$$

with a remainder power series $r(\varepsilon)$ and the proof of (7.32) is accomplished. Note also that the Greenberg function of the power series $L(\varepsilon)$ is given by

$$G(l) = k + l \quad for \quad l \geq 1 \ , \tag{7.38}$$

i.e. by (7.38) we transfer the standard result for the Greenberg function, concerning linear relations in finite dimensions [Rond], to linear equations in Banach spaces.

In addition, it is straightforward to ascertain the following Artin-Rees Lemma type result with respect to infinite dimensions (and noetherianity replaced by stabilization at $k$). For a power series $L(\varepsilon) = \sum_{i=0}^{\infty} \varepsilon^i \cdot L_i$, define the following sets of power series for $i \geq 0$

$$B_i := \{ \, b(\varepsilon) = \sum_{i=0}^{\infty} \varepsilon^i \cdot b_i \, , b_i \in B \ \mid \ L(\varepsilon) \cdot b(\varepsilon) = \varepsilon^i \cdot r(\varepsilon) \, \} \tag{7.39}$$

with $r(\varepsilon)$ a remainder power series, i.e. the set $B_i$ contains all power series approximations $b(\varepsilon)$ of $L(\varepsilon) \cdot b = 0$ of order $i$. Then, we obtain by direct calculation from (7.33)-(7.37) for $l \geq 1$ the following Artin-Rees type inclusion concerning sets of power series with coefficients in $\bar{B}$

$$\{ \, L(\varepsilon) \cdot B_{k+l} \, \} \subset \{ \, \varepsilon^l \cdot L(\varepsilon) \cdot B_0 \, \} \ . \tag{7.40}$$

Thus, the image of an approximation $b_{k+l}(\varepsilon)$ of order $k + l$ under $L(\varepsilon)$ can always be represented by the image of a power series from $B_0$ and scaling factor $\varepsilon^l$ split off, i.e.

$$L(\varepsilon) \cdot b_{k+l}(\varepsilon) = \varepsilon^l \cdot L(\varepsilon) \cdot b_0(\varepsilon) \quad \Leftrightarrow \quad L(\varepsilon) \cdot [\, b_{k+l}(\varepsilon) - \varepsilon^l \cdot b_0(\varepsilon) \,] = 0 \tag{7.41}$$

and the power series solution from (7.36) may also be written according to

$$\hat{b}(\varepsilon) = b_{k+l}(\varepsilon) - \varepsilon^l \cdot b_0(\varepsilon) \, . \tag{7.42}$$

In [Hauser], [Rond] the path over Artin-Rees Lemma is taken to show the existence of power series solutions for linear relations in finite dimensions.

Finally, note that the recursion (3.3)-(3.15) of section 3 can be performed, whenever direct sum decompositions according to (3.9) of the linear spaces $B$ and $\bar{B}$ are possible with respect to kernels and ranges. Hence, most of the investigations in this paper are not restricted to vector spaces and possibly it might be interesting to replace the field $\mathbb{K}$ by a ring $\mathbb{A}$.

*Matthias Stiefenhofer* *matthias.stiefenhofer@hs-kempten.de*
*Fachhochschule Kempten*
*Fakultät Maschinenbau*
*87435 Kempten (Germany)*